\documentclass[11pt, reqno, psamsfonts]{amsart}
\pdfoutput=1

\usepackage{amssymb}
\usepackage{amsthm}
\usepackage{amsmath}
\usepackage{latexsym}
\usepackage[T1]{fontenc}
\usepackage[utf8]{inputenc}
\usepackage[russian, french, english]{babel}

\usepackage{graphicx}
\usepackage{wrapfig}
\usepackage[justification=centering, labelfont=bf]{caption}
\usepackage{mathtools}
\usepackage{tikz}
\usepackage{amsbsy}
\usepackage[inline]{enumitem}
\usepackage{mathrsfs}
\usetikzlibrary{shapes,snakes}
\usetikzlibrary{arrows.meta}
\usetikzlibrary{decorations.pathmorphing}
\usetikzlibrary{patterns}
\usepackage{float}
\usepackage{array}
\usepackage{multicol}
\usepackage{stmaryrd}
\usepackage{eucal,eufrak}
\usepackage{mathabx}\usepackage{lmodern}
\usepackage{upgreek}
\usepackage{titlesec}
\usepackage[spacing=true,kerning=true,babel=true,tracking=true]{microtype}
\usepackage[shortcuts]{extdash}
\usepackage[foot]{amsaddr} 
\usepackage{framed}
\usepackage{algpseudocode}
\usepackage{algorithm}
\definecolor{shadecolor}{gray}{0.95}
\usepackage[left=1in,right=1in,top=1.05in,bottom=1.05in,bindingoffset=0cm]{geometry}
\usepackage{xspace}
\usepackage{bm}

\usepackage[backend=biber, style=alphabetic, sorting=nyt, maxnames=100, backref=true]{biblatex}
\usepackage[hidelinks]{hyperref}
\title{Distributed Algorithms, the Lov\'asz Local Lemma, \\and Descriptive Combinatorics}
\date{}
\author{Anton~Bernshteyn}
\address{\textls{\normalfont{}School of Mathematics, Georgia Institute of Technology, Atlanta, GA, USA}}
\email{bahtoh@gatech.edu}
\thanks{This research was partially supported by the NSF grant DMS-2045412.}

\newtheoremstyle{bfnote}%
{}{}%
{\slshape}{}%
{\bfseries}{\bfseries.}%
{ }%
{\thmname{#1}\thmnumber{ #2}\thmnote{ \ep{\normalfont{}#3}}}

\newtheoremstyle{defbfnote}%
{}{}%
{}{}%
{\bfseries}{.}%
{ }%
{\thmname{#1}\thmnumber{ #2}\thmnote{ (#3)}}

\newtheoremstyle{claim}%
{}{}%
{\slshape}{}%
{\itshape}{.}%
{ }%
{\thmname{#1}\thmnumber{ #2}\thmnote{ \ep{\normalfont{}#3}}}

\newtheoremstyle{smalldefn}%
{}{}%
{}{}%
{\itshape}{.}%
{ }%
{\thmname{#1}\thmnumber{ #2}\thmnote{ \ep{\normalfont{}#3}}}

\theoremstyle{bfnote}

\newtheorem{theo}[equation]{Theorem}
\newtheorem*{theo*}{Theorem}
\newtheorem{prop}[equation]{Proposition}
\newtheorem{lemma}[equation]{Lemma}

\newtheorem{conj}[equation]{Conjecture}

\newtheorem{lemdef}[equation]{Lemma/Definition}
\newtheorem*{claim*}{Claim}

\newtheorem*{corl*}{Corollary}

\theoremstyle{claim}

\newcounter{ForClaims}[section]
\newtheorem{claim}{Claim}[ForClaims]

\theoremstyle{smalldefn}
\newtheorem{smalldefn}[claim]{Definition}

\newcommand*{\myproofname}{Proof}
\newenvironment{claimproof}[1][\myproofname]{\begin{proof}[#1]}{\end{proof}}

\theoremstyle{definition}
\newtheorem{defn}[equation]{Definition}
\newtheorem*{defn*}{Definition}
\newtheorem{exmp}[equation]{Example}

\newtheorem*{exmp*}{Example}

\newtheorem{prob}[equation]{Open Problem}
\newtheorem*{assum*}{Assumptions}

\theoremstyle{remark}
\newtheorem*{ques*}{Question}
\newtheorem*{remk*}{Remark}

\newcommand{\0}{\varnothing}
\newcommand{\set}[1]{\{#1\}}

\newcommand{\dom}{\mathrm{dom}}

\newcommand{\acts}{\curvearrowright}
\newcommand{\N}{{\mathbb{N}}}
\newcommand{\Z}{\mathbb{Z}}

\newcommand{\R}{\mathbb{R}}

\renewcommand{\epsilon}{\varepsilon}
\renewcommand{\phi}{\varphi}
\renewcommand{\theta}{\vartheta}
\renewcommand{\leq}{\leqslant}
\renewcommand{\geq}{\geqslant}

\newcommand{\concat}{{^\smallfrown}}%
\newcommand{\fins}[1]{[#1]^{<\infty}}
\newcommand{\finseq}[1]{#1^{<\infty}}
\newcommand{\finf}[2]{[#1 \rightharpoonup #2]^{<\infty}}
\newcommand{\fun}[2]{[#1 \rightharpoonup #2]}
\newcommand{\G}{\Gamma}

\newcommand{\defeq}{\coloneqq}
\newcommand{\rest}[2]{{{#1}\vert{#2}}}
\newcommand{\emphd}[1]{{\fontseries{b}\selectfont\textsf{#1}}}

\renewcommand{\P}{\mathbb{P}}
\newcommand{\dist}{{\mathrm{dist}}}
\newcommand{\Diff}{\mathrm{d}}
\newcommand{\Nbhd}{\mathrm{N}}

\newcommand{\LOCAL}{$\mathsf{LOCAL}$\xspace}

\newcommand{\pr}{\mathsf{p}}
\newcommand{\de}{\mathsf{d}}
\newcommand{\Det}{\mathsf{Det}}
\newcommand{\Rand}{\mathsf{Rand}}
\newcommand{\class}[1]{\mathsf{#1}}

\newcommand{\bemph}[1]{{\normalfont#1}} %
\newcommand{\ep}[1]{\bemph{(}#1\bemph{)}} %

\newenvironment{scproof}[1][]{\begin{proof}[\textsc{\upshape{Proof}}#1]}{\end{proof}}

\newcommand{\pmp}{$\text{p.m.p.}$\xspace}

\numberwithin{equation}{section}

\renewcommand{\thesubsection}{\arabic{section}.\Alph{subsection}}

\usepackage{etoolbox}

\titleformat{\section}[block]{\scshape\filcenter}{\thesection.}{1ex}{}
\titleformat{\subsection}[block]{\bfseries\filcenter}{\thesubsection.}{1ex}{}
\titleformat{\subsubsection}[runin]{\itshape}{\bfseries\upshape\thesubsubsection.}{1ex}{}[.---]

\titlespacing*{\section}{0pt}{*3}{*1}
\titlespacing*{\subsection}{0pt}{*3}{*1}
\titlespacing*{\subsubsection}{0pt}{*1.5}{*0}

\makeatletter
\newcommand{\neutralize}[1]{\expandafter\let\csname c@#1\endcsname\count@}
\makeatother

\newenvironment{theobis}[1]
{\renewcommand{\theequation}{\ref{#1}$'$}%
	\neutralize{equation}\phantomsection
	\begin{theo}}
	{\end{theo}}

\renewbibmacro{in:}{}

\renewbibmacro*{volume+number+eid}{%
	\printfield{volume}%
	\setunit*{\addnbspace}%
	\printfield{number}%
	\setunit{\addcomma\space}%
	\printfield{eid}}

\DeclareFieldFormat[article]{volume}{\textbf{#1}\space}
\DeclareFieldFormat[article]{number}{\mkbibparens{#1}}

\DeclareFieldFormat{journaltitle}{#1,}
\DeclareFieldFormat[thesis]{title}{\mkbibemph{#1}\addperiod}
\DeclareFieldFormat[article, unpublished, thesis]{title}{\mkbibemph{#1},}
\DeclareFieldFormat[book]{title}{\mkbibemph{#1}\addperiod}
\DeclareFieldFormat[unpublished]{howpublished}{#1, }

\DeclareFieldFormat{pages}{#1}

\DeclareFieldFormat[article]{series}{Ser.~#1\addcomma}

\setlist{topsep=3pt,itemsep=3pt}

\begin{document}
	
	\pagestyle{plain}

	\maketitle
	
	\begin{abstract}
		In this paper we consider coloring problems on graphs and other combinatorial structures on standard Borel spaces. Our goal is to obtain sufficient conditions under which such colorings can be made well-behaved in the sense of topology or measure. To this end, we show that such well-behaved colorings can be produced using certain powerful techniques from finite combinatorics and computer science. First, we prove that efficient distributed coloring algorithms (on finite graphs) yield well-behaved colorings of Borel graphs of bounded degree; roughly speaking, deterministic algorithms produce Borel colorings, while randomized algorithms give measurable and Baire-measurable colorings. Second, we establish measurable and Baire-measurable versions of the Symmetric Lov\'asz Local Lemma (under the assumption $\pr(\de+1)^8 \leq 2^{-15}$, which is stronger than the standard LLL assumption $\pr(\de + 1) \leq e^{-1}$ but still sufficient for many applications). From these general results, we derive a number of consequences in descriptive combinatorics and ergodic theory.
	\end{abstract}
	
	\section{Introduction}\label{sec:intro}
	
	A coloring, broadly construed, is a mapping that assigns to each element of a given structure one of a \ep{typically finite} number of ``colors'' in a way that fulfills a prescribed set of constraints. The prototypical example is a proper $k$-coloring of a graph $G$, i.e., a mapping $f \colon V(G) \to [k]$ such that $f(x) \neq f(y)$ whenever the vertices $x$ and $y$ are adjacent in $G$. Classical problems in combinatorics often seek to identify sufficient conditions under which a coloring of a certain type exists.
	
	In the last twenty or so years, a rich theory has emerged concerning the behavior of colorings and other combinatorial notions under additional regularity requirements. For instance, suppose that $G$ is a graph whose vertex set $V(G)$ is a standard probability space \ep{e.g., the unit interval $[0,1]$ with the usual Lebesgue measure}. Does $G$ admit a proper $k$-coloring $f \colon V(G) \to [k]$ that is measurable, meaning that $f^{-1}(c)$ is a measurable subset of $V(G)$ for every color $c$? Questions of this type are the subject matter of \emph{descriptive combinatorics}, which facilitates the fusion of ideas from combinatorics and descriptive set theory. For a state-of-the-art introduction to this area, see the survey \cite{KechrisMarks} by Kechris and Marks.
	
	Throughout this paper, we shall work under the assumption that the set of all available colors is countable \ep{although there are interesting questions about uncountable colorings as well, see \cite[\S\S3 and 6]{KechrisMarks}}. For concreteness, we can then assume that the colors form a subset of $\N = \set{0,1,2, \ldots}$, so a coloring of a structure $X$ induces a partition $X = X_0 \sqcup X_1 \sqcup X_2 \sqcup \ldots$ of $X$ into countably many labeled pieces, called the \emph{color classes}. Assuming that $X$ is a standard Borel space, we can then classify such colorings according to the regularity properties of their color classes. For instance, one can study \emph{Borel} colorings, i.e., colorings in which every color class $X_i$ is a Borel subset of $X$. This notion turns out to be rather restrictive, so one often considers colorings in which every color class is not necessarily Borel but \emph{measurable} with respect to some probability Borel measure $\mu$ on $X$ or \emph{Baire-measurable} with respect to some compatible Polish topology $\tau$ on $X$. On the other hand, sometimes one can aim for something even stronger than Borel; for instance, one might consider \emph{continuous} colorings, in which every color class is a clopen set.
	
	The goal of this paper is to adapt certain powerful techniques from finite combinatorics to the descriptive setting. Specifically, we establish the following two groups of results:
	
	\begin{enumerate}[label=\ep{\normalfont\Alph*}]
		\item\label{item:A} If a coloring problem on finite graphs can be solved by an efficient distributed algorithm, then on infinite graphs the same problem admits solutions with some regularity properties. \ep{See Theorems~\ref{theo:dist_Borel}, \ref{theo:dist_cont}, \ref{theo:dist_meas}, and \ref{theo:dist_subexp} in \S\ref{subsec:dist_res} for the precise statements.}
		
		\item\label{item:B} There is a measurable/Baire-measurable version of the Symmetric Lov\'asz Local Lemma. \ep{See Theorem~\ref{theo:meas_SLLL} in \S\ref{subsec:LLL} for the precise statement.}
	\end{enumerate}
	
	From these general facts, we derive a variety of new results in descriptive combinatorics, presented in \S\ref{sec:applications}. Here are just two examples \ep{see \S\ref{sec:applications} for the definitions of the terms used}:
	
	\begin{theo}[see Theorem~\ref{theo:Delta}]\label{theo:1}
		There is $\Delta_0 \in \N$ with the following property. Fix integers $\Delta \geq \Delta_0$ and $c$. Let $G$ be a Borel graph of maximum degree at most $\Delta$ and let $\mu$ be a probability Borel measure  on $V(G)$. If %
		$c \leq \sqrt{\Delta + 1/4} - 5/2$, then the following statements are equivalent: %
		\begin{enumerate}[label=\ep{\normalfont\roman*}]
			\item the chromatic number of $G$ is at most $\Delta - c$;
			
			\item the $\mu$-measurable chromatic number of $G$ is at most $\Delta - c$.
		\end{enumerate}
	\end{theo}
	
	\begin{theo}[see Theorem~\ref{theo:triangle_free}]\label{theo:2}
		For every $\epsilon > 0$, there is $\Delta_0 \in \N$ with the following property. Let $G$ be a Borel graph of finite maximum degree $\Delta \geq \Delta_0$ and let $\mu$ be a probability Borel measure on $V(G)$. If $G$ contains no cycles of length $3$ \ep{resp.~at most $4$}, then the $\mu$-measurable chromatic number of $G$ is at most $(4+\epsilon)\Delta/\log \Delta$ \ep{resp.~$(1+\epsilon)\Delta/\log \Delta$}. %
	\end{theo}
	
	Theorem~\ref{theo:1} was previously known only for $c \leq 0$: the $c < 0$ case is due to Kechris, Solecki, and Todorcevic \cite[Proposition 4.6]{KechrisSoleckiTodorcevic} and the $c=0$ case is due to Conley, Marks, and Tucker-Drob \cite{CMTD}. The upper bound on $c$ in the statement of Theorem~\ref{theo:1} is within $2$ of best possible \ep{see Proposition~\ref{prop:sharp}}. The only previously known general upper bound on the measurable chromatic number of Borel graphs without $3$-cycles, or even with no cycles at all, in terms of their maximum degree was $\Delta$, which is implied by the result of Conley, Marks, and Tucker-Drob \cite{CMTD}. Lyons and Nazarov \cite{LyonsNazarov} observed \ep{see also \cite[Theorem 5.46]{KechrisMarks}}
	that there are acyclic Borel graphs with maximum degree $\Delta$ and measurable chromatic number at least $(1/2+o(1))\Delta/\log \Delta$, which means that the bounds in Theorem~\ref{theo:2} are optimal up to a constant factor, even for acyclic graphs.
	
	Let us now give some more details about bullet points \ref{item:A} and \ref{item:B}. Distributed computing is an area of computer science that, among other things, investigates questions of the following form: Given a graph coloring problem, how far in the graph should an individual vertex be allowed to ``see'' in order to be able to compute its own color? This idea is formalized in the \LOCAL model of distributed computation introduced by Linial \cite{Linial}. We describe this model here in a somewhat informal way; a precise definition, in the form needed for our purposes, is given in \S\ref{subsec:dist}.
	
	In the \LOCAL model an $n$-vertex graph $G$ abstracts a communication network where each vertex plays the role of a processor and edges represent communication links. The algorithm proceeds in \emph{rounds}. During each round, the vertices first perform arbitrary local computations and then synchronously broadcast messages to all their neighbors. At the end, each vertex should output its own part of the global solution \ep{i.e., its own color}. There are no restrictions on the complexity of the local computations involved or on the length of the messages that the vertices send to each other, and the only measure of efficiency for such an algorithm is the number of communication rounds required.
	
	We emphasize that in the \LOCAL model, every vertex is executing the same algorithm. This means that, to make this model nontrivial, the vertices must be given a way of breaking symmetry. There are two standard symmetry-breaking approaches, leading to the distinction between deterministic and randomized \LOCAL algorithms:
	
	\begin{itemize}
		\item In the \emph{deterministic} \LOCAL model, each vertex is assigned, as part of its input, a unique $\Theta(\log n)$-bit identifier. The algorithm executed at each vertex is deterministic and must always output a correct solution to the problem, regardless of the specific assignment of the identifiers.
		
		\item In the \emph{randomized} \LOCAL model, each vertex may generate an arbitrarily long finite sequence of independent uniformly distributed random bits. The algorithm may fail to produce a correct solution to the problem, but the probability of failure must not exceed $1/n$.
	\end{itemize}
	We remark that deterministic \LOCAL algorithms can be simulated in the randomized \LOCAL model: each vertex can simply generate a random sequence of $\lceil 3\log_2 n \rceil$ bits and use it as an identifier---the probability that two identifiers generated in this way coincide is less than $1/n$.
	
	Notice that in a \LOCAL algorithm that runs on a graph $G$ for $T$ rounds, each vertex only has access to information in its radius-$T$ neighborhood. Conversely, every $T$-round \LOCAL algorithm can be transformed into one in which every vertex first collects all the information contained in its radius-$T$ neighborhood and then makes a decision, based on this information alone, about its color. This alternative way of thinking about \LOCAL algorithms makes it clear how they measure the ``locality'' of graph coloring problems and is often more convenient to work with. \ep{Indeed, this is the approach we shall use in \S\ref{subsec:dist} to define \LOCAL algorithms formally.} For further background on distributed coloring algorithms, see the book \cite{BE} by Barenboim and Elkin.
	
	Now we can describe some of	our results in group \ref{item:A}; for their precise statements, see \S\ref{subsec:dist_res}. Fix $\Delta \in \N$. We show that if a coloring problem can be solved by a deterministic \LOCAL algorithm that runs in $o(\log n)$ rounds on $n$-vertex graphs of maximum degree $\Delta$, then the same problem for Borel graphs of maximum degree $\Delta$ admits Borel solutions \ep{see Theorem~\ref{theo:dist_Borel}}. Furthermore, under extra topological assumptions, ``Borel'' here may be replaced by ``continuous'' \ep{see Theorem~\ref{theo:dist_cont}}. Similarly, randomized \LOCAL algorithms that run in $o(\log n)$ rounds yield both measurable and Baire-measurable colorings \ep{see Theorem~\ref{theo:dist_meas}}, and ``measurable'' may be upgraded to ``Borel'' under additional assumptions on the growth rate of the underlying graph \ep{see Theorem~\ref{theo:dist_subexp}}. These general facts enable one to prove new results in descriptive combinatorics simply by alluding to known distributed algorithms; for several such examples, see \S\ref{sec:applications}.
	
	While the proofs of Theorems~\ref{theo:dist_Borel} and \ref{theo:dist_cont} \ep{dealing with deterministic algorithms} are relatively straightforward, Theorem~\ref{theo:dist_meas} \ep{concerning randomized algorithms} is more involved. In particular, it relies on a novel measurable version of the Lov\'asz Local Lemma that we establish in this paper, which brings us to bullet point \ref{item:B}.
	
	The \emph{Lov\'asz Local Lemma} \ep{the \emph{LLL} for short} is a powerful tool in probabilistic combinatorics, introduced by Erd\H{o}s and Lov\'asz in the mid-1970s \cite{EL}. The LLL is mostly used to obtain existence results, and it is particularly well-suited for finding colorings that satisfy some ``local'' constraints. Roughly speaking, in order for the LLL to apply in this context, two requirements must be met: First, a random coloring must be ``likely'' to fulfill each individual constraint; second, the constraints must not interact with each other ``too much.'' The precise statement of the LLL requires a few technical definitions, so we postpone it until \S\ref{subsubsec:LLL}.
	
	It has been a matter of interest to determine if the LLL can yield ``constructive'' conclusions \ep{rather than pure existence results}. The first such ``constructive'' version of the LLL was the \emph{algorithmic LLL} due to Beck \cite{Beck}. Beck's result requires somewhat stronger numerical assumptions than the ordinary LLL. This discrepancy has been eventually eliminated in the breakthrough work of Moser and Tardos \cite{MT}; for the long list of intermediate results, see the references in \cite{MT}.
	
	The Moser--Tardos method was later adapted to derive ``constructive'' analogs of the LLL in a variety of different contexts. For example, Rumyantsev and Shen \cite{RSh} proved a \emph{computable} version of the LLL. Here we are interested in the behavior of the LLL in the \emph{descriptive} setting. When can the LLL be used to obtain Borel, measurable, or Baire-measurable colorings? Partial answers to this question have been given in \cite{MLLL} by the present author and in \cite{CGMPT} by Cs\'{o}ka, Grabowski, M\'ath\'e, Pikhurko, and Tyros. %
	The main results of both \cite{MLLL} and \cite{CGMPT} only apply under rather special circumstances: \cite[Theorem~6.6]{MLLL} is a measurable version of the LLL for structures induced by the Bernoulli shift actions $\G \acts [0,1]^\G$ of countable groups $\G$, while \cite[Theorem 1.3]{CGMPT} is a Borel LLL for coloring problems on graphs of subexponential growth.
	
	In this paper we establish the first measurable/Baire-measurable version of the LLL that works without any such special restrictions \ep{see Theorem~\ref{theo:meas_SLLL}}. The drawback is that, as in the seminal algorithmic LLL of Beck, we require stronger numerical bounds than in the ordinary LLL. Nevertheless, just like Beck's result, \hyperref[theo:meas_SLLL]{our measurable LLL} is sufficient for many combinatorial applications. In particular, as mentioned above, it plays the central role in our proof that \hyperref[theo:dist_meas]{sublogarithmic randomized \LOCAL algorithms yield measurable/Baire-measurable colorings}, and hence it lies at the heart of the measurable coloring results presented in \S\ref{sec:applications}. Of course, \hyperref[theo:meas_SLLL]{our measurable LLL} can also be used without any reference to distributed algorithms. We give an example of such a direct application to a problem in ergodic theory in \S\ref{subsec:AW}.
	
	We remark that, while we use \hyperref[theo:meas_SLLL]{the measurable LLL} to establish a relationship between randomized \LOCAL algorithms and measurable colorings, our proof of \hyperref[theo:meas_SLLL]{the measurable LLL} itself relies on the randomized \LOCAL algorithm for the LLL that was recently developed by Fischer and Ghaffari \cite{FG} and sharpened by Ghaffari, Harris, and Kuhn \cite{GHK}. Curiously, both the Fischer--Ghaffari algorithm and our proof of the measurable LLL do not invoke the Moser--Tardos method and instead go back essentially to the original ideas of Beck.
	
	The remainder of this paper is organized as follows. We give the necessary definitions and state our main results precisely in \S\ref{sec:statements}. Next we present a variety of applications in \S\ref{sec:applications}. In \S\ref{sec:dist_to_desc}, we explain how to turn distributed algorithms into colorings with regularity properties. Specifically, \S\ref{sec:dist_to_desc} contains the proofs of our results concerning deterministic \LOCAL algorithms and reduces the results about randomized algorithms to the LLL. Finally, we prove our measurable LLL in \S\ref{sec:LLL_proof}.
	
	\subsubsection*{{Acknowledgments}}
	
	I am very grateful to Clinton Conley for many insightful conversations and to the anonymous referees for reading the manuscript carefully and providing helpful comments.
	
	\section{Main definitions and results}\label{sec:statements}
	
	\subsection{Distributed algorithms for graph coloring problems}\label{subsec:dist}
	
	In this section we formally introduce the terminology pertaining to the \LOCAL model of distributed computation. The nature of our main results requires us to be somewhat more pedantic with our definitions than is standard in the distributed algorithms literature. The reader who is already familiar with the \LOCAL model is encouraged to only skim this section in order to familiarize herself with our notation and quickly move on to \S\ref{subsec:dist_res}, where we connect \LOCAL algorithms to Borel/measurable colorings.
	
	\subsubsection{Graphs and structured graphs}
	
	Given a set $A$, we write $\finseq{A}$ \ep{resp.~$\fins{A}$} to denote the set of all finite sequences \ep{resp.~finite sets} of elements of $A$. Unless otherwise specified, by a ``graph'' we mean a simple undirected graph. Our graph-theoretic terminology and notation are standard; see, e.g., Diestel's book \cite{Die}. In particular, the vertex and edge sets of a graph $G$ are denoted by $V(G)$ and $E(G)$ respectively. As usual, we write $|G| \defeq |V(G)|$ and $\|G\| \defeq |E(G)|$. For $U \subseteq V(G)$, the \emphd{neighborhood} of $U$ in $G$ \ep{i.e., the set of all vertices with a neighbor in $U$} is denoted by $N_G(U)$. For a vertex $x \in V(G)$, we set $N_G(x) \defeq N_G(\set{x})$ and use $\deg_G(x) \defeq |N_G(x)|$ to denote the \emphd{degree} of $x$ in $G$. The \emphd{maximum degree} of $G$ is defined by $\Delta(G) \defeq \sup_{x \in V(G)} \deg_G(x)$. A graph $G$ is said to be \emphd{locally finite} if $\deg_G(x) < \infty$ for all $x \in V(G)$. Throughout this paper we only work with locally finite graphs; in fact, we mostly focus on graph whose maximum degree is finite.
	
	Sometimes one needs to consider graphs that support some additional structure. For instance, one might wish to work with directed graphs, weighted graphs, or graphs with a fixed coloring of the vertices. We capture this idea in the general notion of a ``structured graph'':
	
	\begin{defn}[\textls{Structured graphs}]\label{defn:str_graph}
		A \emphd{structure map} on a graph $G$ is a partial function $\sigma \colon \finseq{V(G)} \rightharpoonup \N$ such that for some $\ell \in \N$, every tuple $\bm{x} \in \dom(\sigma)$ is of length at most $\ell$. A \emphd{structured graph} is a pair $\bm{G} = (G, \sigma)$, where $G$ is a graph and $\sigma$ is a structure map on $G$.
	\end{defn}
	
	In the above definition, we only use $\N$ as the range of $\sigma$ for convenience, as having the range of $\sigma$ fixed once and for all will simplify some of the forthcoming definitions. In applications, any countable set could be used instead of $\N$, and indeed we shall often abuse terminology by referring to structured graphs whose structure maps range over different countable sets. For instance, we can view a graph $G$ equipped with a finite sequence $\sigma_1$, \ldots, $\sigma_k$ of structure maps as a structured graph by replacing the tuple $(\sigma_1, \ldots, \sigma_k)$ with the single function
	\[
	\sigma \colon \bigcup_{i = 1}^k \dom(\sigma_i) \to (\N \cup \set{\ast})^k \colon \bm{x} \mapsto (\sigma_1(\bm{x}), \ldots, \sigma_k(\bm{x})),
	\]
	where $\ast$ is a special symbol distinct from all the elements of $\N$ and $\sigma_i(\bm{x}) = \ast$ means that $\sigma_i(\bm{x})$ is undefined. This is an acceptable construction since the set $(\N \cup \set{\ast})^k$ is countable. As a special case, given a structured graph $\bm{G} = (G, \sigma)$ and another structure map $\tau$ on $G$, we can interpret the pair $(\bm{G}, \tau)$ as a new structured graph \ep{with some ``extra structure'' added to that of $\bm{G}$}. This convention will become important when we formally define \LOCAL algorithms.
	
	\begin{exmp}[\textls{Directed graphs}]
		A directed graph $G$ can naturally be viewed as a structured graph, since we can encode the directions of the edges of $G$ using the function $V(G)^2 \to \set{0,1}$ that sends a pair $(x, y)$ to $1$ if and only if there is a directed edge from $x$ to $y$ in $G$.
	\end{exmp}
	
	\begin{exmp}[\textls{Multigraphs}]
		If the edges of a graph $G$ are allowed to have \ep{finite} multiplicities, then we can view $G$ as a structured graph equipped with the function $V(G)^2 \to \N$ that assigns to each pair $(x,y)$ the multiplicity of the edge $xy$.
	\end{exmp}
	
	\begin{exmp}[\textls{List-coloring}]\label{exmp:list-col}
		Sometimes, instead of assigning to each vertex of a graph $G$ a color from a fixed set, such as $[k]$, one is required to pick a color for every $x \in V(G)$ from its own list $L(x)$ of available colors \ep{this is called the \emph{list-coloring problem}; see \cite[\S5.4]{Die}}. If every list $L(x)$ is a finite subset of $\N$ \ep{or of any other countable set}, then the pair $(G, L)$ can be naturally viewed as a structured graph \ep{since the set $\fins{\N}$ is countable}.
	\end{exmp}
	
	We naturally extend all the standard graph-theoretic notation, such as $V$, $\Delta$, etc., to structured graphs; that is, for a structured graph $\bm{G} = (G, \sigma)$, we write $V(\bm{G}) \defeq V(G)$, $\Delta(\bm{G}) \defeq \Delta(G)$, etc.
	
	Given a graph $G$ and a subset $U \subseteq V(G)$, $G[U]$ denotes the \emphd{subgraph of $G$ induced by $U$}, i.e., the graph with vertex set $U$ in which two vertices $x$, $y \in U$ are adjacent if and only if they are adjacent in $G$. Similarly, if $\bm{G} = (G, \sigma)$ is a structured graph and $U \subseteq V(G)$, then $\bm{G}[U]$ is the structured graph given by $\bm{G}[U] \defeq (G[U], \rest{\sigma}{\finseq{U}})$. \ep{Here and in what follows we use the notation $\rest{f}{X}$ to indicate the restriction of a function $f$ onto the set $X \cap \dom(f)$.}
	
	As usual, the \emphd{distance} $\dist_G(x,y)$ between two vertices $x$, $y \in V(G)$ is the smallest number of edges on a path in $G$ that starts at $x$ and ends at $y$ \ep{if there is no such path, then $\dist_G(x,y)\defeq \infty$}. For a vertex $x \in V(G)$ and a number $R \in \N$, we define $B_G(x, R)$ to be the \emphd{ball} of radius $R$ around $x$ in $G$, i.e., the subgraph of $G$ induced by the set $\set{y \in V(G) \,:\, \dist_G(x,y) \leq R}$. The definition of $B_{\bm{G}}(x, R)$ for structured graphs $\bm{G}$ is the same, \emph{mutatis mutandis}.
	
	We say that two structured graphs $\bm{G}_1 = (G_1, \sigma_1)$ and $\bm{G}_2 = (G_2, \sigma_2)$ are isomorphic if there is an isomorphism between their underlying graphs $G_1$ and $G_2$ that turns the function $\sigma_1$ into $\sigma_2$. More formally, each function $\phi \colon V(G_1) \to V(G_2)$ can be extended to a map $\finseq{V(G_1)} \to \finseq{V(G_2)}$ in the obvious way; namely, given a tuple $\bm{x} = (x_1, \ldots, x_k) \in \finseq{V(G_1)}$, define
	\[
	\phi(\bm{x}) \,\defeq\, (\phi(x_1), \ldots, \phi(x_k)) \,\in\, \finseq{V(G_2)}.
	\]
	An \emphd{isomorphism} between $\bm{G}_1$ and $\bm{G}_2$ is a function $\phi \colon V(G_1) \to V(G_2)$ such that:
	\begin{itemize}
		\item $\phi$ is an isomorphism of the graphs $G_1$ and $G_2$;
		\item $\dom(\sigma_2) = \phi(\dom(\sigma_1))$; and
		\item for all $\bm{x} \in \dom(\sigma_1)$, $\sigma_1(\bm{x}) = \sigma_2(\phi(\bm{x}))$.
	\end{itemize}
	We say that $\bm{G}_1$ and $\bm{G}_2$ are \emphd{isomorphic}, in symbols $\bm{G}_1 \cong \bm{G}_2$, if there is an isomorphism between $\bm{G}_1$ and $\bm{G}_2$. The isomorphism class of a structured graph $\bm{G}$ is denoted by $[\bm{G}]$. The set of all the isomorphism classes of finite structured graphs is denoted by $\class{FSG}$. Notice that $\class{FSG}$ is a countable set \ep{here we exploit the fact that the range of $\sigma$ in Definition~\ref{defn:str_graph} is assumed to be the fixed set $\N$}.
	
	A \emphd{rooted graph} is a pair $(G, x)$, where $G$ is a graph and $x \in V(G)$ is a distinguished vertex, called the \emphd{root}. Isomorphisms between rooted graphs are required to preserve the roots. One similarly defines \emphd{rooted structured graphs} $(\bm{G}, x)$. To simplify the notation, we denote the isomorphism class of a rooted structured graph $(\bm{G}, x)$ by $[\bm{G}, x]$ \ep{instead of $[(\bm{G}, x)]$}. The set of all the isomorphism classes of finite rooted structured graphs is denoted by $\class{FSG}_\bullet$. Again, the set $\class{FSG}_\bullet$ is countable. 
	
	\subsubsection{\LOCAL algorithms}
	
	Now we have enough notation to start defining the \LOCAL model.
	
	\begin{defn}[\textls{\LOCAL algorithms}]\label{defn:LOCAL}
		A \emphd{\LOCAL algorithm} is a function $\mathcal{A} \colon \class{FSG}_\bullet \to \N$. Given a locally finite structured graph $\bm{G}$ and a natural number $T \in \N$, the \emphd{output of $\mathcal{A}$ on $\bm{G}$ after $T$ rounds} is the function $\mathcal{A}(\bm{G}, T) \colon V(\bm{G}) \to \N$ given by
		\[
		\mathcal{A}(\bm{G}, T)(x) \,\defeq\, \mathcal{A}\left([B_{\bm{G}}(x, T), x]\right) \qquad \text{for all } x \in V(\bm{G}).
		\]
		Note that since $\bm{G}$ is locally finite, $B_{\bm{G}}(x,T)$ is finite, so this definition makes sense.
	\end{defn}
	
	Informally, Definition~\ref{defn:LOCAL} says that a \LOCAL algorithm $\mathcal{A}$ operates on a locally finite structured graph $\bm{G}$ as follows. Each vertex $x \in V(\bm{G})$ ``sees'' the isomorphism type of its radius-$T$ ball $B_{\bm{G}}(x, T)$, which is viewed as a structured graph rooted at $x$. The algorithm $\mathcal{A}$ is then a ``rule'' that uses this information to output a ``color'' $\mathcal{A}(\bm{G}, T)(x) \in \N$. The same comment as after Definition~\ref{defn:str_graph} applies here: instead of the set $\N$ in Definition~\ref{defn:LOCAL}, any other countable set of possible outputs could be \ep{and often is} used.
	
	So far we have defined how \LOCAL algorithms operate. Now we need to describe what problems they can solve. The problems we consider are sometimes called ``locally checkable labeling problems'' in the distributed computing literature \ep{this term was introduced in the seminal work of Naor and Stockmeyer \cite{NaorStock}}. However, in order to stay closer to the terminology in descriptive combinatorics and in graph theory, we prefer to use the word ``coloring'' instead of ``labeling.''
	
	Recall that if $\bm{G}$ is a structured graph and $f \colon V(\bm{G}) \rightharpoonup \N$ is a partial function \ep{or, more generally, a structure map}, then the pair $(\bm{G}, f)$ can be thought of as a structured graph in its own right, in which the function $f$ is ``added'' to the structure \ep{see the comments after Definition~\ref{defn:str_graph}}. For brevity, we denote this structured graph by $\bm{G}_f$.
	
	\begin{defn}[\textls{Local coloring problems}]\label{defn:LCL}
		A \emphd{local coloring problem} \ep{or a \emphd{locally checkable labeling problem}} is a pair $\Pi = (t, \mathcal{P})$, where $t \in \N$ and $\mathcal{P} \colon \class{FSG}_\bullet \to \set{0,1}$. Here we interpret $\mathcal{P}$ as a \LOCAL algorithm. Let $\bm{G}$ be a locally finite graph and suppose that we are given a function $f \colon V(\bm{G}) \to \N$. We say that $f$ is a \emphd{$\Pi$-coloring} of $\bm{G}$ if $\mathcal{P}(\bm{G}_f, t)(x) = 1$ for all $x \in V(\bm{G})$, i.e., if the output of the algorithm $\mathcal{P}$ on the structured graph $\bm{G}_f$ after $t$ rounds is the constant $1$ function.
	\end{defn}
	
	What Definition~\ref{defn:LCL} says is that, for a local coloring problem, there is a \LOCAL algorithm, $\mathcal{P}$, that, given a coloring $f$, can verify whether $f$ is ``correct'' in a constant number of rounds, namely $t$. To put this another way, whether or not a coloring $f$ is ``correct'' is completely determined by the restrictions of $f$ to balls of radius $t$.
	
	\begin{exmp}[\textls{Proper coloring}]\label{exmp:prop_col}
		The prototypical example of a local coloring problem is proper $k$-coloring of graphs, since whether or not a coloring is proper is determined by its restrictions to radius-$1$ balls. Explicitly, let $k \in \N$ and define a \LOCAL algorithm $\mathcal{P} \colon \class{FSG}_\bullet \to \set{0,1}$ as follows: Given \ep{the isomorphism type of} a finite rooted structured graph of the form $(G, f, x)$ with $f \colon V(G) \to [k]$ and a root $x$, set $\mathcal{P}([G, f, x]) \defeq 1$ if and only if $f(y) \neq f(z)$ for every pair of adjacent vertices $y$, $z \in V(G)$; in all other cases set $\mathcal{P}$ to $0$. Now if we let $\Pi \defeq (1, \mathcal{P})$, then a $\Pi$-coloring of a locally finite graph $G$ is exactly the same as a proper $k$-coloring of $G$.
	\end{exmp}
	
	\subsubsection{Complexity of local coloring problems}
	
	As mentioned in the \hyperref[sec:intro]{introduction}, we must distinguish between the deterministic and the randomized versions of the \LOCAL model.
	
	\begin{defn}[\textls{Deterministic \LOCAL complexity}]\label{defn:det_comp}
		Let $\Pi$ be a local coloring problem and let $\class{G} \subseteq \class{FSG}$. Given $n \in \N^+$ and $T \in \N$, we write $\Det_{\Pi, \class{G}}(n) \leq T$ if and only if there is a \LOCAL algorithm $\mathcal{A}$ with the following property:

		\begin{leftbar}
			\noindent Let $\bm{G}$ be an $n$-vertex structured graph such that $[\bm{G}] \in \class{G}$ and let $\mathsf{id} \colon V(\bm{G}) \to [n]$ be an arbitrary bijection. Then the function $\mathcal{A}(\bm{G}_\mathsf{id}, T)$ is a $\Pi$-coloring of $\bm{G}$.
		\end{leftbar}

		\noindent For $n \in \N^+$, define $\Det_{\Pi, \class{G}}(n)$ to be the least $T \in \N$ such that $\Det_{\Pi, \class{G}}(n) \leq T$ if such $T$ exists, and $\infty$ otherwise. The function $\Det_{\Pi, \class{G}} \colon \N^+ \to \N \cup \set{\infty}$ is called the \emphd{deterministic \LOCAL complexity} of the problem $\Pi$ on the class of structured graphs $\class{G}$.
	\end{defn}
	
	Let us make a few remarks about the above definition. The set $\class{G}$ in Definition~\ref{defn:det_comp} is the collection of \ep{the isomorphism classes of} the finite structured graphs on which we attempt to solve the given coloring problem $\Pi$. For example, $\class{G}$ may contain all graphs of maximum degree at most $\Delta$, all triangle-free graphs, all trees, all cycles, all directed graphs, etc. If $\Det_{\Pi, \class{G}}(n) \leq T$, it means that there is a \LOCAL algorithm $\mathcal{A}$ that solves the problem $\Pi$ on $n$-vertex graphs $\bm{G}$ from the class $\class{G}$ in $T$ rounds, given as part of its input an arbitrary assignment $\mathsf{id}$ of unique identifiers from the set $[n]$ to the vertices of $\bm{G}$. The reader may recall from the \hyperref[sec:intro]{introduction} that usually the identifiers are sequences of bits of length $\Theta(\log n)$, or, equivalently, elements of $[n^c]$ for some constant $c \geq 1$. In practice it does not matter what constant $c$ one uses, so, for concreteness, we fix $c$ to be $1$ \ep{technically, making $c$ as small as possible only makes our main results stronger}.

	\begin{defn}[\textls{Randomized \LOCAL complexity}]\label{defn:rand_comp}
		Let $\Pi$ be a local coloring problem and let $\class{G} \subseteq \class{FSG}$. Given $n \in \N^+$ and $T \in \N$, we write $\Rand_{\Pi, \class{G}}(n) \leq T$ if and only if there are $m \in \N^+$ and a \LOCAL algorithm $\mathcal{A}$ with the following property:

		\begin{leftbar}
			\noindent Let $\bm{G}$ be an $n$-vertex structured graph such that $[\bm{G}] \in \class{G}$. Pick a function $\theta \colon V(\bm{G}) \to [m]$ uniformly at random \ep{that is, each function $V(\bm{G}) \to [m]$ is chosen with probability $1/m^n$}. Then
			\[
			\P \left[\text{$\mathcal{A}(\bm{G}_\theta, T)$ is a $\Pi$-coloring of $\bm{G}$}\right] \,\geq\, 1 \,-\, \frac{1}{n}.
			\]
		\end{leftbar}

		\noindent For $n \in \N^+$, define $\Rand_{\Pi, \class{G}}(n)$ to be the least $T \in \N$ such that $\Rand_{\Pi, \class{G}}(n) \leq T$ if such $T$ exists, and $\infty$ otherwise. The function $\Rand_{\Pi, \class{G}} \colon \N^+ \to \N \cup \set{\infty}$ is called the \emphd{randomized \LOCAL complexity} of the problem $\Pi$ on the class of structured graphs $\class{G}$.
	\end{defn}
	
	In Definition~\ref{defn:rand_comp}, $\class{G}$ is again the class of finite structured graphs on which we wish to solve a given coloring problem. The algorithm $\mathcal{A}$ takes, as part of its input, an assignment of random numbers from $[m]$ to the vertices of $\bm{G}$, and, while it may fail to output a $\Pi$-coloring, the probability of failure cannot exceed $1/n$. Notice that there is no a priori upper bound on $m$, i.e., the algorithm is allowed to use ``unlimited randomness.'' However, we do require that the same $m$ must work for all $n$-vertex graphs $\bm{G}$ with $[\bm{G}] \in \class{G}$. This restriction is harmless, at least for our purposes, as we shall only apply Definition~\ref{defn:rand_comp} to classes $\class{G}$ that contain finitely many isomorphism types of $n$-vertex graphs for each $n$. 
	
	\subsection{From distributed algorithms to descriptive combinatorics}\label{subsec:dist_res}
	
	\subsubsection{Deterministic algorithms and Borel colorings}
	
	We use standard descriptive set-theoretic terminology; see, e.g., Kechris's book \cite{KechrisDST} or Tserunyan's notes \cite{AnushDST}. Recall that a separable topological space is \emphd{Polish} if its topology is induced by a complete metric. A \emphd{standard Borel space} is a set $X$ equipped with a $\sigma$-algebra $\mathfrak{B}(X)$ of \emphd{Borel sets} generated by a Polish topology on $X$. We say that a Polish topology on a standard Borel space $X$ is \emphd{compatible} if it generates $\mathfrak{B}(X)$.

	By a \emphd{Borel graph} we mean a graph $G$ whose vertex set $V(G)$ is a standard Borel space and such that $\set{(x,y) \,:\, \text{$x$ and $y$ are adjacent in $G$}}$ is a Borel subset of $V(G) \times V(G)$. Similarly, a \emphd{Borel structured graph} is a structured graph $\bm{G} = (G, \sigma)$ such that $G$ is a Borel graph and $\sigma \colon \finseq{V(G)} \rightharpoonup \N$ is a Borel function \ep{which just means that for every $c \in \N$, $\sigma^{-1}(c)$ is a Borel subset of $\finseq{V(G)}$}.
	
	Let $\class{G} \subseteq \class{FSG}$ be a collection of isomorphism types of finite structured graphs and let $R$, $n \in \N$. We say that a structured graph $\bm{G}$ is \emphd{$(R, n)$-locally in $\class{G}$} if for each vertex $x \in V(\bm{G})$, there is an $n$-vertex structured graph $\bm{H}$ with $[\bm{H}] \in \class{G}$ such that $[B_{\bm{G}}(x, R), x] = [B_{\bm{H}}(y, R), y]$ for some $y \in V(\bm{H})$. For instance, suppose that $\class{G}$ is a class of finite graphs \ep{with no other structure} closed under %
	adding isolated vertices. Examples of such classes are the class of all planar graphs, the class of all graphs of maximum degree at most $\Delta$ for some fixed $\Delta \in \N$, the class of all triangle-free graphs, etc. It is not hard to see that in this case a graph $G$ is $(R,n)$-locally in $\class{G}$ provided that %
	for all $x \in V(G)$, $[B_G(x,R)] \in \class{G}$ and $|B_G(x,R)| \leq n$. Now we can state our first result:
	
	\begin{theo}[\textls{Deterministic \LOCAL algorithms yield Borel colorings}]\label{theo:dist_Borel}
		Let $\Pi = (t, \mathcal{P})$ be a local coloring problem and let $\class{G} \subseteq \class{FSG}$. Fix $n \in \N^+$ such that $\Det_{\Pi, \class{G}}(n) < \infty$ and set $R \defeq \Det_{\Pi, \class{G}}(n) + t$. If $\bm{G}$ is a Borel structured graph that is $(R,n)$-locally in $\class{G}$ and such that $|B_{\bm{G}}(x, 2R)| \leq n$ for all $x \in V(\bm{G})$, then $\bm{G}$ has a Borel $\Pi$-coloring. %
	\end{theo} 
	
	Note that a structured graph $\bm{G}$ satisfying the assumptions of Theorem~\ref{theo:dist_Borel} must have finite maximum degree, since $|B_{\bm{G}}(x, R)| \leq n < \infty$ for all $x \in V(\bm{G})$.
	
	Usually the precise value of $\Det_{\Pi, \class{G}}(n)$ is not known and one only has access to asymptotic upper bounds. For instance, assume that $\class{G}$ is a class of finite graphs \ep{with no other structure}, and let $G$ be a Borel graph of finite maximum degree $\Delta$ all of whose finite induced subgraphs are in $\class{G}$. We claim that if $\Pi = (t, \mathcal{P})$ is a local coloring problem such that $\Det_{\Pi, \class{G}}(n) = o(\log n)$, then $G$ has a Borel $\Pi$-coloring. Indeed, we may assume that $G$ is infinite (otherwise $G$ itself is a finite graph in $\class{G}$, and hence it has an---automatically Borel---$\Pi$-coloring by assumption). This implies that $\class{G}$ is closed under adding isolated vertices. Since $t$ is a constant independent of $n$, $R(n) \defeq \Det_{\Pi, \class{G}}(n) + t = o(\log n)$ as well. Therefore, for each $x \in V(G)$,
	\begin{equation}\label{eq:o(1)_power}
		|B_G(x, R(n))| \,\leq\, |B_G(x, 2R(n))| \,\leq\, 1 + \Delta^{2R(n)} \,=\, 1 + \Delta^{o(\log n)} \,=\, n^{o(1)}.
	\end{equation}
	Thus, Theorem~\ref{theo:dist_Borel} may be applied for any large enough value of $n$.
	
	It is possible to incorporate into such calculations extra assumptions on the growth rate of $G$. For example, say that a graph $G$ is of \emphd{subexponential growth} if for each $\epsilon > 0$, there is $R_0 \in \N$ such that if $x \in V(G)$ and $R \geq R_0$, then $|B_G(x, R)| < (1+\epsilon)^R$. As before, let $\class{G}$ be a class of finite graphs and let $G$ be a Borel graph of subexponential growth all of whose finite induced subgraphs are in $\class{G}$. The same calculation as in the last paragraph shows that $G$ has a Borel $\Pi$-coloring whenever $\Pi$ is a local coloring problem with $\Det_{\Pi, \class{G}}(n) = O(\log n)$. Stronger assumptions on the growth rate of $G$ enable one to further relax the necessary bound on $\Det_{\Pi, \class{G}}(n)$.
	
	The reasoning in the preceding two paragraphs will apply, almost verbatim, to all the other	results in this section. %
	
	\subsubsection{Deterministic algorithms and continuous colorings}
	
	Under certain circumstances, the word ``Borel'' in the conclusion of Theorem~\ref{theo:dist_Borel} can be replaced by ``continuous.'' In order to state this result precisely, we require a few definitions. In what follows, we view $\N$ as a discrete topological space. Let $(X, d)$ be a metric space and let $\bm{G}_1$, $\bm{G}_2$ be finite structured graphs with $V(\bm{G}_1)$, $V(\bm{G}_2) \subseteq X$. Given $\epsilon > 0$, we say that $\bm{G}_1$ and $\bm{G}_2$ are \emphd{$\epsilon$-isomorphic} \ep{with respect to the metric $d$} if there is an isomorphism $\phi \colon V(\bm{G}_1) \to V(\bm{G}_2)$ between $\bm{G}_1$ and $\bm{G}_2$ such that $d(x, \phi(x)) < \epsilon$ for all $x \in V(\bm{G}_1)$. Recall that a topological space is \emphd{zero-dimensional} is it has a base consisting of clopen sets.
	
	\begin{defn}[\textls{Topological graphs}]\label{defn:top_graph}
		A \emphd{topological structured graph} is a locally finite structured graph $\bm{G}$ whose vertex set $V(\bm{G})$ is a zero-dimensional Polish space and that has the following property for some \ep{hence any} metric $d$ inducing the topology on $V(\bm{G})$: For each $x \in V(\bm{G})$, $R \in \N$, and $\epsilon > 0$, there is $\delta > 0$ such that if $y \in V(\bm{G})$ satisfies $d(x, y) < \delta$, then the rooted structured graphs $(B_{\bm{G}}(x, R), x)$ and $(B_{\bm{G}}(y, R), y)$ are $\epsilon$-isomorphic.
	\end{defn}
	
	Note that the above definition is indeed independent of the choice of the metric $d$ \ep{we could have stated it purely topologically, by fixing an open neighborhood for each vertex in $B_{\bm{G}}(x, R)$ instead of using the parameter $\epsilon$}. If we are hoping to find continuous colorings of $\bm{G}$, it is natural to assume that the topology on $V(\bm{G})$ is zero-dimensional, so that $V(\bm{G})$ has enough clopen subsets and hence there are many continuous functions $V(\bm{G}) \to \N$. The rest of Definition~\ref{defn:top_graph} describes the interplay between the combinatorics of $\bm{G}$ and the topology on $V(\bm{G})$. In particular, if $\bm{G}$ is a topological graph, then the map $V(\bm{G}) \to \class{FSG}_\bullet \colon x \mapsto [B_{\bm{G}}(x,R), x]$ is continuous for all $R \in \N$, where we view the countable set $\class{FSG}_\bullet$ as a discrete space.
	
	We prove a continuous version of Theorem~\ref{theo:dist_Borel} for topological graphs:
	
	\begin{theo}[\textls{Deterministic \LOCAL algorithms yield continuous colorings}]\label{theo:dist_cont}
		Let $\Pi = (t, \mathcal{P})$ be a local coloring problem and let $\class{G} \subseteq \class{FSG}$. Fix $n \in \N^+$ such that $\Det_{\Pi, \class{G}}(n) < \infty$ and set $R \defeq \Det_{\Pi, \class{G}}(n) + t$. If $\bm{G}$ is a topological structured graph that is $(R,n)$-locally in $\class{G}$ and such that $|B_{\bm{G}}(x, 2R)| \leq n$ for all $x \in V(\bm{G})$, then $\bm{G}$ has a continuous $\Pi$-coloring. %
	\end{theo} 
	
	Observe that Theorem~\ref{theo:dist_cont} is a strengthening of Theorem~\ref{theo:dist_Borel}. Indeed, it follows from standard results in descriptive set theory that if $\bm{G}$ is a locally finite Borel structured graph, then there is a compatible zero-dimensional Polish topology on $V(\bm{G})$ that makes $\bm{G}$ a topological structured graph \cite[\S13]{KechrisDST}. \ep{Nevertheless, we will give a direct proof of Theorem~\ref{theo:dist_Borel} as well.}
	
	An ample supply of natural examples of topological graphs is provided by continuous actions of finitely generated groups. Let $X$ be a zero-dimensional Polish space and let $\G$ be a group generated by a finite set $S \subseteq \G$. Assume additionally that $S$ does not contain the identity element of $\G$ and that $S$ is \emphd{symmetric}, meaning that if $\gamma \in S$, then $\gamma^{-1} \in S $ as well. Given a free action $\alpha \colon \G \acts X$ of $\G$ on $X$ by homeomorphisms, let the \emphd{Schreier graph} of $\alpha$ be the graph $G_\alpha$ with vertex set $X$ in which two vertices $x$, $y \in X$ are adjacent if and only if $x = \gamma \cdot_\alpha y$ for some $\gamma \in S$. Because $\alpha$ is free, it is easy to see that $G_\alpha$ is a topological graph. If $S$ is not assumed to be symmetric, one can similarly define a directed version of the Schreier graph, which is a topological directed graph.
	
	\subsubsection{Randomized algorithms}
	
	Recall that a set or a function is \emphd{measurable} \ep{resp.~\emphd{Baire\-/measurable}} if it differs from a Borel set or function on a null \ep{resp.~meager} set. In order to obtain measurable or Baire-measurable colorings, it is enough to work with randomized rather than deterministic \LOCAL algorithms:
	
	\begin{theo}[\textls{Randomized \LOCAL algorithms yield measurable colorings}]\label{theo:dist_meas}
		Let $\Pi = (t, \mathcal{P})$ be a local coloring problem and let $\class{G} \subseteq \class{FSG}$. Fix $n \in \N^+$ such that $\Rand_{\Pi, \class{G}}(n) < \infty$ and set $R \defeq \Rand_{\Pi, \class{G}}(n) + t$. Let $\bm{G}$ be a Borel structured graph that is $(R,n)$-locally in $\class{G}$ and such that $|B_{\bm{G}}(x, 2R)| \leq n^{1/8}/4$ for all $x \in V(\bm{G})$. Then the following conclusions hold:
		\begin{enumerate}[label=\ep{\normalfont\roman*}]
			\item If $\mu$ is a probability Borel measure on $V(\bm{G})$, then $\bm{G}$ has a $\mu$-measurable $\Pi$-coloring.
			
			\item If $\tau$ is a compatible Polish topology on $V(\bm{G})$, then $\bm{G}$ has a $\tau$-Baire-measurable $\Pi$-coloring.
		\end{enumerate}
	\end{theo}
	
	Although the bound on $|B_{\bm{G}}(x, 2R)|$ required in Theorem~\ref{theo:dist_meas} is stronger than that in Theorem~\ref{theo:dist_Borel}, it remains polynomial in $n$. In particular, inequalities \eqref{eq:o(1)_power} still show that if $\class{G}$ is a class of finite graphs 
    and 
    $G$ is a Borel graph of finite maximum degree all of whose finite induced subgraphs are in $\class{G}$, then $G$ is measurably and Baire-measurably $\Pi$-colorable provided that $\Rand_{\Pi,\class{G}}(n) = o(\log n)$. Furthermore, if $G$ is of subexponential growth, then $\Rand_{\Pi,\class{G}}(n) = O(\log n)$ suffices. Actually, in this case we can even make the coloring Borel:
	
	\begin{theo}[\textls{Randomized \LOCAL algorithms yield Borel colorings of subexponential growth graphs}]\label{theo:dist_subexp}
		Let $\Pi = (t, \mathcal{P})$ be a local coloring problem and let $\class{G} \subseteq \class{FSG}$. Fix $n \in \N^+$ such that $\Rand_{\Pi, \class{G}}(n) < \infty$ and set $R \defeq \Rand_{\Pi, \class{G}}(n) + t$.
		If $\bm{G}$ is a Borel structured graph of subexponential growth that is $(R,n)$-locally in $\class{G}$ and such that $|B_{\bm{G}}(x, 2R)| \leq n/e$ for all $x \in V(\bm{G})$, then $\bm{G}$ has a Borel $\Pi$-coloring. %
	\end{theo}	
	
	In the statement of Theorem~\ref{theo:dist_subexp}, $e= 2.71\ldots$ is the base of the natural logarithm. Theorem~\ref{theo:dist_subexp} contributes to the growing body of results showing that various combinatorial problems on graphs of subexponential growth can be solved in a Borel way, whereas in general one can only hope for a measurable solution; see, e.g., \cite{CGMPT, ConleyTamuz, Thornton}.
	
	\subsection{The Lov\'asz Local Lemma}\label{subsec:LLL}
	
	As mentioned in the \hyperref[sec:intro]{introduction}, the key ingredient in our proof of Theorem~\ref{theo:dist_meas} is a new measurable version of the Lov\'asz Local Lemma \ep{the LLL for short}, which is an interesting result in its own right. Indeed, the original motivation for this paper was to develop a better understanding of the behavior of the LLL in the measurable setting.
	
	\subsubsection{Constraint satisfaction problems and the LLL}\label{subsubsec:LLL}
	
	Recall that for a set $A$, $\fins{A}$ denotes the set of all finite subsets of $A$. Similarly, for sets $A$ and $B$, let $\finf{A}{B}$ be the set of all partial functions $\phi \colon A \rightharpoonup B$ with finite domains.
	
	\begin{defn}[\textls{CSPs}]
		Fix a set $X$ and $m \in \N^+$. An \emphd{$(X,m)$-constraint} \ep{or simply a \emphd{constraint} if $X$ and $m$ are clear from the context} is a set $B \subseteq \finf{X}{[m]}$ such that $\dom(\phi) = \dom(\psi)$ for all $\phi$, $\psi \in B$. If a constraint $B$ is nonempty, then its \emphd{domain} is the set $\dom(B) \coloneqq \dom(\phi)$ for some \ep{hence all} $\phi \in B$; the domain of the empty constraint is defined to be $\0$. A \emphd{constraint satisfaction problem} \ep{a \emphd{CSP} for short} $\mathscr{B}$ on $X$ with range $[m]$, in symbols $\mathscr{B} \colon X \to^? [m]$, is a set of $(X,m)$-constraints. A \emphd{solution} to a CSP $\mathscr{B} \colon X \to^? [m]$ is a function $f \colon X \to [m]$ such that for all $B \in \mathscr{B}$, the restriction $\rest{f}{\dom(B)}$ of $f$ onto $\dom(B)$ is not a member of $B$.
	\end{defn}
	
	In other words, in a CSP $\mathscr{B} \colon X \to^? [m]$, each constraint $B \in \mathscr{B}$ is interpreted as a set of ``forbidden patterns'' that are not allowed to appear in a solution $f \colon X \to [m]$. There are obvious similarities between CSPs in the above sense and local colorings problems in the sense of Definition~\ref{defn:LCL}, and, indeed, a CSP can be viewed as local coloring problem on an auxiliary graph \ep{see \S\ref{subsec:LLL_to_dist} for details}. %
	
	Fix a CSP $\mathscr{B} \colon X \to^? [m]$. The \emphd{probability} $\P[B]$ of a constraint $B \in \mathscr{B}$ is defined by
	\[
	\P[B] \,\defeq\, \frac{|B|}{m^{|\dom(B)|}}.
	\] %
	Notice that if we form a random coloring $f \colon X \to [m]$ by assigning to each $x \in X$ a color from $[m]$ uniformly at random, then $\P[B]$ equals the probability that $\rest{f}{\dom(B)} \in B$, i.e., that the constraint $B$ is {violated} by $f$. The \emphd{neighborhood} of a constraint $B \in \mathscr{B}$ is the set $\Nbhd(B) \subset \mathscr{B}$ given by
	\[
	\Nbhd(B) \,\defeq\, \set{B' \in \mathscr{B} \setminus \set{B} \,:\, \dom(B') \cap \dom(B) \neq \0}.
	\]
	The most widely used form of the LLL, called the \emph{Symmetric LLL}, invokes the parameters
	\[
	\pr(\mathscr{B}) \,\defeq\, \sup_{B \in \mathscr{B}} \P[B] \qquad \text{and} \qquad \de(\mathscr{B}) \,\defeq\, \sup_{B \in \mathscr{B}} |\Nbhd(B)|.
	\]

	\begin{theo}[{\textls{Symmetric LLL}; see \cite[Corollary 5.1.2]{AS}}]\label{theo:SLLL}
		If $\mathscr{B}$ is a CSP such that
		\begin{equation}\label{eq:SLLL}
			\pr(\mathscr{B}) \cdot (\de(\mathscr{B}) + 1) \,\leq\, e^{-1},
		\end{equation}
		where $e = 2.71\ldots$ is the base of the natural logarithm, then $\mathscr{B}$ has a solution.
	\end{theo}
	
	Theorem~\ref{theo:SLLL} is a special case of a stronger result, known as the \emph{General LLL}, in which instead of bounding $\pr(\mathscr{B})$ and $\de(\mathscr{B})$ uniformly, one establishes a more delicate---but somewhat less transparent---relationship between $\P[B]$ and $|\Nbhd(B)|$ for each constraint $B \in \mathscr{B}$:
	
	\begin{theo}[{\textls{General LLL}; see \cite[Theorem 5.1.1]{AS}}]\label{theo:GLLL}
		If $\mathscr{B}$ is a CSP such that there is a function $\eta \colon \mathscr{B} \to [0,1)$ satisfying
		\[
		\P[B] \,\leq\, \eta(B) \prod_{B' \in \Nbhd(B)} (1 - \eta(B')), \qquad \text{for all } B \in \mathscr{B},
		\]
		then $\mathscr{B}$ has a solution.
	\end{theo}
	
	A standard calculation \ep{see \cite[proof of Corollary 5.1.2]{AS}} shows that the bound \eqref{eq:SLLL} implies the existence of a function $\eta$ as in Theorem~\ref{theo:GLLL}, and hence the Symmetric LLL is indeed a special case of the General LLL.
	We remark that, due to its origin in finite combinatorics, the LLL is often stated in the case when the ground set $X$ is finite; however, the case of infinite $X$ follows via a straightforward compactness argument \ep{see, e.g., \cite[proof of Theorem 5.2.2]{AS}}.
	
	Several applications of the LLL in combinatorics and graph theory can be found in \cite{AS, MR}. For most applications, the full power of the General LLL is not needed and the Symmetric LLL is sufficient. Furthermore, in many cases the bound \eqref{eq:SLLL} is quite far from being sharp: one can often prove that $\pr(\mathscr{B})$ is at most $\exp(-\log^{1+\epsilon} \de(\mathscr{B}))$ or even $\exp(-\de(\mathscr{B})^\epsilon)$ for some constant $\epsilon > 0$. For example, this usually happens when the upper bound on $\pr(\mathscr{B})$ is obtained via a concentration of measure argument.
	
	\subsubsection{Measurable Symmetric LLL}
	
	If $X$ is a standard Borel space, then $\fins{X}$ also carries a natural standard Borel structure. Namely, a subset of $\fins{X}$ is said to be Borel if and only if its preimage in $\bigsqcup_{k = 0}^\infty X^k$ under the map $(x_0, \ldots, x_{k-1}) \mapsto \set{x_0, \ldots, x_{k-1}}$ is Borel. Next, for any standard Borel space $Y$, the set $\finf{X}{Y}$ also becomes a standard Borel space when viewed as a Borel subset of $\fins{X\times Y}$ by identifying each partial function $\phi \in \finf{X}{Y}$ with its graph, i.e., with the set $\set{(x, \phi(x)) \,:\, x \in \dom(\phi)}$. Finally, since, by definition, every $(X,m)$-constraint is a finite set of partial functions with a common finite domain, the set of all $(X,m)$-constraints is a Borel subset of $\fins{\finf{X}{[m]}}$. Therefore, we may speak of \emphd{Borel CSPs} $\mathscr{B} \colon X \to^? [m]$, i.e., Borel sets of $(X,m)$-constraints.
	
	Given a Borel CSP, it is natural to ask whether it has a solution with some regularity properties. To that end, we establish the following measurable/Baire-measurable analog of Theorem~\ref{theo:SLLL} under a stronger polynomial bound on $\pr(\mathscr{B})$ and $\de(\mathscr{B})$:
	
	\begin{theo}[\textls{Measurable Symmetric LLL}]\label{theo:meas_SLLL}
		Let $\mathscr{B} \colon X \to^? [m]$ be a Borel CSP such that
		\begin{equation}\label{eq:meas_SLLL}
			\pr(\mathscr{B}) \cdot (\de(\mathscr{B})+1)^8 \,\leq\, 2^{-15}.
		\end{equation}
		Assume additionally that $\sup \set{|\dom(B)| \,:\, B \in \mathscr{B}} < \infty$. Then the following conclusions hold:
		\begin{enumerate}[label=\ep{\upshape\roman*}]
			\item\label{item:SLLLmeas} If $\mu$ is a probability Borel measure on $X$, then $\mathscr{B}$ has a $\mu$-measurable solution.
			
			\item\label{item:SLLLBaire} If $\tau$ is a compatible Polish topology on $X$, then $\mathscr{B}$ has a $\tau$-Baire-measurable solution.
		\end{enumerate}
	\end{theo}
	
	As mentioned earlier, one can often bound $\pr(\mathscr{B})$ from above by a super-polynomially small function of $\de(\mathscr{B})$, making Theorem~\ref{theo:meas_SLLL} applicable. We call CSPs $\mathscr{B}$ such that $\sup \set{|\dom(B)| \,:\, B \in \mathscr{B}} < \infty$ \emphd{bounded}. The boundedness assumption in Theorem~\ref{theo:meas_SLLL} can likely be eliminated by a routine modification of the proof. However, we are not aware of any applications where this assumption is not satisfied, and so we elected to include it in order to make our arguments more transparent.
	
	Theorem~\ref{theo:meas_SLLL} is implied by Theorem~\ref{theo:dist_meas}: a CSP $\mathscr{B}$ can be encoded as a local coloring problem on %
	an auxiliary graph \ep{see \S\ref{subsec:LLL_to_dist} for details}, and Fischer and Ghaffari \cite{FG} designed a randomized \LOCAL algorithm that finds solutions to such local coloring problems in a sublogarithmic number of rounds under a polynomial bound on $\pr(\mathscr{B})$ and $\de(\mathscr{B})$ similar to \eqref{eq:meas_SLLL}. The actual bound that we use %
	comes from a sharpened version of the Fischer--Ghaffari algorithm developed in \cite{GHK} by Ghaffari, Harris, and Kuhn. Together with Theorem~\ref{theo:dist_meas}, their algorithm immediately yields Theorem~\ref{theo:meas_SLLL}. The logic of our argument, however, goes in the opposite direction: We shall {first} prove Theorem~\ref{theo:meas_SLLL} and then derive Theorem~\ref{theo:dist_meas} from it. Nevertheless, our proof of Theorem~\ref{theo:meas_SLLL} does invoke the Ghaffari--Harris--Kuhn algorithm, albeit in a more subtle way.
	
	\subsubsection{Comparison with prior work}

	{
		\renewcommand{\arraystretch}{1.4}
		\begin{table}[t]
			\caption{What is known about the LLL in the descriptive setting.}\label{table:LLL}
			\begin{tabular}{| c || >{\centering\arraybackslash} m{7cm} | >{\centering\arraybackslash} m{5cm}|}
				\hline & Symmetric LLL & General LLL \\\hline\hline
				Borel? & YES for subexponential growth \cite{CGMPT}
				
				NO in general \cite{CJMST-D} & \ep{NO in general}\\\hline\hline
				Measurable? & 
				\emphd{YES when $\pr(\de+1)^8 \leq 2^{-15}$} \ep{this paper} & YES for $[0,1]$-shifts \cite{MLLL}
				
				\ep{NO in general} \\\hline
				With $\epsilon$ error? & YES \cite{MLLL} & NO in general \cite{MLLL}\\\hline\hline
				Baire-measurable? &
				\emphd{YES when $\pr(\de+1)^8 \leq 2^{-15}$} \ep{this paper} & \emphd{???}\\\hline
			\end{tabular}
		\end{table}	
	}
	
	The current state of the knowledge concerning the behavior of the LLL in the descriptive setting is summarized in Table~\ref{table:LLL}. Conley, Jackson, Marks, Seward, and Tucker-Drob \cite[Theorem~1.6]{CJMST-D} constructed examples showing that the Symmetric LLL cannot, in general, produce Borel solutions. On the other hand, Cs\'{o}ka, Grabowski, M\'ath\'e, Pikhurko, and Tyros \cite{CGMPT} showed that a Borel version of the Symmetric LLL holds under certain subexponential growth assumptions on the CSP $\mathscr{B}$ \ep{and we use their result to derive Theorem~\ref{theo:dist_subexp}}.

	In the presence of a measure $\mu$, one can relax the requirements and only ask for a measurable function $f \colon X \to [m]$ that satisfies the constraints on a set of measure $1 - \epsilon$, for any given $\epsilon > 0$. In this regime, the Symmetric LLL always succeeds \cite[Theorem~5.1]{MLLL}, while the General LLL may fail \cite[Theorem 7.1]{MLLL}. In fact, the General LLL may even fail to produce colorings that satisfy the constrains on a set of arbitrarily small
	positive measure. On the other hand, there are certain situations of particular interest in ergodic theory when the General LLL can be used to obtain measurable colorings satisfying all the constraints; namely, this happens when the structure of the CSP $\mathscr{B}$ is, in a certain technical sense, ``induced'' by the Bernoulli shift action $\G \acts [0,1]^\G$ of a countable group $\G$ \cite[Theorem 6.6]{MLLL}.
	
	Theorem~\ref{theo:meas_SLLL} provides the first Baire-measurable variant of the LLL \ep{with the exception of the subexponential Borel LLL of Cs\'{o}ka--Grabowski--M\'ath\'e--Pikhurko--Tyros}. The following question remains open:
	
	\begin{prob}
		Does there exist a Baire-measurable version of the General LLL?
	\end{prob}
	
	\section{Applications}\label{sec:applications}
	
	\subsection{Colorings with the number of colors close to $\Delta$}
	
	Recall that the \emphd{chromatic number} $\chi(G)$ of a graph $G$ is the smallest cardinality of a set $Y$ such that $G$ admits a proper coloring $f \colon V(G) \to Y$. In the descriptive setting, one defines the \emphd{Borel chromatic number} $\chi_\mathrm{B}(G)$ of a Borel graph $G$ as the smallest cardinality of a standard Borel space $Y$ such that $G$ has a Borel proper coloring $f \colon V(G) \to Y$. Similarly, given a probability Borel measure $\mu$ or a compatible Polish topology $\tau$ on $V(G)$, the \emphd{$\mu$-measurable chromatic number} $\chi_\mu(G)$ and the \emphd{$\tau$-Baire-measurable chromatic number} $\chi_\tau(G)$ are defined to be the smallest cardinality of a standard Borel space $Y$ such that $G$ admits a $\mu$-measurable, resp.~$\tau$-Baire-measurable, proper coloring $f \colon V(G) \to Y$. It is clear that $\chi(G) \leq \chi_\mu(G)$, $\chi_\tau(G) \leq \chi_\mathrm{B}(G)$. We shall only work with bounded degree graphs, and for them all these parameters are finite:
	
	\begin{theo}[{Kechris--Solecki--Todorcevic \cite[Proposition 4.6]{KechrisSoleckiTodorcevic}}]\label{theo:KST}
		Let $G$ be a Borel graph of finite maximum degree $\Delta$. Then $\chi_{\mathrm{B}}(G) \leq \Delta+1$.
	\end{theo}
	
	In general, the bound $\chi(G) \leq \Delta(G) + 1$ is best possible, since a graph of maximum degree $\Delta$ may contain $\Delta + 1$ pairwise adjacent vertices \ep{i.e., a \emphd{$(\Delta+1)$-clique}}, all of which would have to receive distinct colors in a proper coloring of $G$. Conversely, a classical theorem of Brooks \ep{see \cite[Theorem 5.2.4]{Die}} asserts that if $\chi(G) = \Delta(G) + 1$ for a graph $G$ with $3 \leq \Delta(G) < \infty$, then $G$ contains a $(\Delta(G)+1)$-clique.
	
	It is natural to ask whether Brooks's theorem can be extended to the setting of Borel, measurable, or Baire-measurable colorings. Marks \cite{Marks} showed that, in the Borel context, Brooks's theorem fails in a very strong sense. A graph $G$ is called \emphd{$\Delta$-regular} if every vertex in $G$ has degree $\Delta$, and \emphd{acyclic} if it contains no cycles \ep{acyclic graphs are also known as \emphd{forests}}. It is easy to see that acyclic graphs are bipartite, i.e., have chromatic number at most $2$. Nevertheless, Marks proved the following:
	
	\begin{theo}[Marks~\cite{Marks}]\label{theo:Marks}
		For each $\Delta \in \N$, there exists an acyclic $\Delta$-regular Borel graph $G$ such that $\chi_{\mathrm{B}}(G) = \Delta + 1$.
	\end{theo}
	
	In contrast to Theorem~\ref{theo:Marks}, Conley, Marks, and Tucker-Drob \cite{CMTD} succeeded in extending Brooks's theorem to the setting of measurable and Baire-measurable colorings:
	
	\begin{theo}[{\textls{Measurable Brooks}; Conley--Marks--Tucker-Drob \cite{CMTD}}]\label{theo:meas_Brooks}
		Let $G$ be a Borel graph of finite maximum degree $\Delta \geq 3$ without a $(\Delta + 1)$-clique.
		\begin{enumerate}[label=\ep{\normalfont\roman*}]
			\item If $\mu$ is a probability Borel measure on $V(G)$, then $\chi_\mu(G) \leq \Delta$.
			
			\item If $\tau$ is a compatible Polish topology on $V(G)$, then $\chi_\tau(G) \leq \Delta$.
		\end{enumerate}
	\end{theo}
	
	In view of Brooks's theorem, an equivalent way of phrasing Theorem~\ref{theo:meas_Brooks} is that a Borel graph $G$ of finite maximum degree $\Delta \geq 3$ is measurably/Baire-measurably $\Delta$-colorable if and only if it is $\Delta$-colorable abstractly, i.e., without any regularity restrictions. Can this result be extended to colorings with fewer colors? We show that the answer is positive for $(\Delta - c)$-colorings, where $c$ can be as large as roughly $\sqrt{\Delta}$:
	
	\begin{theo}\label{theo:Delta}
		There is $\Delta_0 \in \N$ with the following property. Fix integers $\Delta \geq \Delta_0$ and $c$. Let $G$ be a Borel graph of maximum degree at most $\Delta$ and let $\mu$ \ep{resp.~$\tau$} be a probability Borel measure \ep{resp.~a compatible Polish topology} on $V(G)$. If %
		$c \leq \sqrt{\Delta + 1/4} - 5/2$, then the following statements are equivalent: %
		\begin{enumerate}[label=\ep{\normalfont\roman*}]
			\item\label{item:comb} $\chi(G) \leq \Delta - c$;
			
			\item\label{item:meas} $\chi_\mu(G) \leq \Delta - c$;
			
			\item\label{item:BM} $\chi_\tau(G) \leq \Delta - c$.
		\end{enumerate}
	\end{theo}
	\begin{scproof}
		Implications \ref{item:meas}, \ref{item:BM} $\Longrightarrow$ \ref{item:comb} are trivial, so we only need to argue that \ref{item:comb} $\Longrightarrow$ \ref{item:meas}, \ref{item:BM}. For each $k \in \N$, let $\Pi(k)$ denote the local coloring problem that encodes proper $k$-coloring of graphs \ep{see Example~\ref{exmp:prop_col}}. Fix $\Delta \in \N$ and $c \in \Z$ with $c \leq \sqrt{\Delta + 1/4} - 5/2$ and let $\class{G}(\Delta, c)$ be the set of all isomorphism classes of finite graphs $G$ of maximum degree at most $\Delta$ satisfying $\chi(G) \leq \Delta - c$. Assuming $\Delta$ is large enough, Bamas and Esperet \cite[Theorem 1.3]{BamasEsperet} established the following bounds on the randomized \LOCAL complexity of proper $(\Delta-c)$-coloring of graphs in $\class{G}(\Delta, c)$:
		\begin{equation}\label{eq:Delta_bound}
			\Rand_{\Pi(\Delta - c), \class{G}(\Delta, c)}(n) \,\leq\, \exp(O(\sqrt{\log \log n})) \,=\, o(\log n).
		\end{equation}
		\ep{Here the implicit constants in the asymptotic notation may depend on $\Delta$, which we treat as fixed.} If $G$ is a Borel graph with $\Delta(G) \leq \Delta$ and $\chi(G) \leq \Delta -c$, then every finite induced subgraph of $G$ is in $\class{G}(\Delta, c)$. Therefore, Theorem~\ref{theo:dist_meas}, combined with \eqref{eq:Delta_bound}, yields that $G$ is both measurably and Baire-measurably $(\Delta-c)$-colorable, as desired.
	\end{scproof}
	
	Note that if in the statement of Theorem~\ref{theo:Delta} we additionally assume that $G$ is of subexponential growth, then, replacing Theorem~\ref{theo:dist_meas} by Theorem~\ref{theo:dist_subexp}, we may conclude that for such $G$, statements \ref{item:comb}--\ref{item:BM} are also equivalent to
	\begin{enumerate}[label=\ep{\normalfont\roman*}]
		\setcounter{enumi}{3}
		\item $\chi_{\mathrm{B}}(G) \leq \Delta - c$.
	\end{enumerate}
	Marks's Theorem~\ref{theo:Marks} shows that the subexponential growth assumption cannot be removed.
	
	The distributed algorithm of Bamas and Esperet \cite{BamasEsperet} from which we derive Theorem~\ref{theo:Delta} is inspired by the earlier work of Molloy and Reed \cite{MolloyReedDelta} on sequential algorithms for $(\Delta - c)$-coloring. Among other results, Molloy and Reed established the following remarkable extension of Brooks's theorem: If $\Delta$ is large enough, $G$ is a graph of maximum degree at most $\Delta$, and $c \leq \sqrt{\Delta + 1/4} - 5/2$, then $\chi(G) \leq \Delta - c$ if and only if $\chi(B_G(x, 1)) \leq \Delta - c$ for all $x \in V(G)$ \cite[Theorem 5]{MolloyReedDelta}.

	The bound on $c$ in the statement of Theorem~\ref{theo:Delta} is almost sharp, in the sense that it cannot be relaxed to $c \leq \sqrt{\Delta - 3/4} - 1/2$:
	
	\begin{prop}\label{prop:sharp}
		Let $k$, $\Delta \in \N$ be such that $\Delta \geq k \geq 2$ and $\Delta - k \geq \sqrt{\Delta - 3/4} - 1/2$. Then there exist a Borel graph $G$ and a probability Borel measure $\mu$ \ep{resp.~a compatible Polish topology $\tau$} on $V(G)$ such that $\Delta(G) \leq \Delta$ and $\chi(G) \leq k$ but $\chi_\mu(G) > k$ \ep{resp.~$\chi_\tau(G) > k$}.
	\end{prop}
	\begin{scproof}
		We adapt the proof of \cite[Theorem 1.4]{EWHK} due to Embden-Weinert, Hougardy, and Kreuter. We shall give the argument for the measurable chromatic number, the Baire-measurable case being virtually identical. Fix $k \geq 2$. To begin with, observe that there exists a Borel graph $G$ of finite maximum degree with a probability Borel measure $\mu$ on $V(G)$ such that $\chi(G) = k$ but $\chi_\mu(G) > k$. For instance, let $(X, \mu)$ be a standard probability space and let $T \colon X \to X$ be a measure\-/preserving transformation such that for all $n \in \Z\setminus \set{0}$, $T^n$ is ergodic and $T^n(x) \neq x$ for all $x \in X$. An example of such a transformation is the map $T \colon [0,1) \to [0,1) \colon x \mapsto (x+\alpha) \,\mathrm{mod}\, 1$, where $\alpha$ is a fixed irrational number. Let $G$ be the graph with vertex set $X$ in which two distinct vertices $x$ and $y$ are adjacent if and only if $y = T^n(x)$ with $|n| \leq k-1$. It is easy to verify that $\chi(G) = k$ and every proper $k$-coloring $f \colon X \to [k]$ of $G$ is $T^k$-periodic, meaning that $f(T^k(x)) = f(x)$ for all $x \in X$. To show that $\chi_\mu(G) > k$, suppose that $f \colon X \to [k]$ is a measurable proper $k$-coloring. Then every color class of $f$ is a $T^k$-invariant measurable subset of $X$. Since the map $T^k$ is ergodic, this implies that precisely one color class of $f$ is conull, while the other $k-1$ classes are null. This is a contradiction, as the map $T$ cyclically permutes the color classes, showing that they all must have the same measure.
		
		Now let $\Delta \geq k$ satisfy $\Delta - k \geq \sqrt{\Delta - 3/4} - 1/2$ and let $G$ be a Borel graph with a probability Borel measure $\mu$ on $V(G)$ such that $\chi(G) \leq k$, $\chi_\mu(G) > k$, and $\Delta(G)$ is the smallest among all graphs with these properties. Set $d \defeq \Delta(G)$ and $c \defeq d - k$. Suppose, toward a contradiction, that $d > \Delta$. Then $c > \sqrt{d - 3/4} - 1/2$, which implies that $c(c+1) > d - 1$. Since $c$ and $d$ are integers, this yields $c(c+1) \geq d$. Fix a Borel linear order $<$ on $V(G)$ \ep{such an order exists since, by \cite[Theorem 15.6]{KechrisDST}, $V(G)$ is isomorphic to a Borel subset of $\R$}. For each $x \in X$ and $1 \leq i \leq \deg_G(x)$, let $N_i(x)$ be the $i$-th neighbor of $x$ in the order $<$; that is, we have
		\[
		N_G(x) \,=\, \set{N_1(x), \ldots, N_{\deg_G(x)}(x)} \qquad \text{and} \qquad N_1(x) \,<\, \cdots \,<\, N_{\deg_G(x)}(x).
		\]
		For $\alpha \in \Z/(c+1)\Z$, let $N_\alpha(x) \defeq \set{N_i(x) \,:\, \alpha = i \, \mathrm{mod} \, (c+1)}$. Note that $|N_\alpha(x)| \leq \lceil \deg_G(x)/(c+1) \rceil \leq c$, where we are using that $c(c+1) \geq d$. %
		Define a graph $H$ as follows. The vertex set of $H$ is \[V(H) \,\defeq\, V(G) \times ((\Z/(c+1)\Z) \sqcup [k-1]).\]
		For clarity, let $u_\alpha(x) \defeq (x, \alpha)$ and $v_i(x) \defeq (x, i)$ for all $x \in V(G)$, $\alpha \in \Z/(c+1)\Z$, and $i \in [k-1]$. Make the following pairs of vertices adjacent in $H$:
		\begin{itemize}
			\item $v_i(x)$ and $v_j(x)$ for all $x \in V(G)$ and distinct $i$, $j \in [k-1]$; %
			
			\item $v_i(x)$ and $u_\alpha(x)$ for all $x \in V(G)$, $i \in [k-1]$, and $\alpha \in \Z/(c+1)\Z$;
			
			\item $u_\alpha(x)$ and $u_\beta(y)$ for all adjacent $x$, $y \in V(G)$ with $y \in N_\alpha(x)$ and $x \in N_\beta(y)$. %
		\end{itemize}
		Let $\nu$ be the pushforward of $\mu$ under the map $V(G) \to V(H) \colon x \mapsto u_0(x)$. We will show that $H$ satisfies $\chi(H) \leq k$, $\chi_\nu(H) > k$, and $\Delta(H) \leq d-1$, which contradicts the choice of $G$ and thus completes the proof of Proposition~\ref{prop:sharp}.
		
		To see that $H$ is $k$-colorable, let $f \colon V(G) \to [k]$ be any $k$-coloring of $G$ and define $h \colon V(H) \to [k]$ as follows: For all $x \in V(G)$ and $\alpha \in \Z/(c+1)\Z$, set $h(u_\alpha(x)) \defeq f(x)$, and assign to the vertices $v_1(x)$, \ldots, $v_{k-1}(x)$ the $k-1$ colors distinct from $f(x)$. Then $h$ is a proper $k$-coloring of $H$, so $\chi(H) \leq k$. To show that $\chi_\nu(H) > k$, suppose that $h \colon V(H) \to [k]$ is a $\nu$-measurable $k$-coloring of $H$. By the definition of $\nu$, the map $f \colon V(G) \to [k] \colon x \mapsto h(u_0(x))$ is $\mu$-measurable, and we claim that it is a proper $k$-coloring of $G$, which is impossible. Consider any two adjacent vertices $x$, $y \in V(G)$. Since the vertices $v_1(x)$, \ldots, $v_{k-1}(x)$ are pairwise adjacent in $H$, there is precisely one color that is not used by $h$ on any of them, and that must be the color assigned to every vertex of the form $u_\alpha(x)$. Similarly, all the vertices of the form $u_\beta(y)$ have the same color. Let $\alpha$, $\beta \in \Z/(c+1)\Z$ be such that $y \in N_\alpha(x)$ and $x \in N_\beta(y)$. %
		Then the vertices $u_\alpha(x)$ and $u_\beta(y)$ are adjacent, so $f(x) = h(u_0(x)) = h(u_\alpha (x)) \neq h(u_\beta(y)) = h(u_0(y)) = f(y)$, as desired. It remains to verify that $\Delta(H) \leq d-1$. To this end, note that every vertex of the form $v_i(x)$ has degree
		\[
		\deg_H(v_i(x)) \,=\, (k-2) \,+\, (c+1) \,=\, d-1,
		\]
		while every vertex of the form $u_\alpha(x)$ has degree
		\[
		\deg_H(u_\alpha(x)) \,=\, (k-1) \,+\, |N_\alpha(x)| \,\leq\, %
		k-1+ c \,=\, d-1. \qedhere
		\]
	\end{scproof}
	
	There is a small gap between the bounds in Theorem~\ref{theo:Delta} and Proposition~\ref{prop:sharp}. We leave closing this gap as an open problem:
	
	\begin{prob}
		For each $\Delta \in \N$, determine precisely the largest value of $c \in \Z$ for which the conclusion of Theorem~\ref{theo:Delta} holds.
	\end{prob}

	Theorem~\ref{theo:Delta} implies Theorem~\ref{theo:meas_Brooks}, i.e., the measurable version of Brooks's theorem, for all large enough $\Delta$. One can actually use distributed algorithms to deduce Theorem~\ref{theo:meas_Brooks} for all $\Delta \geq 3$. Indeed, Ghaffari, Hirvonen, Kuhn, and Maus \cite{GHKM} developed a randomized \LOCAL algorithm that, given an $n$-vertex graph $G$ of maximum degree $\Delta \geq 3$ and without a complete subgraph on $\Delta + 1$ vertices, finds a proper $\Delta$-coloring of $G$ in $O((\log \log n)^2) = o(\log n)$ rounds \ep{here, as in \eqref{eq:Delta_bound}, the implicit constants in the asymptotic notation may depend on $\Delta$}. Combined with Theorem~\ref{theo:dist_meas}, this yields Theorem~\ref{theo:meas_Brooks}.
	
	\subsection{Graphs without short cycles}\label{subsec:tr_fr}
	
	The distinctions between the three regularity notions---``Borel,'' ``measurable,'' and ``Baire\-/measur\-able''---are clearly demonstrated by colorings of acyclic graphs. Recall that, by Marks's Theorem~\ref{theo:Marks}, the Borel chromatic number of an acyclic Borel graph of maximum degree $\Delta \in \N$ can be as large as $\Delta + 1$. In contrast to this, the Baire-measurable chromatic number of a locally finite acyclic graph is always at most $3$, which is a consequence of the following general result of Conley and Miller:
	
	\begin{theo}[{Conley--Miller \cite{ConleyMiller}}]\label{theo:ConleyMiller}
		Let $G$ be a locally finite Borel graph and let $\tau$ be a compatible Polish topology on $V(G)$. If $\chi(G)$ is finite, then $\chi_\tau(G) \leq 2\chi(G) - 1$. %
	\end{theo}
	
	Since $\chi(G) \leq 2$ for an acyclic graph $G$, Theorem~\ref{theo:ConleyMiller} implies that the Baire-measurable chromatic number of an acyclic locally finite Borel graph is indeed at most $2\cdot 2 - 1 = 3$. It is also not hard to see that this upper bound is best possible; see, e.g., \cite[\S6]{CMTD}.
	
	Now we turn to measurable colorings. Lyons and Nazarov \cite{LyonsNazarov} \ep{see also \cite[Theorem 5.46]{KechrisMarks}} constructed acyclic Borel graphs with maximum degree $\Delta \in \N$ and measurable chromatic number at least $(1/2+o(1))\Delta/\log \Delta$. This shows that acyclic Borel graphs of bounded degree can have arbitrarily large measurable chromatic numbers, contrary to the situation with Baire-measurable colorings. Although the present author showed \cite[Corollary 1.2]{MLLL} that $\Theta(\Delta/\log \Delta)$ is the correct order of magnitude in the Lyons--Nazarov examples, the best heretofore known general upper bound on measurable chromatic numbers of acyclic graphs in terms of their maximum degree was $\Delta$, which is a consequence of the measurable Brooks's Theorem~\ref{theo:meas_Brooks} of Conley--Marks--Tucker-Drob. Here we improve this to $(1+o(1))\Delta/\log \Delta$, which falls within a factor of $2$ of best possible. Furthermore, we do not even require acyclicity---it suffices to only forbid cycles of length at most $4$.

	\begin{theo}\label{theo:triangle_free}
		For every $\epsilon > 0$, there is $\Delta_0 \in \N$ with the following property. Let $G$ be a Borel graph of finite maximum degree $\Delta \geq \Delta_0$ and let $\mu$ be a probability Borel measure on $V(G)$.
		\begin{enumerate}[label=\ep{\normalfont\roman*}]
			\item\label{item:no_4_c} If $G$ contains no cycles of length at most $4$, then $\chi_\mu(G) \leq (1+\epsilon)\Delta/\log \Delta$.
			
			\item\label{item:tr_free} If $G$ contains no cycles of length $3$, then $\chi_\mu(G) \leq (4+\epsilon)\Delta/\log \Delta$.
		\end{enumerate}
	\end{theo}
	\begin{scproof}
		To derive this theorem, we invoke the \LOCAL algorithms for graph coloring developed by Chung, Pettie, and Su in \cite{CPS}. We shall first prove \ref{item:tr_free}. For each $k \in \N$, let $\Pi(k)$ denote the local coloring problem that encodes proper $k$-coloring of graphs \ep{see Example~\ref{exmp:prop_col}}. Fix $\epsilon >0$ and $\Delta \in \N$ and set $k \defeq \lfloor (4+\epsilon) \Delta/\log \Delta \rfloor$. Let $\class{G}(\Delta)$ be the set of all isomorphism classes of finite graphs of maximum degree at most $\Delta$ with no $3$-cycles. It follows from \cite[Theorem 9]{CPS} that there exist positive real numbers $C$ and $\delta$, depending only on $\epsilon$, such that, for all large enough $\Delta$ and $n$,
		\begin{equation}\label{eq:tr_fr_rand_bound}
			\Rand_{\Pi(k), \class{G}(\Delta)}(n) \,\leq\, C \Delta^{-\delta} \log n.
		\end{equation}
		Even though this upper bound is weaker than $o(\log n)$, it still suffices for an application of Theorem~\ref{theo:dist_meas}. Indeed, let $G$ be a Borel graph of maximum degree $\Delta$ without cycles of length $3$ and set $R(n) \defeq \Rand_{\Pi(k), \class{G}(\Delta)}(n) + 1$. If $\Delta$ is large enough and $n$ is much larger than $\Delta$, then for each vertex $x \in V(G)$,
		\[
		|B_G(x, 2R(n))| \,\leq\, 1 + \Delta^{2C\Delta^{-\delta} \log n + 2} \, =\, 1 + \exp(2C\Delta^{-\delta} \log \Delta \log n + 2\log \Delta) \,<\, n^{1/8}/4.  
		\]
		Similarly, $|B_G(x, R(n))| < n$, and, since $\class{G}(\Delta)$ is closed under %
		adding isolated vertices, $G$ is $(R(n), n)$-locally in $\class{G}(\Delta)$.
		Hence, Theorem~\ref{theo:dist_meas} allows us to conclude that $G$ is measurably $k$-colorable, which proves \ref{item:tr_free}. The proof of \ref{item:no_4_c} is virtually the same, except that \cite[Theorem 9]{CPS} is replaced by the randomized \LOCAL algorithm of Chung, Pettie, and Su for $((1+o(1))\Delta / \log \Delta)$-coloring graphs without cycles of length at most $4$ \ep{see the remark in \cite{CPS} after \cite[Theorem 9]{CPS}}.
	\end{scproof}
	
	Again, if $G$ is of subexponential growth, then, due to Theorem~\ref{theo:dist_subexp}, the upper bounds on $\chi_\mu(G)$ given by Theorem~\ref{theo:triangle_free} also hold for $\chi_{\mathrm{B}}(G)$. %
	
	We remark that, although this was not needed in the proof of Theorem~\ref{theo:triangle_free}, the bound 
	\eqref{eq:tr_fr_rand_bound} can actually be improved to $\Rand_{\Pi(k), \class{G}(\Delta)}(n) = o(\log n)$. Indeed, the algorithm developed by Chung, Pettie, and Su in order to prove \eqref{eq:tr_fr_rand_bound} invokes as a subroutine a certain distributed version of the LLL. Fischer and Ghaffari \cite{FG} later designed a more efficient distributed version of the LLL, and using their result reduces the complexity of the Chung--Pettie--Su algorithm to $o(\log n)$.
	
	The upper bound $\chi(G) \leq (1+o(1))\Delta/ \log \Delta$ for finite graphs $G$ without cycles of length at most $4$ is due to Kim \cite{Kim}. The existence of a constant $C > 0$ such that 
	$\chi(G) \leq (C + o(1))\Delta/\log \Delta$ for finite graphs $G$ without $3$-cycles was first established by Johansson \cite{Johansson} \ep{Johansson's original paper is hard to access, but a detailed presentation of his argument can be found in \cite[\S{}13]{MR}}. Johansson's original proof gave the value $C=9$. Pettie and Su \cite{PS15} improved this to $C = 4$ \ep{which is the value appearing in Theorem~\ref{theo:triangle_free}}. Recently, Molloy \cite{Molloy} further reduced the constant to $C = 1$. We leave the question of whether Molloy's result also holds for measurable colorings as an open problem:
	
	\begin{prob}\label{prob:tr_fr}
		Is it true that for every $\epsilon > 0$, there is $\Delta_0 \in \N$ with the following property? Let $G$ be a Borel graph of finite maximum degree $\Delta \geq \Delta_0$ and let $\mu$ be a probability Borel measure on $V(G)$. If $G$ contains no cycles of length $3$, then $\chi_\mu(G) \leq (1+\epsilon) \Delta/\log \Delta$.
	\end{prob}
	
	One way to solve Problem~\ref{prob:tr_fr} would be to develop a sublogarithmic randomized \LOCAL algorithm for $((1+o(1))\Delta/\log \Delta)$-coloring graphs without $3$-cycles.

	\subsection{A result on list-coloring}
	
	In this section we apply Theorem~\ref{theo:dist_meas} to obtain a useful result concerning list-colorings of Borel graphs. For an introduction to the theory of list-coloring, see \cite[\S5.4]{Die}. A \emphd{list assignment} for a graph $G$ is a function $L \colon V(G) \to \fins{\N}$. For each vertex $x \in V(G)$, the set $L(x)$ is called the \emphd{list} of $x$, and the elements of $L(x)$ are the colors \emphd{available} to $x$. An \emphd{$L$-coloring} of $G$ is a map $f \colon V(G) \to \N$ such that $f(x) \in L(x)$ for all $x \in V(G)$. An $L$-coloring $f$ is \emphd{proper} if $f(x) \neq f(y)$ whenever $x$ and $y$ are adjacent in $G$. We say that $L$ is an \emphd{$(\ell, d)$-list assignment} if it has the following two properties:
	\begin{itemize}
		\item for all $x \in V(G)$, $|L(x)| \geq \ell$; and
		\item for all $x \in V(G)$ and $\alpha \in L(x)$, $|\set{y \in N_G(x) \,:\, \alpha \in L(y)}| \leq d$.
	\end{itemize}
	The following is a result of Reed and Sudakov:
	
	\begin{theo}[{Reed--Sudakov \cite{ReedSud}}]\label{theo:ReedSud}
		For every $\epsilon > 0$, there is $d_0 \in \N$ with the following property. Let $G$ be a graph and let $L$ be an $(\ell, d)$-list assignment for $G$, where $d \geq d_0$ and $\ell \geq (1+\epsilon)d$. Then $G$ admits a proper $L$-coloring.
	\end{theo}
	
	A version of Theorem~\ref{theo:ReedSud} with the bound $\ell \geq (1+\epsilon)d$ replaced by $\ell \geq 2e d$ was first proved by Reed \cite{Reed}; this was then improved by Haxell \cite{Haxell} to $\ell \geq 2d$. While it is weaker than Theorem~\ref{theo:ReedSud} for large $d$, Haxell's result holds for all positive integers $d$ and not only for sufficiently large ones. Reed \cite{Reed} conjectured that in fact $\ell \geq d+1$ should suffice, but this conjecture was refuted by Bohman and Holzman \cite{BohHol}.
	
	Theorem~\ref{theo:ReedSud} \ep{or its weaker versions mentioned in the previous paragraph} is a somewhat technical but rather useful fact that plays a crucial role in the proofs of many graph coloring results \ep{see \cite{MR} for a number of examples}. Here we establish a version of Theorem~\ref{theo:ReedSud} for measurable and Baire-measurable colorings: %
	
	\begin{theo}\label{theo:meas_RS}
		For every $\epsilon > 0$, there is $d_0 \in \N$ with the following property. Let $G$ be a Borel graph and let $L$ be a Borel $(\ell, d)$-list assignment for $G$, where $d \geq d_0$ and $\ell \geq (1+\epsilon)d$. Let $\mu$ \ep{resp.~$\tau$} be a probability Borel measure \ep{resp.~a compatible Polish topology} on $V(G)$. Then $G$ admits a $\mu$-measurable \ep{resp.~$\tau$-Baire-measurable} proper $L$-coloring.
	\end{theo}
	\begin{scproof}
		We shall use another \LOCAL algorithm due to Chung, Pettie, and Su \cite{CPS}. As explained in Example \ref{exmp:list-col}, we can interpret pairs $(G, L)$, where $G$ is a graph and $L$ a list assignment for $G$, as structured graphs in the sense of Definition~\ref{defn:str_graph}. Furthermore, proper list-coloring can naturally be encoded as a local coloring problem $\Pi$. Explicitly, let $\mathcal{P} \colon \class{FSG}_\bullet \to \set{0,1}$ be the \LOCAL algorithm defined as follows: Given \ep{the isomorphism type of} a finite rooted structured graph of the form $(G, L, f, x)$ with $L$ a list assignment for $G$ and $f \colon V(G) \to \N$, set $\mathcal{P}([G, L, f, x]) \defeq 1$ if and only if $f$ is a proper $L$-coloring of $G$; in all other cases set $\mathcal{P}$ to $0$. Now if we let $\Pi \defeq (1, \mathcal{P})$, then a $\Pi$-coloring of $(G,L)$ is precisely the same as a proper $L$-coloring of $G$.
		
		Let $\class{G}(\ell, d) \subset \class{FSG}$ denote the set of all isomorphism types of finite structured graphs corresponding to pairs of the form $(G, L)$, where $G$ is a finite graph and $L$ is an $(\ell, d)$-list assignment for $G$. The result of Chung, Pettie, and Su presented in \cite[\S4.4]{CPS} implies that for every $\epsilon > 0$, there exist positive real numbers $C$ and $\delta$ such that, for all large enough $d$ and $n$,
		\begin{equation}\label{eq:list_bound}
			\Rand_{\Pi, \class{G}((1+\epsilon)d, d)}(n) \,\leq\, C d^{-\delta} \log n.
		\end{equation}
		As in the proof of Theorem~\ref{theo:triangle_free}, this upper bound is weaker than $o(\log n)$ but sufficient for an application of Theorem~\ref{theo:dist_meas}. Indeed, let $G$ be a Borel graph and let $L$ be a Borel $((1+\epsilon)d, d)$-list assignment for $G$, where $d$ is a large positive integer. Without loss of generality, we may assume that $|L(x)| = \ell$ for all $x \in V(G)$, and, by removing from $G$ all the edges $xy$ such that $L(x) \cap L(y) = \0$, we may arrange that $\Delta(G) \leq (1+\epsilon) d^2 \leq d^3$. Set $R(n) \defeq \Rand_{\Pi, \class{G}((1+\epsilon)d, d)}(n) + 1$. For all $x \in V(G)$,
		\[
		|B_G(x, 2R(n))| \,\leq\, 1 + d^{6C d^{-\delta} \log n + 6} \,=\, 1 + \exp(6Cd^{-\delta} \log d \log n + 6 \log d) \,<\, n^{1/8}/4,
		\]
		where the last inequality holds whenever $d$ and $n$ are large enough. Similarly, $|B_G(x, R(n))| < n$, and hence the pair $(G,L)$ is $(R(n),n)$-locally in $\class{G}((1+\epsilon)d, d)$. Therefore, by Theorem~\ref{theo:dist_meas}, $G$ admits a measurable/Baire-measurable proper $L$-coloring, as desired.
	\end{scproof}
	
	Again, for graphs of subexponential growth, we can use Theorem~\ref{theo:dist_subexp} to upgrade the conclusion of Theorem~\ref{theo:meas_RS} to a Borel proper $L$-coloring. Like \eqref{eq:tr_fr_rand_bound}, the bound \eqref{eq:list_bound} can be improved to $\Rand_{\Pi, \class{G}((1+\epsilon)d, d)}(n) = o(\log n)$ by using the Fischer--Ghaffari distributed LLL \cite{FG}.
	
	\subsection{Sparse graphs}
	
	Among the oldest topics in graph theory is studying colorings of planar graphs. For the purposes of this paper, we say that an infinite graph is planar if all its finite subgraphs are planar. One of the most celebrated results in graph theory is the \emph{Four Color Theorem} of Appel and Haken \cite[Theorem 5.1.1]{Die} that asserts that $\chi(G) \leq 4$ for all planar graphs $G$. Furthermore, if $G$ is a planar graph without $3$-cycles, then $\chi(G) \leq 3$---this is a theorem of Gr\"otzsch \cite[Theorem 5.1.3]{Die}. On the other hand, the Lyons--Nazarov examples \cite{LyonsNazarov} mentioned in \S\ref{subsec:tr_fr} \ep{see also \cite[Theorem 5.46]{KechrisMarks}} show that bounded degree acyclic \ep{hence planar} Borel graphs can have arbitrarily large measurable chromatic numbers. Nevertheless, we show that under the additional assumption of subexponential growth, Borel chromatic numbers of planar Borel graphs are bounded by $5$, and can be further lowered for graphs without short cycles:
	
	\begin{theo}\label{theo:planar}
		Let $G$ be a planar Borel graph of subexponential growth and define
		\[
		k \,\defeq\, \begin{cases}
			3 & \text{if $G$ contains no cycles of length at most $4$};\\
			4 & \text{if $G$ contains a $4$-cycle but no $3$-cycles};\\
			5 & \text{otherwise}.
		\end{cases}
		\]
		Then $\chi_{\mathrm{B}}(G) \leq k$. Furthermore, if $G$ is also a topological graph \ep{in the sense of Definition~\ref{defn:top_graph}}, then $G$ admits a continuous proper $k$-coloring.
	\end{theo}
	\begin{scproof}
		This is a consequence of the recent work of Postle \cite{Postle_planar}.
		Let $\class{P}$ denote the set of all isomorphism classes of finite planar graphs and let $\class{P}(g) \subseteq \class{P}$ be the set of all isomorphism classes of finite planar graphs without cycles of length strictly less than $g$. Letting $\Pi(k)$ denote the local coloring problem corresponding to proper $k$-coloring \ep{see Example~\ref{exmp:prop_col}}, Postle \cite[Theorem 1.3]{Postle_planar} established the following bounds:
		\[
		\max \set{\Det_{\Pi(5), \class{P}}(n),\, \Det_{\Pi(4), \class{P}(4)}(n), \, \Det_{\Pi(3), \class{P}(5)}(n)} \,=\, O(\log n).
		\]
		Combined with Theorems~\ref{theo:dist_Borel} and \ref{theo:dist_cont}, this yields the desired results. \ep{Theorems~\ref{theo:dist_Borel} and \ref{theo:dist_cont} may be applied as graphs of subexponential growth have finite maximum degree.}
	\end{scproof}
	
	An important feature of planar graphs is their sparsity: by Euler's formula, a planar graph $G$ with $n \geq 3$ vertices can have at most $3n-6$ edges. A convenient measure of sparsity for an arbitrary graph is its \emph{arboricity}, defined as follows. Recall that $|G|$ and $\|G\|$ are the cardinalities of the vertex and the edge sets of $G$, respectively. The \emphd{arboricity} of a graph $G$ with $|G| \geq 2$ is the quantity
	\[
	a(G) \,\defeq\, \sup_H \left\lceil \frac{\|H\|}{|H| - 1} \right\rceil, 
	\]
	where the supremum is taken over all the finite subgraphs $H$ of $G$ with $|H| \geq 2$. For graphs $G$ with $|G| \leq 1$, $a(G) \defeq 0$ by definition. A theorem of Nash-Williams \cite[Theorem 2.4.4]{Die} asserts that if $G$ is a finite graph, then $a(G)$ is equal to the smallest $k \in \N$ such that $G$ has $k$ acyclic subgraphs $F_1$, \ldots, $F_k$ with $E(G) = E(F_1) \cup \ldots \cup E(F_k)$. Since acyclic graphs are also called ``forests,'' this explains the term ``arboricity.''
	
	Many classes of graphs have bounded arboricity. For instance, all planar graphs have arboricity at most $3$. It is not hard to see that graphs $G$ of finite arboricity satisfy $\chi(G) \leq 2 a(G)$ \ep{and this bound is, in general, best possible}. On the other hand, the measurable chromatic number of a bounded degree acyclic \ep{i.e., arboricity $1$} Borel graph can be arbitrarily large. As in the case of planar graphs, the situation improves under the subexponential growth assumption:
	
	\begin{theo}\label{theo:arb_subexp}
		If $G$ is a Borel graph of subexponential growth, then $\chi_{\mathrm{B}}(G) \leq 2a(G) + 1$. Furthermore, if $G$ is also a topological graph, then $G$ admits a continuous proper $(2a(G) +1)$-coloring.
	\end{theo}
	\begin{scproof}
		Let $\class{G}(a)$ denote the set of all isomorphism classes of finite graphs of arboricity at most $a$, and let $\Pi(k)$ be the local coloring problem corresponding to proper $k$-coloring. Barenboim and Elkin \cite[\S4]{BarEl_arb} proved that for all $a \in \N$,
		\[
		\Det_{\Pi(2a+1), \class{G}(a)} \,=\, O(\log n),
		\]
		where the implicit constants in the asymptotic notation may depend on $a$. \ep{The statement in \cite[\S4]{BarEl_arb} is more general and involves an additional positive parameter $\epsilon$; the bound that we need is obtained by setting $\epsilon = 1/(a+1)$.} It remains to apply Theorems~\ref{theo:dist_Borel} and \ref{theo:dist_cont}. %
	\end{scproof}
	
	The bound $\chi_{\mathrm{B}}(G) \leq 2a(G) + 1$ for Borel graphs $G$ of subexponential growth is, in general, sharp, since there exist acyclic $2$-regular Borel graphs $G$ with $\chi_{\mathrm{B}}(G) = 3$, and such graphs are of subexponential \ep{in fact, linear} growth; see, e.g., \cite[\S6]{CMTD}. Nevertheless, we conjecture that acyclic graphs are an exception and the bound $\chi_{\mathrm{B}}(G) \leq 2a(G)$ should hold whenever $a(G) \geq 2$:
	
	\begin{conj}\label{conj:arb}
		If $G$ is a Borel graph of subexponential growth, then $\chi_{\mathrm{B}}(G) \leq \max \set{2a(G), 3}$.
	\end{conj}
	
	As evidence for Conjecture~\ref{conj:arb}, we show that it holds for graphs whose order of growth is somewhat lower than just subexponential:
	
	\begin{theo}
		Let $G$ be a Borel graph and suppose that for each $\epsilon > 0$, there is $R_0 \in \N$ such that if $x \in V(G)$ and $R \geq R_0$, then $|B_G(x, R)| < \exp(\epsilon R^{1/3})$. Then $\chi_{\mathrm{B}}(G) \leq \max \set{2a(G), 3}$. Furthermore, if $G$ is also a topological graph, then $G$ admits a continuous proper $\max \set{2a(G), 3}$-coloring.
	\end{theo}
	\begin{scproof}
		As in the proof of Theorem~\ref{theo:arb_subexp}, let $\class{G}(a)$ be the set of all isomorphism classes of finite graphs of arboricity at most $a$, and let $\Pi(k)$ be the local coloring problem encoding proper $k$-coloring. Aboulker, Bonamy, Bousquet, and Esperet \cite[Corollary 1.4]{ABBE} showed that for $a \geq 2$,
		\[
		\Det_{\Pi(2a), \class{G}(a)} = O((\log n)^3)
		\]
		where the implicit constants in the asymptotic notation again depend on $a$. Now a straightforward computation shows that Theorems~\ref{theo:dist_Borel} and \ref{theo:dist_cont} yield the desired results.
	\end{scproof}
	
	Another result about graphs of subexponential growth that should be mentioned here is a theorem of Gao and Jackson \cite[Theorem 4.2]{GaoJack} that asserts that for every $d \geq 2$, the Schreier graph of the free part of the Bernoulli shift action $\Z^d \acts \set{0,1}^{\Z^d}$ admits a continuous proper $4$-coloring. This fact can be alternatively derived by using Theorem~\ref{theo:dist_cont} in combination with the deterministic \LOCAL algorithm for $4$-coloring grid graphs designed by Brandt et al.~\cite[Theorem 4]{grids}.
	
	\subsection{A pointwise version of the Ab\'ert--Weiss theorem}\label{subsec:AW}
	
	So far we have described a number of results in descriptive combinatorics that can be obtained using distributed algorithms. Now we present a direct application of the \hyperref[theo:meas_SLLL]{Measurable Symmetric LLL} to a question in ergodic theory.  
	
	Throughout \S\ref{subsec:AW}, $\G$ shall denote a countably infinite group. We are interested in \emphd{probability measure\-/preserving \ep{\pmp} actions} of $\G$, i.e., actions of the form $\alpha \colon \G \acts (X, \mu)$, where $(X, \mu)$ is a standard probability space and the measure $\mu$ is $\alpha$-invariant. More generally, we consider \emphd{Borel actions} $\alpha \colon \G \acts X$, i.e., actions of $\G$ on a standard Borel space $X$ by Borel automorphisms. An action $\alpha \colon \G \acts X$ is \emphd{free} if the $\alpha$-stabilizer of every point $x \in X$ is trivial.
	
	An important example of a \pmp action is the {Bernoulli shift action} \[\sigma \colon \G\acts ([0,1]^\G, \lambda^\G),\] where $([0,1], \lambda)$ is the unit interval equipped with the Lebesgue probability measure $\lambda$ \ep{owing to the measure isomorphism theorem \cite[Theorem~17.41]{KechrisDST}, any other atomless standard probability space could be used instead}. For brevity, we write \[\Omega \defeq [0,1]^\G \qquad \text{and} \qquad \bm{\lambda} \defeq \lambda^\G\] \ep{this notation will only be used in \S\ref{subsec:AW}}. Our starting point is a result of Ab\'ert and Weiss:
	
	\begin{theo}[{Ab\'ert--Weiss \cite{AW}}]\label{theo:AW}
		Fix the following data:
		\begin{itemize}
			\item a partition $\Omega = A_1 \sqcup \ldots \sqcup A_k$ of $\Omega$ into finitely many Borel pieces;
			\item a finite set $F \in \fins{\G}$; and
			\item $\epsilon > 0$.
		\end{itemize}
		For every free \pmp action $\alpha \colon \G \acts (X, \mu)$, there is a Borel partition $X = B_1 \sqcup \ldots \sqcup B_k$ such that %
		\[
		\mu(B_i \cap \gamma \cdot B_j) \,\approx_\epsilon \, \bm{\lambda}(A_i \cap \gamma \cdot A_j),
		\]
		for all  $1 \leq i$, $j \leq k$ and all $\gamma \in F$.
	\end{theo}
	
	Here and in what follows, we write $a \approx_\epsilon b$ to mean $|a-b| < \epsilon$. Theorem~\ref{theo:AW} can be stated briefly as ``The shift action $\sigma \colon \G\acts (\Omega, \bm{\lambda})$ is \emph{weakly contained} in every free \pmp action $\alpha \colon \G\acts (X, \mu)$.'' The relation of weak containment was introduced by Kechris in \cite[\S10(C)]{K_book}. For more details and further background on this topic, see the survey \cite{BK} by Burton and Kechris. 
	
	Here we strengthen Theorem~\ref{theo:AW} by replacing the quantities $\mu(B_i \cap \gamma \cdot B_j)$ with certain pointwise averages almost everywhere:
	
	\begin{theo}\label{theo:point_AW}
		Fix the following data:
		\begin{itemize}
			\item a partition $\Omega = A_1 \sqcup \ldots \sqcup A_k$ of $\Omega$ into finitely many Borel pieces;
			\item a finite set $F \in \fins{\G}$; and
			\item $\epsilon > 0$.
		\end{itemize}
		Then there is $n_0 \in \N^+$ with the following property.  %
		Fix a finite set $D \in \fins{\G}$ of size $|D| \geq n_0$ and let $\alpha \colon \G \acts X$ be a free Borel action of $\G$. For every probability Borel measure $\mu$ on $X$, there exists a Borel partition $X = B_1 \sqcup \ldots \sqcup B_k$ such that
		\begin{equation}\label{eq:point_approx}
			\frac{|\set{\delta \in D \,:\, \delta \cdot x \in B_i \cap \gamma \cdot B_j}|}{|D|} \,\approx_\epsilon \, \bm{\lambda}(A_i \cap \gamma \cdot A_j),
		\end{equation}
		for all  $1 \leq i$, $j \leq k$, all $\gamma \in F$, and $\mu$-almost all $x \in X$.
	\end{theo}
	
	Note that in Theorem~\ref{theo:point_AW}, the measure $\mu$ is not required to be $\alpha$-invariant. If $\mu$ is $\alpha$-invariant, then the partition $X = B_1 \sqcup \ldots \sqcup B_k$ given by Theorem~\ref{theo:point_AW} also witnesses the conclusion of Theorem~\ref{theo:AW}, as for an $\alpha$-invariant measure $\mu$, we have
	\[
	\mu(Y) \,=\, \int_X \frac{|\set{\delta \in D \,:\, \delta \cdot x \in Y}|}{|D|} \,\Diff \mu(x),
	\]
	for all Borel $Y \subseteq X$ and any nonempty finite set $D \in \fins{\G}$. This shows that Theorem~\ref{theo:point_AW} is indeed a strengthening of Theorem~\ref{theo:AW}. Statements in the spirit of Theorem~\ref{theo:point_AW} were first considered by the present author in \cite{Ber_AW}. There, a weaker version of Theorem~\ref{theo:point_AW} was established, with \eqref{eq:point_approx} satisfied not for $\mu$-almost all $x \in X$, but only on a set of $x \in X$ of measure at least $1 - \delta$, for any given $\delta > 0$ \cite[Theorem 2.11]{Ber_AW}. The question of whether Theorem~\ref{theo:point_AW} is true was left there as an open problem \cite[Problem 8.2]{Ber_AW}.
	
	\begin{scproof}[ of Theorem~\ref{theo:point_AW}]
		Only minimal modifications to the proof of \cite[Theorem 2.11]{Ber_AW} are needed to obtain Theorem~\ref{theo:point_AW}. Specifically, in \cite[\S\S4.B and 7.A]{Ber_AW} the construction of a partition $X = B_1 \sqcup \ldots \sqcup B_k$ with the desired property is reduced to finding a measurable solution to a certain Borel CSP $\mathscr{B}$ on $X$. This CSP depends on the original partition $\Omega = A_1 \sqcup \ldots \sqcup A_k$, the finite set $F$, the parameter $\epsilon$, and the choice of the averaging set $D$. It is immediate from the construction of $\mathscr{B}$ that $\sup \set{|\dom(B)| \,:\, B \in \mathscr{B}} < \infty$. Furthermore, the calculations given in the proof of \cite[Lemma 7.2]{Ber_AW} show that there exist positive reals $a$, $b$, and $c$ depending on the partition $\Omega = A_1 \sqcup \ldots \sqcup A_k$, the finite set $F$, and the parameter $\epsilon$ but not on $D$ such that
		\begin{equation}\label{eq:point_AW_bounds}
			\pr(\mathscr{B}) \,\leq\, a \exp (-b |D|) \qquad \text{and} \qquad \de(\mathscr{B}) \,\leq\, c |D|^2 -1.
		\end{equation}
		At this point the proof of \cite[Theorem 2.11]{Ber_AW} invokes the approximate Symmetric LLL from \cite{MLLL} \ep{Theorem 6.5 in \cite{Ber_AW}}, which yields a measurable function that satisfies the constraints of the CSP $\mathscr{B}$ away from a set of measure less than $\delta$, for given $\delta > 0$. To obtain Theorem~\ref{theo:point_AW}, we instead use Theorem~\ref{theo:meas_SLLL}. By \eqref{eq:point_AW_bounds},
		\[
		\pr(\mathscr{B}) \cdot (\de(\mathscr{B}) + 1)^8 \,\leq\, a \exp(-b|D|) \cdot c^8 |D|^{16}.	
		\]
		The latter quantity approaches $0$ as $|D| \to \infty$; in particular, it is less than $2^{-15}$ whenever $|D|$ is large enough. Hence, by Theorem~\ref{theo:meas_SLLL}, assuming $|D|$ is sufficiently large, the CSP $\mathscr{B}$ has a measurable solution, and we are done.
	\end{scproof}
	
	\section{Using distributed algorithms}\label{sec:dist_to_desc}
	
	\subsection{Proof of Theorem~\ref{theo:dist_Borel}}\label{subsec:dist_Borel}
	
	We start with the following simple observation, which will be used repeatedly throughout \S\S\ref{sec:dist_to_desc} and \ref{sec:LLL_proof}:
	
	\begin{lemma}\label{lemma:BorelBall}
		If $\bm{G}$ is a locally finite Borel structured graph, then the function \[V(\bm{G}) \times \N \to \class{FSG}_\bullet \colon (x, R) \mapsto [B_{\bm{G}}(x,R), x]\] is Borel \ep{here $\class{FSG}_\bullet$ is viewed as a discrete countable space}.
	\end{lemma}
	\begin{scproof}
		We need to argue that for every finite rooted structured graph $(\bm{H}, v)$, the set
		\begin{equation}\label{eq:xR}
			\set{(x, R) \in V(\bm{G}) \times \N \,:\, (B_{\bm{G}}(x,R), x) \cong (\bm{H}, v)}
		\end{equation}
		is Borel. Since $\bm{G}$ is locally finite, by the Feldman--Moore theorem \cite[Theorem 1.3]{KechrisMiller}, there exist Borel involutions $\gamma_i \colon V(\bm{G}) \to V(\bm{G})$, $i \in \N$, such that $\dist_{\bm{G}}(x, y) < \infty$ if and only if $y = \gamma_i(x)$ for some $i \in \N$. Without loss of generality, we may assume that $\gamma_0$ is the identity map on $V(\bm{G})$. For all $x$, $y \in V(\bm{G})$ such that $\dist_{\bm{G}}(x,y) < \infty$, let $\mathsf{index}(x, y)$ denote the minimum $i \in \N$ with $\gamma_{i}(x) = y$ (in particular, $\mathsf{index}(x,x) = 0$). Now, for any $x \in V(\bm{G})$ and $R \in \N$, we have that $(B_{\bm{G}}(x,R), x) \cong (\bm{H}, v)$ if and only if there exists a partial function $\phi \in \finf{\N}{V(\bm{H})}$ such that:
		\begin{itemize}
			\item $\phi$ is a bijection between $\dom(\phi) \in \fins{\N}$ and $V(\bm{H})$;
			\item for all $i \in \N$, we have $i \in \dom(\phi)$ if and only if $\dist_{\bm{G}}(x, \gamma_i(x)) \leq R$ and $i = \mathsf{index}(x,\gamma_i(x))$;
			\item the mapping $\gamma_i(x) \mapsto \phi(i)$ is an isomorphism of the structured graphs $B_{\bm{G}}(x,R)$ and $\bm{H}$;
			\item $\phi(0) = v$.
		\end{itemize}
		Since the set $\finf{\N}{V(\bm{H})}$ is countable, it follows that the set \eqref{eq:xR} is a countable union of Borel sets (one for each $\phi \in \finf{\N}{V(\bm{H})}$), and hence it is itself Borel.
	\end{scproof}
	
	A useful consequence of Lemma~\ref{lemma:BorelBall} is that if $\mathcal{A}$ is a \LOCAL algorithm and $\bm{G}$ is a locally finite Borel structured graph, then for each $T \in \N$, the map $\mathcal{A}(\bm{G}, T) \colon V(\bm{G}) \to \N$ is Borel. We will make extensive use of this observation throughout the remainder of the paper.
	
	After these preliminary remarks, we are ready to prove Theorem~\ref{theo:dist_Borel}.
	
	\begin{scproof}[ of Theorem~\ref{theo:dist_Borel}]
		Our proof of Theorem~\ref{theo:dist_Borel} is similar to the argument used by Chang, Kopelowitz, and Pettie to prove that no local coloring problem has deterministic \LOCAL complexity in the range $\omega(\log^\ast n)$ and $o(\log n)$ \cite[Corollary 3]{CKP}. Let $\Pi = (t, \mathcal{P})$ be a local coloring problem and let $\class{G} \subseteq \class{FSG}$. Fix $n \in \N^+$ such that $T \defeq \Det_{\Pi, \class{G}}(n)$ is finite and let $\mathcal{A}$ be a \LOCAL algorithm witnessing the bound $\Det_{\Pi, \class{G}}(n) \leq T$. Set $R \defeq T + t$. Now let $\bm{G}$ be a Borel structured graph that is $(R,n)$-locally in $\class{G}$ and such that $|B_{\bm{G}}(x, 2R)| \leq n$ for all $x \in V(\bm{G})$. Our goal is to show that $\bm{G}$ has a Borel $\Pi$-coloring. To this end, let $G'$ be the graph with $V(G') \defeq V(\bm{G})$ in which two distinct vertices $x$, $y$ are adjacent if and only if $\dist_{\bm{G}}(x,y) \leq 2R$. The graph $G'$ is Borel and satisfies $\Delta(G') = \sup \set{|B_{\bm{G}}(x, 2R)| - 1 \,:\, x \in V(G)} \leq n - 1$ \ep{we are subtracting $1$ since a vertex is never adjacent to itself}. Hence, by the Kechris--Solecki--Todorcevic Theorem~\ref{theo:KST}%
		, $G'$ has a Borel proper coloring $c \colon V(\bm{G}) \to [n]$. Define a function $f \colon V(\bm{G}) \to \N$ by $f \defeq \mathcal{A}(\bm{G}_c, T)$. The function $f$ is Borel, and we claim that it is a $\Pi$-coloring of $\bm{G}$. %
		In other words, we claim that %
		\[
		\mathcal{P}(\bm{G}_f, t)(x) \,=\, 1 \quad \text{for all } x \in V(\bm{G}).
		\]
		Fix any $x \in V(\bm{G})$. Since $\bm{G}$ is $(R,n)$-locally in $\class{G}$, there exist an $n$-vertex structured graph $\bm{H}$ with $[\bm{H}] \in \class{G}$ and a vertex $y \in V(\bm{H})$ such that $[B_{\bm{G}}(x, R), x] = [B_{\bm{H}}(y, R), y]$. Let $\phi$ be an isomorphism between $B_{\bm{H}}(y, R)$ and $B_{\bm{G}}(x, R)$ sending $y$ to $x$. Notice that the vertices of $B_{\bm{G}}(x, R)$ are pairwise adjacent in $G'$, so they are assigned distinct colors by $c$. Hence, we can extend the function $c \circ \phi$ to a bijection $\mathsf{id} \colon V(\bm{H}) \to [n]$. Since $\mathcal{A}$ is a \LOCAL algorithm witnessing the bound $\Det_{\Pi, \class{G}}(n) \leq T$, the function $g \defeq \mathcal{A}(\bm{H}_{\mathsf{id}}, T)$ is a $\Pi$-coloring of $\bm{H}$. In particular,
		\[
		\mathcal{P}(\bm{H}_g, t)(y) \,=\, 1.
		\] 
		It remains to observe that the $t$-ball around $x$ in $\bm{G}_f$ is isomorphic to the $t$-ball around $y$ in $\bm{H}_g$, and hence $\mathcal{P}(\bm{G}_f, t)(x) = \mathcal{P}(\bm{H}_g, t)(y) = 1$, as desired.
	\end{scproof}
	
	\subsection{Proof of Theorem~\ref{theo:dist_cont}}
	
	The proof of Theorem~\ref{theo:dist_cont} is virtually the same as the proof of Theorem~\ref{theo:dist_Borel} given in \S\ref{subsec:dist_Borel}. We just have to make sure that the construction described there produces a continuous coloring. To this end, we start with a few simple observations about topological graphs. First, it is an immediate consequence of Definition~\ref{defn:top_graph} that if $\bm{G}$ is a topological graph, then the function \[V(\bm{G}) \times \N \to \class{FSG}_\bullet \colon (x, R) \mapsto [B_{\bm{G}}(x,R), x]\] is continuous \ep{here the countable set $\class{FSG}_\bullet$ is viewed as a discrete space}. In particular, this implies that if $\mathcal{A}$ is a \LOCAL algorithm and $\bm{G}$ is a topological structured graph, then for each $T \in \N$, the map $\mathcal{A}(\bm{G}, T) \colon V(\bm{G}) \to \N$ is continuous. Next, we shall need the following fact:
	
	\begin{lemma}
		Let $\bm{G}$ be a topological structured graph. If $f \colon V(\bm{G}) \to \N$ is a continuous function, then $\bm{G}_f$ is also a topological structured graph.
	\end{lemma}
	\begin{scproof}
		Let $d$ be a metric inducing the topology on $V(\bm{G})$. Fix $x \in V(\bm{G})$, $R \in \N$, and $\epsilon > 0$. Since $\bm{G}$ is locally finite, we have $|B_{\bm{G}}(x, R)| < \infty$, and hence, by making $\epsilon$ smaller if necessary, we may arrange that $f$ is constant on the $\epsilon$-neighborhood of each vertex $z \in V(B_{\bm{G}}(x, R))$. Since $\bm{G}$ is a topological structured graph, there is $\delta > 0$ be such that for every $y$ in the $\delta$-neighborhood of $x$, the rooted structured graphs $(B_{\bm{G}}(x, R), x)$ and $(B_{\bm{G}}(y, R), y)$ are $\epsilon$-isomorphic. Let $\phi$ be an isomorphism between $(B_{\bm{G}}(x, R), x)$ and $(B_{\bm{G}}(y, R), y)$ such that $d(z, \phi(z)) < \epsilon$ for all $z \in V(B_{\bm{G}}(x, R))$. Then, by the assumption on $\epsilon$, $f(z) = f(\phi(z))$ for all $z \in V(B_{\bm{G}}(x, R))$ as well, so the rooted structured graphs $(B_{\bm{G}_f}(x, R), x)$ and $(B_{\bm{G}_f}(y, R), y)$ are $\epsilon$-isomorphic, as desired.
	\end{scproof}
	
	We also require a couple results about continuous colorings of topological graphs.
	
	\begin{lemma}\label{lemma:cont_count}
		If $G$ is a topological graph, then $G$ admits a continuous proper coloring $c \colon V(G) \to \N$.
	\end{lemma}
	\begin{scproof}
		Since $V(G)$ is a zero-dimensional Polish space, there is a countable base $(U_i)_{i=0}^\infty$ for the topology on $V(G)$ consisting of clopen sets. For each $i \in \N$, define a set $V_i \subseteq V(G)$ by setting
		\[
		x \in V_i \quad \vcentcolon\Longleftrightarrow \quad x \in U_i \text{ and } N_G(x) \cap U_i = \0.
		\]
		By construction, each set $V_i$ is independent in $G$ \ep{i.e., no two vertices in $V_i$ are adjacent}. Since $(U_i)_{i=0}^\infty$ is a base for the topology on $V(G)$ and $G$ is locally finite, for each $x \in V(G)$ there is some $i \in \N$ such that $x \in V_i$. In other words, $V(G) = \bigcup_{i=0}^\infty V_i$. Note that whether or not $x \in V_i$ is determined by the isomorphism type of the rooted radius-$1$ ball around $x$ in the structured graph $(G, 1_{U_i})$, where $1_{U_i} \colon V(G) \to \set{0,1}$ is the indicator function of $U_i$. Since $U_i$ is clopen, $1_{U_i}$ is continuous, so $(G, 1_{U_i})$ is a topological structured graph, and hence the set $V_i$ is clopen as well. To summarize, we have expressed $V(G)$ as a countable union of clopen sets that are independent in $G$. Now define $W_i \defeq V_i \setminus \bigcup_{j=0}^{i-1} V_j$. The sets $W_i$ are again clopen, independent in $G$, and satisfy $V(G) = \bigcup_{i=0}^\infty W_i$. Additionally, they are pairwise disjoint. Thus, we may define a function $c \colon V(G) \to \N$ via
		\[
		c(x) = i \quad \vcentcolon\Longleftrightarrow \quad x \in W_i.
		\]
		This $c$ is a desired continuous proper coloring.
	\end{scproof}
	
	\begin{lemma}\label{lemma:KST_cont}
		If $G$ is a topological graph of finite maximum degree $\Delta$, then $G$ admits a continuous proper coloring $c \colon V(G) \to [\Delta + 1]$.
	\end{lemma}
	\begin{scproof}
		Let $r \colon V(G) \to \N$ be a continuous proper coloring of $G$ that exists by Lemma~\ref{lemma:cont_count}. For each $i \in \N$, set $W_i \defeq r^{-1}(i)$ and define functions $c_i \colon W_i \to [\Delta + 1]$ recursively by
		\[
		c_i(x) \,\defeq\, \min \set{j \in [\Delta + 1] \,:\, \text{there is no $y \in N_G(x)$ with $r(y) < i$ and $c_{r(y)}(y) = j$}}.
		\]
		The fact that $|N_G(x)| \leq \Delta$ for all $x \in V(G)$ ensures that $c_i(x)$ is well-defined. Set $c \defeq \bigcup_{i=0}^\infty c_i$. By construction, $c$ is a proper $(\Delta+1)$-coloring of $G$. It remains to observe that $c$ is continuous, since for each $x \in W_i$, the value $c(x)$ is determined by the isomorphism type of the rooted radius-$i$ ball around $x$ in the topological structured graph $G_r$.
	\end{scproof}
	
	Now we can easily verify that the construction from \S\ref{subsec:dist_Borel} produces a continuous coloring. Let $\Pi = (t, \mathcal{P})$ be a local coloring problem and let $\class{G} \subseteq \class{FSG}$. Fix $n \in \N^+$ such that $T \defeq \Det_{\Pi, \class{G}}(n) < \infty$ and let $\mathcal{A}$ be a \LOCAL algorithm witnessing the inequality $\Det_{\Pi, \class{G}}(n) \leq T$. Set $R \defeq T + t$. Let $\bm{G}$ be a topological structured graph that is $(R,n)$-locally in $\class{G}$ and such that $|B_{\bm{G}}(x, 2R)| \leq n$ for all $x \in V(\bm{G})$. Define $G'$ to be the graph with $V(G') \defeq V(\bm{G})$ in which two distinct vertices $x$, $y$ are adjacent if and only if $\dist_{\bm{G}}(x,y) \leq 2R$. It is clear that $G'$ is a topological graph of maximum degree at most $n-1$, so, by Lemma~\ref{lemma:KST_cont}, $G'$ has a continuous proper coloring $c \colon V(\bm{G}) \to [n]$. Define a function $f \colon V(\bm{G}) \to \N$ by $f \defeq \mathcal{A}(\bm{G}_c, T)$. Since $\bm{G}_c$ is a topological structured graph, $f$ is continuous, and the argument from \S\ref{subsec:dist_Borel} shows that $f$ is a $\Pi$-coloring of $\bm{G}$, as desired.
	
	\subsection{Reduction from randomized \LOCAL algorithms to the LLL}
	
	Recall that a CSP $\mathscr{B}$ is called {bounded} if $\sup \set{|\dom(B)| \,:\, B \in \mathscr{B}} < \infty$ \ep{this is one of the assumptions in Theorem~\ref{theo:meas_SLLL}}. The main result of this section is the following lemma:
	
	\begin{lemma}\label{lemma:rand_to_LLL}
		Let $\Pi = (t, \mathcal{P})$ be a local coloring problem and let $\class{G} \subseteq \class{FSG}$. Fix $n \in \N^+$ such that $T \defeq \Rand_{\Pi, \class{G}}(n)$ is finite and set $R \defeq T + t$. Let $m \in \N^+$ and a \LOCAL algorithm $\mathcal{A}$ witness the bound $\Rand_{\Pi, \class{G}}(n) \leq T$. If $\bm{G}$ is a Borel structured graph of finite maximum degree that is $(R, n)$-locally in $\class{G}$, then there exists a bounded Borel CSP $\mathscr{B} \colon V(\bm{G}) \to^? [m]$ such that: %
		\begin{enumerate}[label=\ep{\normalfont\roman*}]
			\item\label{item:CSP_sol} for every solution $\theta \colon V(\bm{G}) \to [m]$ to $\mathscr{B}$, the function $\mathcal{A}(\bm{G}_\theta, T)$ is a $\Pi$-coloring of $\bm{G}$;
			
			\item $\pr(\mathscr{B}) \leq 1/n$ and $\de(\mathscr{B}) \leq \sup \set{|B_{\bm{G}}(x, 2R)| - 1 \,:\, x \in V(\bm{G})}$.
		\end{enumerate}
	\end{lemma}
	\begin{scproof}
		For each $\theta \colon V(\bm{G}) \to [m]$, let $f_\theta \defeq \mathcal{A}(\bm{G}_\theta, T)$. Note that for every $x \in V(\bm{G})$, the value $\mathcal{P}(\bm{G}_{f_\theta}, t)(x)$ is determined by the isomorphism type of the rooted $R$-ball around $x$ in $\bm{G}_\theta$. Therefore, we may define a set $B_x$ of functions $\phi \colon V(B_{\bm{G}}(x, R)) \to [m]$ via
		\[
		\phi \in B_x \quad \vcentcolon\Longleftrightarrow \quad \mathcal{P}(\bm{G}_{f_\theta}, t)(x) = 0 \text{ for some \ep{hence all} } \theta \colon V(\bm{G}) \to [m] \text{ such that } \theta \supseteq \phi.
		\]
		Each $B_x$ is a $(V(\bm{G}), m)$-constraint such that either $\dom(B_x) = V(B_{\bm{G}}(x, R))$ or $B_x = \0$. Let \[\mathscr{B} \,\defeq\, \set{B_x \,:\, x \in V(\bm{G})}.\] Then $\mathscr{B}$ is a Borel CSP on $V(\bm{G})$ with range $[m]$, and we claim that it has all the desired properties. Since the maximum degree of $\bm{G}$ is finite, $\mathscr{B}$ is bounded. The definition of each constraint $B_x$ implies that if $\rest{\theta}{\dom(B_x)} \not \in B_x$, then $\mathcal{P}(\bm{G}_{f_\theta}, t)(x) = 1$, and therefore \ref{item:CSP_sol} holds. The bound \[\de(\mathscr{B}) \,\leq\, \sup \set{|B_{\bm{G}}(x, 2R)| - 1 \,:\, x \in V(\bm{G})}\] follows from the observation that if $B_y \in \Nbhd(B_x)$, then $V(B_{\bm{G}}(x, R)) \cap V(B_{\bm{G}}(y, R)) \neq \0$ and thus $\dist_{\bm{G}}(x,y) \leq 2R$. \ep{Here we are subtracting $1$ because $\dist_{\bm{G}}(x, x) \leq 2R$ but $B_x \not \in \Nbhd(B_x)$.} It remains to verify that $\pr(\mathscr{B}) \leq 1/n$. To this end, fix a vertex $x \in V(\bm{G})$. We wish to show that $\P[B_x] \leq 1/n$. Since $\bm{G}$ is $(R,n)$-locally in $\class{G}$, there exist an $n$-vertex structured graph $\bm{H}$ with $[\bm{H}] \in \class{G}$ and a vertex $y \in V(\bm{H})$ such that $[B_{\bm{G}}(x, R), x] = [B_{\bm{H}}(y, R), y]$. To simplify the notation we assume, without loss of generality, that $x = y$ and $B_{\bm{G}}(x, R) = B_{\bm{H}}(x, R)$. Pick a function $\theta \colon V(\bm{H}) \to [m]$ uniformly at random and define $g \defeq \mathcal{A}(\bm{H}_\theta, T)$ and $\phi \defeq \rest{\theta}{\dom(B_x)}$
		\ep{thus, $g$ and $\phi$ are also random functions}. By the construction of $B_x$ and since $B_{\bm{G}}(x, R) = B_{\bm{H}}(x, R)$, if $\phi \in B_x$, then $\mathcal{P}(\bm{H}_g, t)(x) = 0$ and, in particular, $g$ is not a $\Pi$-coloring of $\bm{H}$. Therefore,
		\[
		\P[B_x] \,=\, \P[\phi \in B_x] \,\leq\, \P[\text{$g$ is not a $\Pi$-coloring of $\bm{H}$}] \,\leq\, \frac{1}{n},
		\]
		where the last inequality holds since $\mathcal{A}$ witnesses the bound $\Rand_{\Pi, \class{G}}(n) \leq T$.
	\end{scproof}
	
	Next we apply Lemma~\ref{lemma:rand_to_LLL} to derive Theorem~\ref{theo:dist_meas} from the \hyperref[theo:meas_SLLL]{Measurable Symmetric LLL}. For technical reasons,  a somewhat stronger version of the \hyperref[theo:meas_SLLL]{Measurable Symmetric LLL} is required:
	
	\begin{theobis}{theo:meas_SLLL}\label{theo:meas_SLLL_bis}
		Let $\mathscr{B} \colon X \to^? [m]$ be a bounded Borel CSP such that
		\[
		\pr(\mathscr{B}) \cdot (\de(\mathscr{B})+1)^8 \,\leq\, 2^{-15}.
		\]
		Then the following conclusions hold:
		\begin{enumerate}[label=\ep{\upshape\roman*}]
			\item\label{item:SLLLmeasbis} If $\mathcal{M}$ is a countable set of probability Borel measures on $X$, then $\mathscr{B}$ has a solution that is simultaneously $\mu$-measurable for all $\mu \in \mathcal{M}$.
			
			\item\label{item:SLLLBairebis} If $\mathcal{T}$ is a countable set of compatible Polish topologies on $X$, then $\mathscr{B}$ has a solution that is simultaneously $\tau$-Baire-measurable for all $\tau \in \mathcal{T}$.
		\end{enumerate}
	\end{theobis}
	
	While it appears more general, Theorem~\ref{theo:meas_SLLL_bis}\ref{item:SLLLmeasbis} is actually an easy consequence of Theorem~\ref{theo:meas_SLLL}\ref{item:SLLLmeas}. Indeed, if $\mathcal{M} = \set{\mu_n \,:\, n \in \N}$ is a countable set of probability Borel measures on $X$, then the conclusion of Theorem~\ref{theo:meas_SLLL_bis}\ref{item:SLLLmeasbis} is obtained by applying Theorem~\ref{theo:meas_SLLL}\ref{item:SLLLmeas} to the measure $\mu \defeq \sum_{n=0}^\infty 2^{-n-1} \mu_n$. On the other hand, Theorem~\ref{theo:meas_SLLL_bis}\ref{item:SLLLBairebis} requires an independent argument, which we give in \S\ref{subsec:proof_Baire}.
	
	\begin{scproof}[ of Theorem~\ref{theo:dist_meas} assuming Theorem~\ref{theo:meas_SLLL_bis}]
		Let $\Pi = (t, \mathcal{P})$ be a local coloring problem and let $\class{G} \subseteq \class{FSG}$. Fix $n \in \N^+$ such that $T \defeq \Rand_{\Pi, \class{G}}(n)$ is finite and set $R \defeq T + t$. Let $m \in \N^+$ and a \LOCAL algorithm $\mathcal{A}$ witness the bound $\Rand_{\Pi, \class{G}}(n) \leq T$. Let $\bm{G}$ be a Borel structured graph that is $(R,n)$-locally in $\class{G}$ and such that \[|B_{\bm{G}}(x, 2R)| \leq n^{1/8}/4 \quad \text{for all } x \in V(\bm{G}).\] Given a Borel probability measure $\mu$ on $V(\bm{G})$, we wish to find a $\mu$-measurable $\Pi$-coloring of $\bm{G}$ \ep{the argument for Baire-measurable colorings is the same, \emph{mutatis mutandis}}. By the Feldman--Moore theorem \cite[Theorem 1.3]{KechrisMiller}, there exist Borel involutions $\gamma_i \colon V(\bm{G}) \to V(\bm{G})$, $i \in \N$, such that $\dist_{\bm{G}}(x, y) < \infty$ if and only if $y = \gamma_i(x)$ for some $i \in \N$. Define
		\[
		\mathcal{M} \,\defeq\, \set{(\gamma_i)_\ast(\mu) \,:\, i \in \N}.
		\]
		Letting $\mathscr{B} \colon V(\bm{G}) \to^? [m]$ be a Borel CSP given by Lemma~\ref{lemma:rand_to_LLL}, we have
		\[
		\pr(\mathscr{B}) \cdot (\de(\mathscr{B}) + 1)^8 \,\leq\, \frac{1}{n} \cdot (n^{1/8}/4)^8 \,=\, 2^{-16}.
		\]
		Therefore, we may apply Theorem~\ref{theo:meas_SLLL_bis}\ref{item:SLLLmeasbis} to obtain a solution $\theta \colon V(\bm{G}) \to \N$ to $\mathscr{B}$ that is measurable with respect to every measure $(\gamma_i)_\ast(\mu)$. This means that $\theta$ agrees with a Borel function away from a set $S \subseteq V(\bm{G})$ that is $(\gamma_i)_\ast(\mu)$-null for all $i \in \N$. Let $S'$ be the \emphd{$\bm{G}$-saturation} of $S$, i.e., the set of all the vertices of $\bm{G}$ whose distance to an element of $S$ is finite. Then the set $S'$ is \emphd{$\bm{G}$-invariant}, i.e., no edges of $\bm{G}$ join $S'$ to $V(G) \setminus S'$. Notice that $S' = \bigcup_{n=0}^\infty \gamma_i (S)$, and hence
		\[
		\mu(S') \,\leq\, \sum_{n=0}^\infty \mu(\gamma_i(S)) \,=\, \sum_{n=0}^\infty ((\gamma_i)_{\ast}(\mu))(S) \,=\, 0.
		\]
		Since $\theta$ is a solution to $\mathscr{B}$, $\mathcal{A}(\bm{G}_\theta, T)$ is a $\Pi$-coloring of $\bm{G}$. Also, since $\theta$ agrees with a Borel function away from $S'$ and $S'$ is $\bm{G}$-invariant, we conclude that so does $\mathcal{A}(\bm{G}_\theta, T)$. The set $S'$ is $\mu$-null, so this implies that $\mathcal{A}(\bm{G}_\theta, T)$ is $\mu$-measurable, as desired.
	\end{scproof}
	
	We can use Lemma~\ref{lemma:rand_to_LLL} in a similar fashion to prove Theorem~\ref{theo:dist_subexp}. Instead of Theorem~\ref{theo:meas_SLLL}, this argument involves the subexponential Borel LLL of Cs\'{o}ka, Grabowski, M\'ath\'e, Pikhurko, and Tyros:
	
	\begin{theo}[{Cs\'{o}ka--Grabowski--M\'ath\'e--Pikhurko--Tyros \cite{CGMPT}}]\label{theo:subexp_LLL}
		Let $G$ be a Borel graph of subexponential growth. Fix $m \in \N^+$ and assign to each vertex $x \in V(G)$ a set $B_x$ of functions from $V(B_G(x, 1))$ to $[m]$. Viewing every $B_x$ as a $(V(G),m)$-constraint, form a CSP $\mathscr{B} \defeq \set{B_x \,:\, x \in V(G)}$. Suppose that the CSP $\mathscr{B}$ is Borel. If
		\[
		\pr(\mathscr{B}) \cdot \sup \set{|B_G(x,2)| \,:\, x \in V(G)} \,\leq\, e^{-1},
		\]
		then $\mathscr{B}$ has a Borel solution.
	\end{theo}
	
	\begin{scproof}[ of Theorem~\ref{theo:dist_subexp}]
		Let $\Pi = (t, \mathcal{P})$ be a local coloring problem and let $\class{G} \subseteq \class{FSG}$. Fix $n \in \N^+$ such that $T \defeq \Rand_{\Pi, \class{G}}(n)$ is finite and set $R \defeq T + t$. Let $m \in \N^+$ and a \LOCAL algorithm $\mathcal{A}$ witness the bound $\Rand_{\Pi, \class{G}}(n) \leq T$. Let $\bm{G}$ be a Borel structured graph of subexponential growth that is $(R,n)$-locally in $\class{G}$ and such that $|B_{\bm{G}}(x, 2R)| \leq n/e$ for all $x \in V(\bm{G})$. We wish to find a Borel $\Pi$-coloring of $\bm{G}$. To this end, let $\mathscr{B} \colon V(\bm{G}) \to^? [m]$ be a Borel CSP given by Lemma~\ref{lemma:rand_to_LLL} and let $G'$ be the graph with $V(G') \defeq V(\bm{G})$ in which two distinct vertices $x$, $y$ are adjacent if and only if $\dist_{\bm{G}}(x,y) \leq R$. It is straightforward to verify that the construction of $\mathscr{B}$ presented in the proof of Lemma~\ref{lemma:rand_to_LLL} fulfills the requirements of Theorem~\ref{theo:subexp_LLL} with the role of $G$ played by $G'$. Furthermore,
		\[
		\pr(\mathscr{B}) \cdot \sup \set{|B_{G'}(x,2)| \,:\, x \in V(G')} \,=\, \pr(\mathscr{B}) \cdot \sup \set{|B_{\bm{G}}(x,2R)| \,:\, x \in V(\bm{G})} \,\leq\, \frac{1}{n} \cdot \frac{n}{e} \,=\, e^{-1}.
		\]
		Since $G'$ is a Borel graph of subexponential growth, Theorem~\ref{theo:subexp_LLL} yields a Borel solution $\theta \colon V(\bm{G}) \to \N$ to $\mathscr{B}$. Then $\mathcal{A}(\bm{G}_\theta, T)$ is a Borel $\Pi$-coloring of $\bm{G}$, and we are done.
	\end{scproof}
	
	\section{Proof of the \hyperref[theo:meas_SLLL]{Measurable Symmetric LLL}}\label{sec:LLL_proof}
	
	\subsection{Proof outline}\label{subsec:outline}
	
	In this section we outline our strategy for proving Theorem~\ref{theo:meas_SLLL}. Say that $\mathscr{B}$ is an \emph{$(N, \epsilon)$-CSP} if
	\[
	\pr(\mathscr{B}) \cdot (\de(\mathscr{B}) + 1)^{N} \,\leq\, \epsilon.
	\]
	We shall describe here the main steps in proving that Borel $(8,2^{-15})$-CSPs have measurable solutions; the Baire-measurable case is similar but differs in some technical aspects.
	
	An important role in our argument is played by the notion of a \emph{reduction} between CSPs. Roughly speaking, a CSP $\mathscr{B}$ is reducible to a CSP $\mathscr{C}$ if there is a ``local rule'' that transforms any solution to $\mathscr{C}$ into a solution to $\mathscr{B}$; thus, to solve $\mathscr{B}$, it suffices to solve $\mathscr{C}$. The precise definition of what we mean by a ``local rule'' here is given in Definitions \ref{defn:conn} and \ref{defn:reduction}; it has some similarities and is closely related to the notion of a \LOCAL algorithm.
	The first key ingredient in the proof of Theorem~\ref{theo:meas_SLLL} is a ``bootstrapping lemma'' \ep{Lemma~\ref{lemma:boot}}, which asserts that if $\mathscr{B}$ is an $(8,2^{-15})$-CSP, then $\mathscr{B}$ is reducible to an $(N,\epsilon)$-CSP $\mathscr{C}$ for \emph{any} $N \in \N$ and $\epsilon > 0$. The proof of this fact involves the bootstrapping technique developed by Fischer and Ghaffari in \cite[\S3.2]{FG}, which in turn is based on the work of Chang and Pettie \cite{CP} on the randomized time hierarchy for the \LOCAL model. %
	A CSP $\mathscr{B}$ can be interpreted as a local coloring problem on an auxiliary graph \ep{see \S\ref{subsec:LLL_to_dist} for details}. Ghaffari, Harris, and Kuhn \cite{GHK} showed that if $\mathscr{B}$ is an $(8,2^{-15})$-CSP, then the randomized \LOCAL complexity of the corresponding local coloring problem is $o(\log n)$ \ep{the work of Ghaffari, Harris, and Kuhn builds on the earlier breakthrough of Fischer and Ghaffari \cite{FG} who obtained the same conclusion for $(32, e^{-32})$-CSPs}. Our ``bootstrapping lemma'' follows by combining the Ghaffari--Harris--Kuhn result with Lemma~\ref{lemma:rand_to_LLL}.
	
	Next we need to discuss partial solutions to CSPs. Given a CSP $\mathscr{B} \colon X \to^? [m]$, a \emph{partial solution} to $\mathscr{B}$ is a partial map $g \colon X \rightharpoonup [m]$ that can be extended to a full solution $f \colon X \to [m]$. Given a partial solution $g$ to $\mathscr{B}$, the problem of extending $g$ to a full solution can naturally be encoded as a CSP on $X \setminus \dom(g)$; we denote this CSP by $\mathscr{B}/g$ \ep{see \S\ref{subsec:part} for the definition}. The second key ingredient in our proof is Lemma~\ref{lemma:partial}, which says, roughly, that if $N$ and $1/\epsilon$ are large enough and $\mathscr{B}$ is a Borel $(N,\epsilon)$-CSP on a standard probability space $(X, \mu)$, then $\mathscr{B}$ has a Borel partial solution $g$ such that $\dom(g)$ has measure at least, say, $1/2$ and the CSP $\mathscr{B}/g$ satisfies $\pr(\mathscr{B}/g) \lessapprox \sqrt{\pr(\mathscr{B})}$. Arguments similar to Lemma~\ref{lemma:partial} have been used by Fischer and Ghaffari \cite{FG} and Molloy and Reed \cite{MR_alg} and ultimately go back to the seminal work of Beck \cite{Beck}.
	
	Now we can sketch the overall flow of the proof. We are given a Borel $(8, 2^{-15})$-CSP $\mathscr{B}$ on a standard probability space $(X,\mu)$. Using Lemma~\ref{lemma:boot}, we ``bootstrap'' $\mathscr{B}$ to an $(N, \epsilon)$-CSP, where $N$, $1/\epsilon \gg 1$. We then apply Lemma~\ref{lemma:partial} to get a Borel partial solution $g_0$ to $\mathscr{B}$ such that $\mu(\dom(g_0)) \geq 1/2$ and $\pr(\mathscr{B}/g_0)$ is still very small; in particular, by making $N$ and $1/\epsilon$ sufficiently large, we may arrange $\mathscr{B}/g_0$ to be {reducible to} an $(8, 2^{-15})$-CSP \ep{see Lemma~\ref{lemma:step}}. We then repeat the same steps with $\mathscr{B}/g_0$ in place of $\mathscr{B}$ and obtain a Borel partial solution $g_1$ to $\mathscr{B}/g_0$ such that $\mu(\dom(g_1)) \geq (1 - \mu(\dom(g_0)))/2$ and $\mathscr{B}/(g_0\cup g_1)$ is again reducible to an $(8,2^{-15})$-CSP. After countably many such iterations, we will have constructed a sequence $g_0$, $g_1$, \ldots{} of Borel functions such that $g \defeq g_0 \cup g_1 \cup \ldots$ is a partial solution to $\mathscr{B}$ with $\mu(\dom(g)) = 1$. Since $g$ is a partial solution, it can be extended to a full solution $f$, and any such $f$ is $\mu$-measurable, as desired.
	
	\subsection{Reductions between CSPs}\label{subsec:red}
	
	For sets $A$ and $B$, we use $\fun{A}{B}$ to denote the set of all partial functions $A \rightharpoonup B$. For functions $f$, $g$, the notation $g \subseteq f$ means that $f$ is an extension of $g$, i.e., $\dom(g) \subseteq \dom(f)$ and $g = \rest{f}{\dom(g)}$. The following definitions are crucial for our proof of Theorem~\ref{theo:meas_SLLL}.
	
	\begin{defn}[\textls{Connections}]\label{defn:conn}
		An \emphd{$(X,Y)$-connection}, where $X$ and $Y$ are sets, is a mapping $\rho \colon \fun{Y}{\N} \to \fun{X}{\N}$ that is monotone in the sense that if $f \colon Y \rightharpoonup \N$ and $g \subseteq f$, then $\rho(g) \subseteq \rho(f)$ as well. We say that a set $S \subseteq Y$ \emphd{$\rho$-determines} an element $x \in X$ if for all $f \colon Y \rightharpoonup \N$,
		\begin{equation}\label{eq:det}
			\rho(f)(x) \,=\, \rho(\rest{f}{S})(x).
		\end{equation}
		When one side of \eqref{eq:det} is undefined, we interpret the equality to mean that the other side is undefined as well. The \emphd{width} of an $(X,Y)$-connection $\rho$ is the quantity
		\[
		\mathsf{w}(\rho) \,\defeq\, \sup_{x \in X} \, \inf \set{|S| \,:\, \text{$S \subseteq Y$ is a set that $\rho$-determines $x$}}.
		\]
		An $(X,Y)$-connection $\rho$ is \emphd{local} if $\mathsf{w}(\rho)$ is finite. %
	\end{defn}
	
	\begin{lemdef}[\textls{$S_\rho(x)$}]
		Let $\rho$ be a local $(X,Y)$-connection. Then for each $x \in X$, there is a unique set $S_\rho(x) \subseteq Y$ such that $S \subseteq Y$ $\rho$-determines $x$ if and only if $S \supseteq S_\rho(x)$.
	\end{lemdef}
	\begin{scproof}
		This is a consequence of the fact that if $S_1$, $S_2 \subseteq Y$ are sets that $\rho$-determine $x$, then $S_1 \cap S_2$ also $\rho$-determines $x$, which holds since for any $f \colon Y \rightharpoonup \N$,
		\[
		\rho(f)(x) \,=\, \rho(\rest{f}{S_1})(x) \,=\, \rho(\rest{(\rest{f}{S_1})}{S_2})(x) \,=\, \rho(\rest{f}{(S_1 \cap S_2)})(x). \qedhere
		\]
	\end{scproof}
	
	\begin{defn}[\textls{Reductions}]\label{defn:reduction}
		Let $\mathscr{B} \colon X \to^? [m]$ and $\mathscr{C} \colon Y \to^? [n]$ be CSPs. A \emphd{reduction} from $\mathscr{B}$ to $\mathscr{C}$ is a local $(X,Y)$-connection $\rho$ such that for all $f \colon Y \to [n]$,
		\[
		\text{$f$ is a solution to $\mathscr{C}$} \quad \Longrightarrow \quad \text{$\rho(f)$ is a solution to $\mathscr{B}$}.
		\]
		The \emphd{degree} of a reduction $\rho$ is the quantity
		\[
		\de(\rho) \,\defeq\, \sup_{x \in X} \,|\set{C \in \mathscr{C} \,:\, \dom(C) \cap S_\rho(x) \neq \0}|.
		\]
	\end{defn}
	
	If there is a reduction from $\mathscr{B}$ to $\mathscr{C}$, then we say that $\mathscr{B}$ is \emphd{reducible} to $\mathscr{C}$ %
	and write $\rho \colon \mathscr{B} \rightsquigarrow \mathscr{C}$ to indicate that $\rho$ is a reduction from $\mathscr{B}$ to $\mathscr{C}$. Note that every reduction $\rho \colon \mathscr{B} \rightsquigarrow \mathscr{C}$ satisfies
	\[
	\de(\rho) \,\leq\, \mathsf{w}(\rho) \cdot (\de(\mathscr{C})+1).
	\] 
	If $X$ and $Y$ are standard Borel spaces, we say that an $(X,Y)$-connection $\rho$ is \emphd{Borel} if $\rho(f) \colon X \rightharpoonup \N$ is a Borel function whenever $f \colon Y \rightharpoonup \N$ is Borel. If $\mathscr{B}$ and $\mathscr{C}$ are Borel CSPs and there is a Borel reduction $\rho \colon \mathscr{B} \rightsquigarrow \mathscr{C}$, we say that $\mathscr{B}$ is \emphd{Borel\-/reducible} to $\mathscr{C}$. %
	
	It is clear that the \ep{Borel-}reducibility relation is reflexive: the identity map $\fun{X}{\N} \to \fun{X}{\N}$ is a reduction from any CSP $\mathscr{B}$ on $X$ to itself \ep{the width of this reduction is $1$}. It is not hard to see that this relation is also transitive. Indeed, we have the following:
	
	\begin{lemma}[\textls{Transitivity}]\label{lemma:trans}
		If $\rho \colon \mathscr{B} \rightsquigarrow \mathscr{C}$ and $\sigma \colon \mathscr{C} \rightsquigarrow \mathscr{D}$ are {Borel} reductions between {Borel} CSPs, then the composition $\rho \circ \sigma$ is a Borel reduction from $\mathscr{B}$ to $\mathscr{D}$ such that
		\begin{equation}\label{eq:comp}
			\mathsf{w}(\rho \circ \sigma) \,\leq\, \mathsf{w}(\rho)  \mathsf{w}(\sigma) \qquad \text{and} \qquad \de(\rho \circ \sigma) \,\leq\, \mathsf{w}(\rho)  \de(\sigma).
		\end{equation}
	\end{lemma}
	\begin{scproof}
		Let the given Borel CSPs $\mathscr{B}$, $\mathscr{C}$, and $\mathscr{D}$ be on spaces $X$, $Y$, and $Z$ respectively. It is clear that $\rho \circ \sigma$ is a Borel $(X,Z)$-connection that sends solutions to $\mathscr{D}$ to solutions to $\mathscr{B}$. It remains to verify inequalities \eqref{eq:comp} \ep{the first of which implies that $\rho \circ \sigma$ is local}. To this end, notice that for any $x \in X$, the union $\bigcup\set{S_\sigma(y) \,:\, y \in S_\rho(x)}$ $(\rho \circ \sigma)$-determines $x$, has size at most $\mathsf{w}(\rho) \mathsf{w}(\sigma)$, and intersects at most $\mathsf{w}(\rho)\de(\sigma)$ sets of the form $\dom(D)$ for $D \in \mathscr{D}$.
	\end{scproof}
	
	\subsection{Partial solutions}\label{subsec:part}
	
	Let $\mathscr{B} \colon X \to^? [m]$ be a CSP and let $g \colon X \rightharpoonup [m]$ be a partial function. We say that $g$ is a \emphd{partial solution} to $\mathscr{B}$ if $g$ can be extended to a solution $f \colon X \to [m]$. Given a partial map $g \colon X \rightharpoonup [m]$ and $B \in \mathscr{B}$, let $B/g$ be the constraint with domain $\dom(B/g) \defeq \dom(B) \setminus \dom(g)$ given by
	\[
	B/g \,\defeq\, \set{\phi \colon \dom(B) \setminus \dom(g) \to [m] \,:\, \phi \cup (\rest{g}{\dom(B)}) \in B}.
	\]
	In other words, $\phi \in B/g$ if and only if $\phi \cup g$ violates the constraint $B$. Note that if $\dom(B) \subseteq \dom(g)$, then $\dom(B/g) = \0$; specifically, $B/g = \set{\0}$ if $g$ violates $B$, and $B/g = \0$ otherwise. Let
	\[
	\mathscr{B}/g \,\defeq\, \set{B/g \,:\, B \in \mathscr{B}}. %
	\]
	We view $\mathscr{B}/g$ as a CSP on $X \setminus \dom(g)$. By construction, $\de(\mathscr{B}/g) \leq \de(\mathscr{B})$.
	If the CSP $\mathscr{B}$ and the partial map $g$ are Borel, then the CSP $\mathscr{B}/g$ is Borel as well. %
	The following observation is immediate:
	
	\begin{lemma}\label{lemma:part}
		Let $\mathscr{B} \colon X \to^? [m]$ be a CSP and let $g \colon X \rightharpoonup [m]$ be a partial function. Then $g$ is a partial solution to $\mathscr{B}$ if and only if $\mathscr{B}/g$ has a solution.
	\end{lemma}
	\begin{scproof}
		If $h \colon X \setminus \dom(g) \to [m]$ is a solution to $\mathscr{B}/g$, then $g \cup h$ is a solution to $\mathscr{B}$ extending $g$. Conversely, if $f$ is a solution to $\mathscr{B}$ extending $g$, then $\rest{f}{(X \setminus \dom(g))}$ is a solution to $\mathscr{B}/g$.
	\end{scproof}
	
	Next we notice that a countable union of partial solutions is also a partial solution:
	
	\begin{lemma}\label{lemma:union}
		Let $\mathscr{B} \colon X \to^? [m]$ be a CSP. If $(g_i)_{i = 0}^\infty$ is a sequence such that for all $i \in \N$, $g_i$ is a partial solution to $\mathscr{B}/(g_0 \cup \ldots \cup g_{i-1})$, then $\bigcup_{i=0}^\infty g_i$ is a partial solution to $\mathscr{B}$.
	\end{lemma}
	\begin{scproof}
		It is clear that each finite union $g_0 \cup \ldots \cup g_i$ is a partial solution to $\mathscr{B}$, so let $f_i \colon X \to [m]$ be a solution to $\mathscr{B}$ extending $g_0 \cup \ldots \cup g_i$. Consider the product space $[m]^X$, where the topology on $[m]$ is discrete. %
		By Tychonoff's theorem, $[m]^X$ is compact. The set $\mathfrak{S}$ of all solutions to $\mathscr{B}$ is closed in $[m]^X$, and, for each $i \in \N$, the set $\mathfrak{E}_i$ of all $f \colon X \to [m]$ extending $g_i$ is also closed. The functions $f_i$, $i \in \N$, certify that the collection $\set{\mathfrak{S}} \cup \set{\mathfrak{E}_i \,:\, i \in \N}$ has the finite intersection property, so %
		there is some $f \in \mathfrak{S} \cap \bigcap \set{\mathfrak{E}_i \,:\, i \in \N}$. This $f$ is a solution to $\mathscr{B}$ extending $\bigcup_{i=0}^\infty g_i$, and we are done.
	\end{scproof}
	
	Now we discuss the interplay between partial solutions and reductions:
	
	\begin{lemma}\label{lemma:pull_sol}
		Let $\rho \colon \mathscr{B} \rightsquigarrow \mathscr{C}$ be a Borel reduction between Borel CSPs. If $g$ is a Borel partial solution to $\mathscr{C}$, then $\rho(g)$ is a Borel partial solution to $\mathscr{B}$ and $\mathscr{B}/\rho(g)$ is Borel-reducible to $\mathscr{C}/g$.
	\end{lemma}
	\begin{scproof}
		Let the given Borel CSPs $\mathscr{B}$ and $\mathscr{C}$ be on spaces $X$ and $Y$ respectively. For brevity, set $Y' \defeq \dom(g)$ and $X' \defeq \dom(\rho(g))$. The fact that $\rho(g)$ is Borel follows from the Borelness of $\rho$. If $f$ is a solution to $\mathscr{C}$ extending $g$, then $\rho(f)$ is a solution to $\mathscr{B}$ and, since $\rho$ is monotone, $\rho(f)$ extends $\rho(g)$, proving that $\rho(g)$ is a partial solution to $\mathscr{B}$. Finally, the map
		\[
		\fun{Y\setminus Y'}{\N} \to \fun{X \setminus X'}{\N} \colon h \mapsto \rest{\rho(g \cup h)}{(X \setminus X')}
		\]
		is a Borel reduction from $\mathscr{B}/\rho(g)$ to $\mathscr{C}/g$.
	\end{scproof}

	\subsection{The LLL as a distributed problem}\label{subsec:LLL_to_dist}
	
	There are several natural and essentially equivalent ways of encoding a CSP as a local coloring problem on an auxiliary graph. For our purposes, the following formalism will be most convenient. A \emphd{graph-CSP} with range $m \in \N^+$ is a pair $(G, \mathscr{B})$, where $G$ is a graph and $\mathscr{B} \colon V(G) \to^? [m]$ is a bounded CSP on $V(G)$ such that $\de(\mathscr{B}) < \infty$, with the following property: If $x$, $y \in V(G)$ are distinct vertices such that $x$, $y \in \dom(B)$ for some $B \in \mathscr{B}$, then $x$ and $y$ are adjacent in $G$. A graph-CSP $(G, \mathscr{B})$ can be naturally interpreted as a structured graph, where, for each tuple of vertices $(x_1, \ldots, x_k) \in \finseq{V(G)}$, the structure contains the information about all the constraints $B \in \mathscr{B}$ such that $\dom(B) = \set{x_1, \ldots, x_k}$. Formally, we construct a structure map $\sigma$ on $G$ as follows. A tuple $\bm{x} = (x_1, \ldots, x_k) \in \finseq{V(G)}$ is in the domain of $\sigma$ if and only if:
	\begin{itemize}
		\item the vertices $x_1$, \ldots, $x_k$ are pairwise distinct; and
		\item there is a constraint $B \in \mathscr{B}$ with $\dom(B) = \set{x_1, \ldots, x_k}$.
	\end{itemize}
	To compute $\sigma(\bm{x})$ for $\bm{x} \in \dom(\sigma)$, let $\iota \colon \set{x_1, \ldots, x_k} \to [k]$ be given by $\iota(x_i) \defeq i$. For each $B \in \mathscr{B}$ with $\dom(B) = \set{x_1, \ldots, x_k}$, let $B^\ast$ be the set of all maps $\phi \colon [k] \to [m]$ such that $\phi \circ \iota \in B$, and let
	\[
	\sigma(\bm{x}) \,\defeq\, \set{B^\ast \,:\, B \in \mathscr{B} \text{ and } \dom(B) = \set{x_1, \ldots, x_k}}.
	\]
	Since $\mathscr{B}$ is bounded, $\sigma$ is defined on tuples of bounded length. Furthermore, since $\de(\mathscr{B}) < \infty$, the set $\sigma(\bm{x})$ is finite for every $\bm{x}$, which means that the range of $\sigma$ is contained in the countable set $\fins{\fins{\finf{\N}{\N}}}$. As explained after Definition~\ref{defn:str_graph}, this means that we can indeed view $\sigma$ as a structure map on $G$. We then identify $(G, \mathscr{B})$ with the structured graph $(G, \sigma, m)$. \ep{Here $m$ is interpreted as a global parameter ``known'' to every vertex; this can be realized, formally, by viewing it as a function $V(G) \to \N$ mapping each vertex $x \in V(G)$ to $m$.}
	
	Let $(G,\mathscr{B})$ be a graph-CSP. The requirement that $x$ and $y$ are adjacent whenever $x$, $y \in \dom(B)$ for some $B \in \mathscr{B}$ means that in one round of the \LOCAL model, each vertex can ``learn'' about all constraints that involve it. Because of this, the problem of solving $\mathscr{B}$ can be naturally encoded as a local coloring problem $\Pi_{\mathrm{CSP}} = (1, \mathcal{P}_{\mathrm{CSP}})$. Given parameters $m$, $k$, $p$, and $d$, let $\class{CSP}(m,k,p,d)$ denote the set of all isomorphism types of finite graph-CSPs $(G, \mathscr{B})$ such that:
	\begin{itemize}
		\item the range of $\mathscr{B}$ is $m$;
		\item $\sup \set{|\dom(B)| \,:\, B \in \mathscr{B}} \leq k$;
		\item $\pr(\mathscr{B}) \leq p$ and $\de(\mathscr{B}) \leq d$.
	\end{itemize}
	We shall need the following bound on the randomized \LOCAL complexity of solving graph-CSPs:
	
	\begin{theo}[{Ghaffari--Harris--Kuhn \cite{GHK}}]\label{theo:dist_LLL}
		If $m$, $k$, $p$, and $d$ satisfy $p(d+1)^8 \leq 2^{-15}$, then
		\[
		\Rand_{\Pi_{\mathrm{CSP}}, \class{CSP}(m,k,p,d)}(n) \,=\, \exp(O(\sqrt{\log \log n})) \,=\, o(\log n).
		\]
		\ep{Here the implicit constants in the asymptotic notation may depend on $m$, $k$, and $d$.}
	\end{theo}
	
	We say that a bounded CSP $\mathscr{B}$ is an \emphd{$(N, \epsilon)$-CSP} if $\pr(\mathscr{B}) (\de(\mathscr{B}) + 1)^{N} \leq \epsilon$. A bounded Borel CSP $\mathscr{B}$ is a \emphd{potential Borel $(N,\epsilon)$-CSP} if it is Borel-reducible to a Borel $(N,\epsilon)$-CSP. Using Theorem~\ref{theo:dist_LLL}, we derive the following ``bootstrapping lemma,'' inspired by \cite[\S3.2]{FG}:
	
	\begin{lemma}[\textls{Bootstrapping}]\label{lemma:boot}
		Let $\mathscr{B}$ be a potential Borel $(8,2^{-15})$-CSP. Then $\mathscr{B}$ is in fact a potential Borel $(N,\epsilon)$-CSP for any $N \in \N$ and $\epsilon > 0$. Moreover, there exist a Borel $(N,\epsilon)$-CSP $\mathscr{C}$ and a Borel reduction $\rho \colon \mathscr{B} \rightsquigarrow \mathscr{C}$ such that $\pr(\mathscr{C}) \de(\rho)^N \leq \epsilon$.
	\end{lemma}
	\begin{scproof}
		Let $\mathscr{B} \colon X \to^? [m]$ be the given Borel CSP and let $\sigma \colon \mathscr{B} \rightsquigarrow \mathscr{D}$ be a Borel reduction from $\mathscr{B}$ to a Borel $(8, 2^{-15})$-CSP $\mathscr{D} \colon Y \to^? [\ell]$. Set $k \defeq \sup \set{|\dom(D)| \,:\, D \in \mathscr{D}}$, $p \defeq \pr(\mathscr{D})$, and $d \defeq \de(\mathscr{D})$. Define a graph $G$ with $V(G) \defeq Y$ by making distinct vertices $x$ and $y$ adjacent if and only if $x$, $y \in \dom(D)$ for some $D \in \mathscr{D}$. Then $\bm{G} \defeq (G, \mathscr{D})$ is a graph-CSP. It is routine to check that, as a structured graph, $\bm{G}$ is Borel. Note that $\Delta \defeq \Delta(\bm{G}) = \Delta(G) \leq (k-1)(d + 1)$; in particular, $\Delta$ is finite. Consider the local coloring problem $\Pi_{\mathrm{CSP}}$. For $n \in \N^+$, let
		\[
		T(n) \,\defeq\, \Rand_{\Pi_{\mathrm{CSP}}, \class{CSP}(\ell,k,p,d)}(n) \qquad \text{and} \qquad R(n) \,\defeq\, T(n) + 1.
		\]
		By Theorem~\ref{theo:dist_LLL}, $T(n) = o(\log n)$, and hence $R(n) = o(\log n)$. Let $m_n \in \N^+$ and a \LOCAL algorithm $\mathcal{A}_n$ witness the bound $\Rand_{\Pi_{\mathrm{CSP}}, \class{CSP}(\ell,k,p,d)}(n) \leq T(n)$. For all large enough $n$, $\bm{G}$ is $(R(n),n)$-locally in $\class{CSP}(\ell,k,p,d)$. Indeed, the class $\class{CSP}(\ell,k,p,d)$ is closed under %
		adding isolated vertices, and every finite induced subgraph of $\bm{G}$ is in $\class{CSP}(\ell,k,p,d)$. Thus, we only need to verify that, for large $n$, $|B_{\bm{G}}(x, R(n))| \leq n$ for all $x \in Y$, which holds as
		\[
		|B_{\bm{G}}(x, R(n))| \,\leq\, 1 + \Delta^{o(\log n)} \,=\, n^{o(1)}.
		\]
		By Lemma~\ref{lemma:rand_to_LLL}, for every large enough $n$, there is a bounded Borel CSP $\mathscr{C}_n \colon Y \to^? [m_n]$ such that:
		\begin{enumerate}[label=\ep{\normalfont\roman*}]
			\item\label{item:sol_to_sol} for every solution $\theta \colon Y \to [m_n]$ to $\mathscr{C}_n$, the function $\mathcal{A}_n(\bm{G}_\theta, T(n))$ is a $\Pi_{\mathrm{CSP}}$-coloring of $\bm{G}$, or, equivalently, a solution to $\mathscr{D}$;
			
			\item $\pr(\mathscr{C}_n) \leq 1/n$ and $\de(\mathscr{C}_n) \leq \sup \set{|B_{\bm{G}}(x, 2R(n))| - 1 \,:\, x \in Y} \leq 1 + \Delta^{o(\log n)} = n^{o(1)}$.
		\end{enumerate}
		Using \ref{item:sol_to_sol}, we can define a Borel reduction $\tau_n \colon \mathscr{D} \rightsquigarrow \mathscr{C}_n$ by setting, for all $\theta \colon Y \rightharpoonup \N$ and $x \in Y$,
		\[
		\tau(\theta)(x) \,\defeq\, \begin{cases}
			\mathcal{A}_n(\bm{G}_\theta, T(n))(x) &\text{if $V(B_{\bm{G}}(x, R(n))) \subseteq \dom(\theta)$};\\
			\text{undefined} &\text{otherwise}.
		\end{cases}
		\]
		Each $x \in Y$ is $\tau_n$-determined by the vertex set of the radius-$R(n)$ ball around $x$ in $\bm{G}$. Thus,
		\[
		\mathsf{w}(\tau_n) \,\leq\, \sup\set{|B_{\bm{G}}(x, R(n))| \,:\, x \in Y} \,\leq\, 1 + \Delta^{o(\log n)} \,=\, n^{o(1)},
		\]
		and hence $\de(\tau_n) \leq \mathsf{w}(\tau_n)(\de(\mathscr{C}_n)+1) \leq n^{o(1)}$ as well. By Lemma~\ref{lemma:trans}, $\rho_n \defeq \sigma \circ \tau_n$ is a Borel reduction from $\mathscr{B}$ to $\mathscr{C}_n$ such that $\mathsf{d}(\rho_n) \leq \mathsf{w}(\sigma) \de(\tau_n) \leq n^{o(1)}$ \ep{since $\mathsf{w}(\sigma)$ is a constant independent of $n$}. Now fix $N \in \N$ and $\epsilon > 0$. For sufficiently large $n$, we can write
		\[
		\pr(\mathscr{C}_n) \cdot (\de(\mathscr{C}_n) + 1)^N \,\leq\, \frac{1}{n} \cdot (n^{o(1)})^N \,=\, n^{-1 + o(1)} \,<\,\epsilon,
		\]
		and similarly $\pr(\mathscr{C}_n) \de(\rho_n)^N < \epsilon$. Therefore, for all sufficiently large $n$, the CSP $\mathscr{C} \defeq \mathscr{C}_n$ and the reduction $\rho \defeq \rho_n$ are as desired.
	\end{scproof}
	
	\subsection{Proof of Theorem~\ref{theo:meas_SLLL}\ref{item:SLLLmeas}}\label{subsec:proof_meas}
	
	In order to build partial solutions to CSPs, we use the following lemma, which is an adaptation of \cite[Lemma 8]{FG} to the measurable setting:
	
	\begin{lemma}[\textls{Partial solutions}]\label{lemma:partial}
		Let $\mathscr{B} \colon X \to^? [m]$ be a Borel CSP and let $\rho \colon \mathscr{B} \rightsquigarrow \mathscr{C}$ be a Borel reduction from $\mathscr{B}$ to a Borel CSP $\mathscr{C} \colon Y \to^? [n]$. Assume that $\mathscr{C}$ is a $(2, e^{-2}/n^2)$-CSP. If $\mu$ is a probability Borel measure on $X$, then there is a Borel partial solution $h \colon Y \rightharpoonup [n]$ to $\mathscr{C}$ such that
		\[
		\pr(\mathscr{C}/h) \,\leq\, n\sqrt{\pr(\mathscr{C})} \qquad \text{and} \qquad \mu(\dom(\rho(h))) \,\geq\, 1 - \de(\rho) \sqrt{\pr(\mathscr{C})}.
		\]
	\end{lemma}
	
	\begin{scproof}\stepcounter{ForClaims} \renewcommand{\theForClaims}{\ref{lemma:partial}}
		Set $p\defeq \pr(\mathscr{C})$ and $d \defeq \de(\mathscr{C})$. For a set $A$ and $i \in \N$, let $\mathsf{const}(A,i)$ denote the constant function on $A$ sending every $x \in A$ to $i$. An important observation is that if $C$ is a $(Y,n)$-constraint and $A \subseteq Y$ is a subset such that $|A \cap \dom(C)| \leq 1$, then
		\begin{equation}\label{eq:two}
			\P[C] \,=\, \frac{1}{n}\sum_{i = 1}^n \P[C/\mathsf{const}(A,i)].
		\end{equation}
		Indeed, if $A \cap \dom(C) = \0$, then $C/\mathsf{const}(A,i) = C$ for all $i \in [n]$ and \eqref{eq:two} is trivial. Otherwise, there is a single element $x \in A \cap \dom(C)$, and $\P[C/\mathsf{const}(A,i)]$ is the probability that a uniformly random function $\phi \colon \dom(C) \to [n]$ is in $C$ conditioned on the event $\phi(x) = i$. Since each of these events has probability $1/n$, equation \eqref{eq:two} follows.
		
		We say that a subset $A \subseteq X$ is \emphd{$\mathscr{C}$-discrete} if $|A \cap \dom(C)| \leq 1$ for all $C \in \mathscr{C}$. The first step in the construction of a desired partial solution $h$ is to fix a partition $Y = A_1 \sqcup \ldots \sqcup A_N$ of $Y$ into finitely many Borel $\mathscr{C}$-discrete sets. To see that such a partition exists, let $G$ be the graph with vertex set $Y$ in which distinct vertices $x$ and $y$ adjacent if and only if $x$, $y \in \dom(C)$ for some $C \in \mathscr{C}$. Since $\mathscr{C}$ is bounded and $d < \infty$, $G$ has finite maximum degree, so, by Theorem~\ref{theo:KST}, %
		$G$ admits a Borel proper coloring $c \colon Y \to [N]$ for some $N \in \N$. Letting $A_k \defeq c^{-1}(k)$, the partition $Y = A_1 \sqcup \ldots \sqcup A_N$ is as desired. Next, we need the following definition:
		
		\begin{smalldefn}[\textls{Dangerous constraints and elements}]
			Given a partial map $h \colon Y \rightharpoonup [n]$, we say that a constraint $C \in \mathscr{C}$ is \emphd{$h$-dangerous} if $\P[C/h] > \sqrt{p}$, and an element $x \in Y$ is \emphd{$h$-dangerous} if $x \in \dom(C)$ for some $h$-dangerous constraint $C \in \mathscr{C}$. Let $D(h) \subseteq Y$ denote the set of all $h$-dangerous elements; thus, $D(h) = \bigcup \set{\dom(C) \,:\, C \in \mathscr{C} \text{ is $h$-dangerous}}$.
		\end{smalldefn}
		
		We use $\concat$ to denote concatenation of finite sequences and write $u \sqsubseteq w$ to mean that a sequence $u$ is an initial segment of a sequence $w$. For $\ell \in \N$, let $W_\ell$ denote the set of all sequences of elements of $[n]$ of length at most $\ell$. %
		To each sequence $w \in W_N$ of length $k$, we associate a Borel partial map $h_{w} \colon A_1 \sqcup \ldots \sqcup A_k \rightharpoonup [n]$  as follows. If $k = 0$ \ep{i.e, $w = \0$}, then we let $h_\0 \defeq \0$. If $k > 0$, we can express $w$ as $w = u \concat i$, where $i \in [n]$ is the last entry of $w$ and $u$ is the initial segment of $w$ of length $k-1$. Then we recursively define
		\[
		h_{u \concat i} \,\defeq\, h_{u} \sqcup \mathsf{const}(A_k\setminus D(h_{u}), i).
		\]
		In other words, if $w = (w_1, w_2,  \ldots, w_k)$, then $h_w$ assigns the value $w_t$ to the points in $A_t$, except that it does not assign any value to the elements of $A_t$ that are dangerous with respect to the partial map $h_{(w_1, \ldots,  w_{t-1})}$. Our goal now is to argue that that there exists a sequence $w \in [n]^N$ such that $h_{w}$ satisfies the conclusion of Lemma~\ref{lemma:partial}.
		
		For brevity, given $C \in \mathscr{C}$ and $w \in W_N$, we write $C_w \defeq C/h_w$ and $\mathscr{C}_w \defeq \mathscr{C}/h_w$. We also say that a constraint or an element is \emphd{$w$-dangerous} to mean that it is $h_w$-dangerous, and define $D(w) \defeq D(h_w)$. If $u \in W_{N - 1}$ is a sequence of length $k - 1$, $i \in [n]$, and $C \in \mathscr{C}$, then, by construction,
		\[
		C_{u \concat i} \,=\, C_u/\mathsf{const}(A_k\setminus D(u), i).
		\]
		The set $A_k\setminus D(u)$ is $\mathscr{C}$-discrete, so we may apply \eqref{eq:two} to conclude that
		\begin{equation}\label{eq:concat}
			\P[C_u] \,=\, \frac{1}{n}\sum_{i = 1}^n \P[C_{u \concat i}].
		\end{equation}
		
		\begin{claim}\label{claim:p_bound}
			For all $w \in W_N$, $\pr(\mathscr{C}_w) \leq n\sqrt{p}$ and $\de(\mathscr{C}_w) \leq d$.
		\end{claim}
		\begin{claimproof}
			The bound $\de(\mathscr{C}_w) \leq d$ is clear. To prove $\pr(\mathscr{C}_w) \leq n\sqrt{p}$, we proceed by induction on the length of $w$. If $w = \0$, then $\mathscr{C}_\0 = \mathscr{C}$, so $\pr(\mathscr{C}_\0) = p \leq n\sqrt{p}$. Now suppose that $w = u \concat i$ for some $i \in [n]$. Consider an arbitrary constraint $C \in \mathscr{C}$. If $\P[C_u] > \sqrt{p}$, i.e., $C$ is $u$-dangerous, then $C_{u \concat i} = C_u$, so $\P[C_{u \concat i}] = \P[C_u] \leq n\sqrt{p}$ by the inductive hypothesis. If, on the other hand, $\P[C_u] \leq \sqrt{p}$, then \[\P[C_{u \concat i}] \,\leq\, \sum_{j = 1}^n \P[C_{u \concat j}] \,\stackrel{\text{\eqref{eq:concat}}}{=}\, n \P[C_u] \,\leq\, n \sqrt{p}.\] In either case, $\P[C_{u \concat i}] \leq n \sqrt{p}$, and thus $\pr(\mathscr{C}_{u \concat i}) \leq n\sqrt{p}$, as desired.
		\end{claimproof}
		
		\begin{claim}\label{claim:ps}
			If $w \in W_N$, then $h_w$ is a partial solution to $\mathscr{C}$.
		\end{claim}
		\begin{claimproof}
			Using Claim~\ref{claim:p_bound} and the fact that $\mathscr{C}$ is a $(2, e^{-2}/n^2)$-CSP, we can write
			\[
			\pr(\mathscr{C}_w) \cdot (\de(\mathscr{C}_w)+1) \,\leq\, n\sqrt{p}\cdot  (d+1) \,\leq\, e^{-1}.
			\]
			Hence, $\mathscr{C}_w$ has a solution by the LLL \ep{see Theorem~\ref{theo:SLLL}}. By Lemma~\ref{lemma:part}, this precisely means that $h_w$ is a partial solution to $\mathscr{C}$.
		\end{claimproof}
		
		In view of Claims~\ref{claim:p_bound} and \ref{claim:ps}, it only remains to argue that there is some $w \in [n]^N$ such that $\mu(\dom(\rho(h_w))) \geq 1 - \de(\rho) \sqrt{p}$.
		
		\begin{claim}\label{claim:monotone}
			The following statements are valid.
			
			\begin{enumerate}[label=\ep{\normalfont\roman*}]
				\item\label{item:mon1} Let $u$, $w \in W_N$ be such that $u \sqsubseteq w$. Then $D(u) \subseteq D(w)$.
				
				\item\label{item:mon2} If $w \in [n]^N$, then $\dom(h_w) \supseteq Y \setminus D(w)$.
				
				\item\label{item:mon3} If $w \in [n]^N$ and $x \in X$ are such that $S_\rho(x) \cap D(w) = \0$, then $x \in \dom(\rho(h_w))$.
			\end{enumerate}
		\end{claim}
		\begin{claimproof}
			\ref{item:mon1} It is enough to prove the claim when $w = u \concat i$ for some $i \in [n]$. Suppose $x \in D(u)$, i.e., $x \in \dom(C)$ for a $u$-dangerous constraint $C \in \mathscr{C}$. Then $C_{u \concat i} = C_u$, so $\P[C_{u \concat i}] = \P[C_u] > \sqrt{p}$ and $C$ is $(u \concat i)$-dangerous as well, which implies that $x \in D(u \concat i)$.
			
			\ref{item:mon2} We need to show that for each $k \in [N]$, $A_k \setminus \dom(h_w) \subseteq D(w)$. To this end, let $u \sqsubseteq w$ be the initial segment of $w$ of length $k-1$. Then, by construction, $A_k \setminus \dom(h_w) = D(u) \subseteq D(w)$.
			
			\ref{item:mon3} Let $f$ be a solution to $\mathscr{C}$ extending $h_w$ \ep{such $f$ exists by Claim~\ref{claim:ps}}. Then $\rho(f)$ is a solution to $\mathscr{B}$ and, in particular, $\rho(f)(x)$ is defined. By \ref{item:mon2}, $\dom(h_w) \supseteq S_\rho(x)$, so $f$ agrees with $h_w$ on $S_\rho(x)$. Therefore, $\rho(h_w)(x)$ is defined and equal to $\rho(f)(x)$.
		\end{claimproof}
		
		We are now ready to finish the proof. By applying \eqref{eq:concat} $N$ times, we see that for all $C \in \mathscr{C}$,
		\begin{equation}\label{eq:average}
			\P[C] \,=\, \frac{1}{n^N} \sum_{w \in [n]^N} \P[C_w].
		\end{equation}
		Pick $w \in [n]^N$ uniformly at random. We use $\P_w$ and $\mathbb{E}_w$ to denote probability and expectation with respect to this random choice of $w$. With this notation, \eqref{eq:average} can be rewritten as $\mathbb{E}_w[\P[C_w]] = \P[C]$. By Markov's inequality applied to the random variable $\P[C_w]$, this yields
		\[
		\P_w[\text{$C$ is $w$-dangerous}] \,=\, \P_w[\P[C_w] > \sqrt{p}] \,\leq\, \frac{\mathbb{E}_w[\P[C_w]]}{\sqrt{p}} \,=\, \frac{\P[C]}{\sqrt{p}} \,\leq\, \sqrt{p}.
		\]
		Consider any $x \in X$. By Claim~\ref{claim:monotone}\ref{item:mon3},
		\begin{align}
			\P_w[x \in \dom(\rho(h_w))] \,&\geq\, 1 - \P_w[S_\rho(x) \cap D(w) \neq \0] \nonumber\\
			&\geq\, 1 - \sum_C \P_w[\text{$C$ is $w$-dangerous}] \,\geq\, 1 - \de(\rho)\sqrt{p}, \label{eq:fiber}
		\end{align}
		where the sum is over all $C \in \mathscr{C}$ such that $\dom(C) \cap S_\rho(x) \neq \0$. Since \eqref{eq:fiber} holds for every $x \in X$, we conclude that, by Fubini's theorem, there is a choice of $w \in [n]^N$ such that $\mu(\dom(\rho(h_w))) \geq 1 - \de(\rho) \sqrt{p}$, and the proof is complete
	\end{scproof}
	
	The upper bound on $\pr(\mathscr{C}/g)$ given by Lemma~\ref{lemma:partial} depends on $n$, the cardinality of the range of $\mathscr{C}$. In order to control it, we show that every bounded CSP is reducible to a \emphd{binary} CSP, i.e., one with range $[2]$:
	
	\begin{lemma}\label{lemma:binary}
		Let $\rho \colon \mathscr{B} \rightsquigarrow \mathscr{C}$ be a Borel reduction between bounded Borel CSPs and let $\epsilon > 0$. Then there exist a bounded binary {Borel} CSP $\mathscr{D}$ and a Borel reduction $\sigma \colon \mathscr{B} \rightsquigarrow \mathscr{D}$ such that
		\[
		\pr(\mathscr{D}) \,\leq\, (1+\epsilon)\pr(\mathscr{C}), \qquad \de(\mathscr{D}) \,=\, \de(\mathscr{C}), \qquad \text{and} \qquad \de(\sigma) \,=\, \de(\rho).
		\]
	\end{lemma}
	\begin{scproof}
		Let the given Borel CSPs be $\mathscr{B} \colon X \to^? [m]$ and $\mathscr{C} \colon Y \to^? [n]$. We need to somehow reduce $n$ to just $2$. The idea is to represent each $i \in [n]$ by a binary sequence of an appropriate finite length $N$. We can then replace every element $y \in Y$ by $N$ copies $(y,1)$, \ldots, $(y, N)$ and identify a solution $f \colon Y \to [n]$ to $\mathscr{C}$ with a function $f' \colon Y \times [N] \to [2]$ where $f'(y, i)$ is the $i$-th digit in a binary representation of $f(y)$. This will allow us to reduce $\mathscr{C}$ to a CSP $\mathscr{D} \colon Y \times [N] \to^? [2]$. We remark that slight technical complications arise from the fact that $n$ may not be a power of 2 (this is why we lose a factor of $1 + \epsilon$ in the bound on $\pr(\mathscr{D})$).
		
		Let us now describe the construction formally. Since the CSP $\mathscr{C}$ is bounded by assumption, the value $k \defeq \sup \set{|\dom(C)| \,:\, C \in \mathscr{C}}$ is finite. Fix $\delta > 0$ such that $(1+\delta)^k \leq 1 + \epsilon$. Let $N \in \N$ be so large that it is possible to express $2^N$ as a sum of the form $2^N = s_1 + s_2 + \cdots + s_n$, where each $s_i$ is a positive integer and $\max \set{s_i \,:\, i \in [n]} \leq (1+\delta)2^N/n$ (so $2^N$ is ``approximately divisible by $n$''). Fix an arbitrary mapping $\xi \colon [2]^N \to [n]$ such that $|\xi^{-1}(i)| = s_i$ for each $i \in [n]$. We think of each tuple $(c_1, \ldots, c_N) \in [2]^N$ as a ``binary code'' representing the value $\xi(c_1, \ldots, c_N) \in [n]$.
		
		Let $Z \defeq Y \times [N]$ and for each $f \colon Z \rightharpoonup \N$, define $\tau(f) \colon Y \rightharpoonup \N$ as follows:
		\[
		\tau(f)(y) \,\defeq\, \begin{cases}
			\xi(f(y,1), \ldots, f(y, N)) &\text{if $(y,i) \in \dom(f)$ and $f(y,i) \in [2]$ for all $i \in [N]$};\\
			\text{undefined} &\text{otherwise}.
		\end{cases}
		\]
		It is clear that $\tau$ is a Borel $(Y, Z)$-connection such that $\mathsf{w}(\tau) = N$, with the set $\set{(y,1), \ldots, (y,N)}$ $\tau$-determining each $y \in Y$. Next we define a binary Borel CSP $\mathscr{D}$ on $Z$ so that $\tau$ is a reduction from $\mathscr{C}$ to $\mathscr{D}$, in the obvious way. Namely, for each $C \in \mathscr{C}$, we let $C^\ast$ be the $(Z, 2)$-constraint with $\dom(C^\ast) \defeq \dom(C) \times [N]$ given by
		$
		C^\ast \defeq \set{\phi \colon \dom(C) \times [N] \to [2] \,:\, \tau(\phi) \in C}
		$,
		and set $\mathscr{D} \defeq \set{C^\ast \,:\, C \in \mathscr{C}}$. Let $\sigma \defeq \rho \circ \tau$. By Lemma~\ref{lemma:trans}, $\sigma$ is a Borel reduction from $\mathscr{B}$ to $\mathscr{D}$. We claim that the CSP $\mathscr{D}$ and the reduction $\sigma$ are as desired.

		The equality $\de(\mathscr{D}) = \de(\mathscr{C})$ follows since for all $C \in \mathscr{C}$, $|\Nbhd(C^\ast)| = |\Nbhd(C)|$. To see that $\de(\sigma) = \de(\rho)$, observe that for each $x \in X$, the set $S_\rho(x) \times [N]$ $\sigma$-determines $x$ and $\dom(C^\ast) \cap (S_\rho(x) \times [N]) \neq \0$ if and only if $\dom(C) \cap S_\rho(x) \neq \0$. Finally, to bound $\pr(\mathscr{D})$, let $C \in \mathscr{C}$ and consider an arbitrary function $\psi \in C$. If $\phi \colon \dom(C^\ast) \to [2]$ satisfies $\tau(\phi) = \psi$, then $(\phi(y, 1), \ldots, \phi(y, N)) \in \xi^{-1}(\psi(y))$ for all $y \in \dom(C)$. Hence there are at most $(1+\delta)2^N/n$ possible values for the tuple $(\phi(y, 1), \ldots, \phi(y, N))$, and so the number of such functions $\phi$ cannot exceed $((1+\delta)2^N/n)^{|\dom(C)|}$. Therefore,
		\begin{align*}
			\P[C^\ast] \,=\, \frac{|C^\ast|}{2^{|\dom(C^\ast)|}} \,&=\, \frac{\sum_{\psi \in C} |\set{\phi \,:\, \tau(\phi) = \psi}|}{2^{N|\dom(C)|}} \,\leq\, \frac{1}{2^{N|\dom(C)|}} \cdot |C| \cdot \left(\frac{(1+\delta)2^N}{n}\right)^{|\dom(C)|} \\
			&=\, (1+\delta)^{|\dom(C)|} \cdot \frac{|C|}{n^{|\dom(C)|}} \,=\, (1+\delta)^{|\dom(C)|} \P[C] \,\leq\, (1+\epsilon)\P[C],
		\end{align*}
		and thus $\pr(\mathscr{D}) \leq (1+\epsilon) \pr(\mathscr{C})$, as desired.
	\end{scproof}
	
	Now we can combine Lemmas~\ref{lemma:boot}, \ref{lemma:partial}, and \ref{lemma:binary} to obtain the following:
	
	\begin{lemma}\label{lemma:step}
		Let $\mathscr{B} \colon X \to^? [m]$ be a potential Borel $(8,2^{-15})$-CSP. If $\mu$ is a probability Borel measure on $X$, then $\mathscr{B}$ admits a Borel partial solution $g \colon X \rightharpoonup [m]$ such that $\mu(\dom(g)) \geq 1/2$ and $\mathscr{B}/g$ is again a potential Borel $(8,2^{-15})$-CSP.
	\end{lemma}
	\begin{scproof}
		Using Lemmas~\ref{lemma:boot} and \ref{lemma:binary}, we can find a binary Borel $(16, 2^{-32})$-CSP $\mathscr{C}$ and a Borel reduction $\rho \colon \mathscr{B} \rightsquigarrow \mathscr{C}$ such that $\pr(\mathscr{C})\de(\rho)^2 \leq 1/4$. Let $h$ be a partial Borel solution to $\mathscr{C}$ satisfying the conclusion of Lemma~\ref{lemma:partial} and define $g \defeq \rho(h)$. By Lemma~\ref{lemma:pull_sol}, $g$ is a Borel partial solution to $\mathscr{B}$ and $\mathscr{B}/g$ is Borel-reducible to $\mathscr{C}/h$. Thus, it suffices to verify that $\mu(\dom(g)) \geq 1/2$ and $\mathscr{C}/h$ is an $(8,2^{-15})$-CSP. By the choice of $h$, we have
		\[
		\mu(\dom(g)) \,\geq\, 1- \de(\rho)\sqrt{\pr(\mathscr{C})} \,\geq\, 1 - 1/2 \,=\, 1/2,
		\]
		as desired. Furthermore,
		\[
		\pr(\mathscr{C}/h) \cdot (\de(\mathscr{C}/h) + 1)^8 \,\leq\, 2\sqrt{\pr(\mathscr{C})} \cdot (\de(\mathscr{C}) + 1)^8 \,\leq\, 2 \cdot \sqrt{2^{-32}} \,=\, 2^{-15}.
		\]
		Thus, $\mathscr{C}/h$ is indeed an $(8,2^{-15})$-CSP, and we are done.
	\end{scproof}
	
	With Lemma~\ref{lemma:step} in hand, it is easy to finish the proof of Theorem~\ref{theo:meas_SLLL}\ref{item:SLLLmeas}. Let $\mathscr{B} \colon X \to^? [m]$ be a Borel $(8, 2^{-15})$-CSP and let $\mu$ be a probability Borel measure on $X$. Repeated applications of Lemma~\ref{lemma:step} produce a sequence $(g_i)_{i=0}^\infty$ such that for all $i \in \N$:
	\begin{itemize}
		\item $\mathscr{B}/(g_0 \cup \ldots \cup g_{i-1})$ is a potential Borel $(8,2^{-15})$-CSP;
		
		\item $g_i$ is a Borel partial solution to $\mathscr{B}/(g_0 \cup \ldots \cup g_{i-1})$; and
		
		\item $\mu(\dom(g_{i})) \geq (1 - \mu(\dom(g_0)) - \cdots - \mu(\dom(g_{i-1})))/2$.
	\end{itemize}
	Let $g \defeq \bigcup_{i=0}^\infty g_i$. Then $\mu(\dom(g)) = 1$. By Lemma~\ref{lemma:union}, $g$ is a Borel partial solution to $\mathscr{B}$, so let $f$ be any solution to $\mathscr{B}$ extending $g$. Since $f$ agrees with the Borel function $g$ on a $\mu$-conull set, $f$ is $\mu$-measurable, and the proof is complete.

	\subsection{Proof of Theorem~\ref{theo:meas_SLLL_bis}\ref{item:SLLLBairebis}}\label{subsec:proof_Baire}
	
	In the Baire category setting, instead of Lemma~\ref{lemma:step} we use the following:
	
	\begin{lemma}\label{lemma:Baire_step}
		Let $\mathscr{B} \colon X \to^? [m]$ be a potential Borel $(8,2^{-15})$-CSP. Then there is a finite set $\mathcal{G}$ of Borel partial solutions to $\mathscr{B}$ such that:
		\begin{itemize}
			\item $\bigcup \set{\dom(g) \,:\, g \in \mathcal{G}} = X$; and
			\item for every $g \in \mathcal{G}$, $\mathscr{B}/g$ is a potential Borel $(8,2^{-15})$-CSP.
		\end{itemize}
	\end{lemma}
	\begin{scproof}
		Using Lemmas~\ref{lemma:boot} and \ref{lemma:binary}, we can find a binary Borel $(16, 2^{-32})$-CSP $\mathscr{C}$ and a Borel reduction $\rho \colon \mathscr{B} \rightsquigarrow \mathscr{C}$ such that $\pr(\mathscr{C})\de(\rho)^2 \leq 1/4$. Using the construction from the proof of Lemma~\ref{lemma:partial}, we obtain a natural number $N \in \N$ and an assignment to each sequence $w \in [2]^N$ of a Borel partial solution $h_w$ to $\mathscr{C}$ with the following properties:
		\begin{enumerate}[label=\ep{\normalfont\roman*}]
			\item\label{item:1} for all $w \in [2]^N$, $\pr(\mathscr{C}/h_w) \leq 2 \sqrt{\pr(\mathscr{C})}$ \ep{Claim~\ref{claim:p_bound}}; and
			
			\item\label{item:2} for every $x \in X$, $\P_w[x \in \dom(\rho(h_w))] \geq 1 - \de(\rho) \sqrt{\pr(\mathscr{C})}$, where $\P_w$ denotes probability with respect to a uniformly random choice of $w \in [2]^N$ \ep{equation \eqref{eq:fiber}}.
		\end{enumerate}
		From \ref{item:1} it follows that each $\mathscr{C}/h_w$ is an $(8,2^{-15})$-CSP, since
		\[
		\pr(\mathscr{C}/h_w) \cdot (\de(\mathscr{C}/h_w) + 1)^8 \,\leq\, 2 \sqrt{\pr(\mathscr{C})} \cdot (\de(\mathscr{C}) + 1)^8 \,\leq\, 2 \cdot \sqrt{2^{-32}} \,=\, 2^{-15}.
		\]
		Furthermore, from \ref{item:2} it follows that for all $x \in X$,
		\[
		\P_w[x \in \dom(\rho(h_w))] \,\geq\, 1 - \de(\rho) \sqrt{\pr(\mathscr{C})} \,\geq\, 1 - 1/2 \,=\, 1/2.
		\]
		In particular, for each $x \in X$, there is some $w \in [2]^N$ such that $x \in \dom(\rho(h_w))$. By Lemma~\ref{lemma:pull_sol}, each $\rho(h_w)$ is a Borel partial solution to $\mathscr{B}$ and $\mathscr{B}/\rho(h_w)$ is Borel-reducible to $\mathscr{C}/h_w$. Therefore, the set $\mathcal{G} \defeq \set{\rho(h_w) \,:\, w \in [2]^N}$ is as desired.
	\end{scproof}
	
	We are now ready to prove Theorem~\ref{theo:meas_SLLL_bis}\ref{item:SLLLBairebis}. Let $\mathscr{B} \colon X \to^? [m]$ be a Borel $(8, 2^{-15})$-CSP and let $\mathcal{T}$ be a countable set of compatible Polish topologies on $X$. %
	For each $\tau \in \mathcal{T}$, fix a countable base $\mathcal{U}_\tau$ consisting of nonempty open sets, and let $(\tau_i, U_i)_{i=0}^\infty$ be an enumeration of all pairs $(\tau, U)$ with $\tau \in \mathcal{T}$ and $U \in \mathcal{U}_\tau$. Using Lemma~\ref{lemma:Baire_step}, we can recursively build a sequence $(g_i)_{i=0}^\infty$ such that for all $i \in \N$:
	\begin{itemize}
		\item $\mathscr{B}/(g_0 \cup \ldots \cup g_{i-1})$ is a potential Borel $(8,2^{-15})$-CSP;
		
		\item $g_i$ is a Borel partial solution to $\mathscr{B}/(g_0 \cup \ldots \cup g_{i-1})$; and
		
		\item the set $(\dom(g_0) \cup \ldots \cup \dom(g_i)) \cap U_i$ is $\tau_i$-nonmeager. %
	\end{itemize}
	Specifically, once $g_0$, \ldots, $g_{i-1}$ have been constructed, we define $g_i$ as follows. Let $\mathcal{G}$ be a finite set of Borel partial solutions to $\mathscr{B}/(g_0 \cup \ldots \cup g_{i-1})$ given by Lemma~\ref{lemma:Baire_step}. Then
	\[
	U_i \,\subseteq \, \dom(g_0) \cup \ldots \cup \dom(g_{i-1}) \cup \bigcup \set{\dom(g) \,:\, g \in \mathcal{G}}.
	\]
	Since $U_i$ is nonempty and $\tau_i$-open, it is $\tau_i$-nonmeager, so $\dom(g) \cap U_i$ must be $\tau_i$-nonmeager for some $g \in \set{g_0, \ldots, g_{i-1}} \cup \mathcal{G}$. If $\dom(g_j) \cap U_i$ is $\tau_i$-nonmeager for some $0 \leq j \leq i -1$, then
	we can let $g_i \defeq \0$, and otherwise
	we can make $g_i$ be any $g \in \mathcal{G}$ such that $\dom(g) \cap U_i$ is $\tau_i$-nonmeager. %
	
	Once we have built a sequence $(g_i)_{i=0}^\infty$ as above, we let $g \defeq \bigcup_{i=0}^{\infty} g_i$. By Lemma~\ref{lemma:union}, $g$ is a Borel partial solution to $\mathscr{B}$. Furthermore, for each $\tau \in \mathcal{T}$, the set $\dom(g) \cap U$ is $\tau$-nonmeager for every $U \in \mathcal{U}_\tau$, which, by \cite[Proposition 8.26]{KechrisDST}, implies that $\dom(g)$ is $\tau$-comeager. %
	Now let $f$ be any solution to $\mathscr{B}$ extending $g$. For each $\tau \in \mathcal{T}$, $f$ agrees with the Borel function $g$ on a $\tau$-comeager set, and therefore $f$ is $\tau$-Baire-measurable, which finishes the proof.

	\printbibliography%
	
\end{document}